\documentclass[reqno]{amsart}

\usepackage[utf8]{inputenc}
\usepackage[english]{babel}
\usepackage{amsmath, amssymb, amsfonts}
\usepackage{mathrsfs}
\usepackage{mathtools}
\usepackage{bm}
\usepackage{booktabs}
\usepackage{graphicx}
\usepackage{caption}
\usepackage{subfigure}
\usepackage[table]{xcolor}
\usepackage{multirow}
\usepackage{hyperref}
\usepackage{pifont}
\usepackage{tikz}
\usepackage[foot]{amsaddr}

\usetikzlibrary{positioning, shapes.geometric, arrows.meta, calc, decorations.pathreplacing, patterns, shadows, fit, backgrounds}
\usepackage[top=1in, bottom=1.25in, left=1.25in, right=1.25in]{geometry}
\usepackage{algorithm}
\usepackage{algorithmic}
\usepackage{doi}
\usepackage[numbers,sort&compress]{natbib}

\usepackage[autostyle, english = american]{csquotes}
\MakeOuterQuote{"}

\usepackage[colorinlistoftodos, textsize=tiny]{todonotes}
\newcommand{\eg}{{\it e.g.}}

\def\u{{\bm u}}

\def\x{{\bm x}}
\def\y{{\bm y}}
\def\r{{\bm r}}

\newcommand{\uHF}{\u_{\mathrm{HF}}}
\newcommand{\uLF}{\u_{\mathrm{LF}}}

\def\Pcal{\mathcal{P}}
\def\Mcal{\mathcal{M}}

\def\R{\mathbb{R}}

\newtheorem{rem}{Remark}[section]

\hypersetup{
    colorlinks=true,
    linkcolor=blue,
    filecolor=magenta,
    urlcolor=cyan,
    citecolor=red,
    pdftitle={A Multi-Fidelity and Parametric Framework for Reduced-Order Modeling using Optimal Transport-based Interpolation},
    pdfauthor={Moaad Khamlich, et al.},
}

\begin{document}

\title[Multi-Fidelity and Parametric OT-ROMs]{A Multi-Fidelity and Parametric Reduced-Order Modeling Framework with Optimal Transport-based Interpolation: Applications to Diffused-Interface Two-Phase Flows}

\author{Moaad Khamlich$^{1}$}
\author{Niccolò Tonicello$^{1}$}
\author{Federico Pichi$^{1}$}
\author{Gianluigi Rozza$^{1}$}
\address{$^1$ mathLab, Mathematics Area, SISSA, via Bonomea 265, I-34136 Trieste, Italy}

\begin{abstract}
This work introduces a data-driven, non-intrusive reduced-order modeling (ROM) framework that leverages Optimal Transport (OT) for multi-fidelity and parametric problems in two-phase flows modelling. Building upon the success of displacement interpolation for data augmentation in handling nonlinear dynamics, we extend its application to more complex and practical scenarios. The framework is designed to correct a computationally inexpensive low-fidelity (LF) model to match an accurate high-fidelity (HF) one by capturing its temporal evolution via displacement interpolation while preserving the problem's physical consistency.
The framework is further extended to address systems dependent on a physical parameter, for which we construct a surrogate model using a hierarchical, two-level interpolation strategy. First, it creates synthetic HF checkpoints via displacement interpolation in the parameter space. Second, the residual between these synthetic HF checkpoints and a true LF solution is interpolated in the time domain using the multi-fidelity OT-based methodology. This strategy provides a robust and efficient way to explore the parameter space and to obtain a refined description of the dynamical system.
The potential of the method is discussed in the context of complex and computationally expensive diffuse-interface methods for two-phase flow simulations, which are characterized by moving interfaces and nonlinear evolution, and challenging to be dealt with traditional ROM techniques.

\end{abstract}

\maketitle

\section{Introduction}
\label{sec:intro}

Many engineering and physical phenomena of practical interest can be accurately described by mathematical models based on partial differential equations (PDEs). The numerical solution of such PDEs is generally performed through spatio-temporal discretization, using methods like the Finite Element (FE) or Finite Volume (FV) methods, which transform the continuous problem into a system with a finite number $N_h$ of degrees of freedom~\cite{QuarteroniNumericalApproximationPartial1994}. An accurate discretization of general PDEs often yields high-dimensional solution vectors, meaning that a proper description of their behavior requires a computationally intractable number of degrees of freedom ($N_h \gg 1$). For time-dependent problems, further discretization in time results in a sequence of large algebraic systems to be solved, compounding the computational burden.

Although these traditional high-fidelity (HF) models can provide robust and accurate predictions, their intrinsic computational cost often represents a bottleneck. This complexity has a substantial impact on their feasibility to analyze physical phenomena, especially when the HF strategy is part of a more general framework where its output serves as input for other tasks. In many practical scenarios, such as real-time control~\cite{ManzoniOptimalControlPartial2021,PichiDrivingBifurcatingParametrized2022}, uncertainty quantification~\cite{BabuskaStochasticCollocationMethod2007,ChenReducedBasisMethods2017}, design optimization~\cite{GiacominiSurrogateModelTopology2026,BoulleControlBifurcationStructures2022}, and inverse problems~\cite{isakov2017inverse,MulaInverseProblemsDeterministic2023}, the computational demands of HF methods become prohibitive. The need to perform faster simulations, without excessively compromising their accuracy, has driven extensive research into the field of model order reduction~\cite{Benner2015, Quarteroni2016, Hesthaven2016}.

Reduced-order models (ROMs) have emerged as a surrogate strategy exploiting data to alleviate the computational burden posed by accurate but costly HF methods, providing rapid and reliable evaluations. The key concept behind ROMs is that many practical problems exhibit dominant patterns that can be described by a low-dimensional set of features. This allows to discard redundant variables in the description of the phenomena, motivating the use of dimensionality reduction techniques to compress the dataset constituted by physical simulation while retaining its essential properties. Typically, the construction of a ROM involves two main stages: an \emph{offline stage}, where the computationally intensive work of generating high-fidelity solution snapshots and extracting dominant features occurs, and an \emph{online stage}, where the pre-computed reduced representation is used to rapidly approximate the system's solution for new parameter values or time instances. This offline-online decomposition makes ROMs particularly useful in the many-query and real-time settings.

A popular class of ROMs relies on projecting the governing equations onto a low-dimensional linear subspace, often constructed using Proper Orthogonal Decomposition (POD)~\cite{SirovichTurbulenceDynamicsCoherent1987a,Chatterjee2000}. While effective for many diffusion-dominated problems~\cite{Lassila2013,BuffaPrioriConvergenceGreedy2012,RomorROMViscousIncompressible2025}, these linear-subspace methods encounter significant challenges when dealing with systems that exhibit strong nonlinearities or advection-dominated phenomena~\cite{HesthavenNonlinearModelReduction2026,GlasReducedBasisMethod2020,ArbesKolmogorovNwidthLinear2025}. The mathematical root of this difficulty lies in the slow decay of the Kolmogorov $n$-width, which measures the degree of compactness of the solution manifold in terms of approximability by $n$-dimensional linear subspaces~\cite{Cohen15,OhlbergerReducedBasisMethods2016,PeherstorferBreakingKolmogorovBarrier2022}. In particular, advection-dominated problems featuring steep gradients or traveling waves that move across the domain are very difficult to represent efficiently by linear methods with a fixed global basis. The moving nature of such features often leads to poor support intersection between snapshots at different time instances, fundamentally limiting the effectiveness of linear approximation and often resulting in Gibbs phenomena or requiring an impractically large reduced basis to maintain accuracy.

These fundamental limitations have spurred extensive research into more advanced ROM techniques. One direction involves adapting the data representation to better suit linear methods through shifted or registered POD approaches~\cite{Taddei2020,NoninoCalibrationBasedALEModel2024a,Blickhan23}. Another direction concerns the development of nonlinear reduction methodologies using deep learning frameworks~\cite{LeeModelReductionDynamical2020, fresca2021,PichiGraphConvolutionalAutoencoder2024}, which seek to learn nonlinear mappings between the high-dimensional state space and a low-dimensional latent space.

In our previous work~\cite{Khamlich2025JCPDI}, we introduced a ROM framework that addresses these challenges by incorporating the concept of Displacement Interpolation (DI) from Optimal Transport (OT) theory~\cite{villani2021topics, Peyre2019,SantambrogioOptimalTransportApplied2015} to enhance the representation of nonlinear dynamics in complex systems. OT, originally developed by Monge~\cite{Monge1781} and Kantorovich~\cite{Kantorovich1942}, concerns the matching and transport of probability distributions by identifying the most cost-effective way to transform one into another. Unlike classical $L^2$-based euclidean distance metrics that compare distributions pointwise, OT considers the ``work'' required to move a source mass distribution to match a target distribution, inherently accounting for the underlying geometry of the space.

The fundamental concept behind this is the so-called Wasserstein distance, which quantifies the minimum cost of transporting one distribution to another~\cite{EhrlacherNonlinearModelReduction2020,BattistiWassersteinModelReduction2023,DoSparseWassersteinBarycenters2025}. This distance defines a metric in the space of probability measures, endowing it with a Riemannian-like structure that allows for geometric operations such as the definition of geodesics and barycenters~\cite{agueh2011barycenters, mccan-interpolation}. The geometric perspective provided by OT is particularly attractive for ROMs concerned with the treatment of features that translate or deform, as the Wasserstein distance can capture these transformations in a more natural way compared to traditional Euclidean-based metrics. By interpreting solution snapshots as probability distributions and leveraging displacement interpolation---the construction of geodesic paths in Wasserstein space---our method generates physically plausible synthetic snapshots between sparse high-fidelity checkpoints, addressing the challenge of data scarcity while enabling continuous-time prediction~\cite{IolloMappingCoherentStructures2022,CucchiaraModelOrderReduction2024,Khamlich2025JCPDI}.

In this paper, we extend this OT-based methodology to address two challenges in computational modeling that arise frequently in practice, namely the \emph{multi-fidelity} and the \emph{parametric} modeling. 
Indeed, in the first context, the goal is to correct computationally inexpensive low-fidelity (LF) approximations with limited accuracy, for instance due to a coarser mesh or simplified physics, using information from a few expensive high-fidelity (HF) simulations~\cite{Peherstorfer2018,HowardMultifidelityDeepOperator2023,ContiMultifidelityReducedorderSurrogate2023}. Here, we propose a new approach where this information is retrieved by interpolating in time the residual field between HF and LF solutions via OT, leveraging the geometric structure of Wasserstein space to track the movement of discrepancy features.

Moreover, many physical systems depend on a set of parameters $\bm{\mu}$, such as fluids and material properties or boundary/initial conditions, and the goal is to build a ROM capable of rapid prediction for any $\bm{\mu}$ within a given parameter domain. To investigate efficiently different system's configurations, we introduce a hierarchical, two-level interpolation strategy to simultaneously deal with both parameter and time dependency. In the first level, we use the OT-based displacement interpolation across the parameter space to generate synthetic HF snapshots for unseen parameter value from HF checkpoints computed at nearby training points. In the second level, we compute the residual between these physically consistent HF synthetic data and true LF simulations at the test parameter, and then apply our multi-fidelity framework to interpolate this residual in time. This approach avoids the construction of a monolithic ROM over the joint time-parameter space, which can be prohibitively complex.
These two extensions, the Multi-Fidelity OT-ROM (MF-OT-ROM) and the Parametric Multi-Fidelity OT-ROM (PMF-OT-ROM), broaden the applicability of our original method.

The remainder of the paper is structured as follows. Section \ref{sec:intro_2p} presents a brief review of diffuse-interface methods for two-phase flows, with a particular focus on surrogate modelling, serving as the application context for the proposed methodologies. Section \ref{sec:preliminaries} provides the necessary mathematical background on reduced-order modeling and optimal transport theory. Section \ref{sec:ot_rom} summarizes the core OT-ROM methodology from our previous work~\cite{Khamlich2025JCPDI}, which serves as the foundation for the extensions presented here. Section \ref{sec:mf_ot_rom} introduces multi-fidelity residual interpolation, while Section \ref{sec:pmf_ot_rom} presents parametric displacement interpolation, applicable to both single- and multi-fidelity settings. Section \ref{sec:problem_def} introduces the conservative Allen-Cahn equation coupled within the five-equation model, while Section \ref{sec:numerical_results} is reserved for numerical experiments. Finally, Section \ref{sec:conclusions} provides concluding remarks and directions for future research.
\section{Diffuse-interface two-phase modelling}\label{sec:intro_2p}
We frame our discussion within the context of two-phase flow modelling using the conservative Allen-Cahn equation~\cite{chiu2011conservative} coupled with the five-equation model~\cite{allaire2002five}. More broadly, diffuse-interface methods have gained significant attention in the computational fluid dynamics community due to their versatility across different numerical schemes, mesh anisotropies, governing equations and flow regimes (from incompressible to compressible, from laminar to turbulent), as well as their relative ease of implementation (see~\cite{saurel2018diffuse} for an extensive review). In these approaches, the interface between the two fluids is represented as a thin but finite transition layer described by a continuous phase-field variable, thereby avoiding the need for explicit interface tracking or reconstruction as it custom in Volume-of-Fluid (VOF) or Level-Set (LS) methods. This modelling paradigm offers significant advantages in terms of numerical robustness and flexibility when dealing with complex interfacial topological changes such as breakup, coalescence, or strong deformation.

The widespread adoption of diffuse-interface methods has enabled the simulation of highly complex two-phase flows, which are often characterized by strongly coupled multi-physics phenomena. These include interfacial mass and heat transfer~\cite{mirjalili2022computational}, phase change~\cite{haghani2021phase,huang2022consistent,salimi2025low}, boiling processes~\cite{roccon2025boiling,ibrahim2025conservative,weber2026consistent}, turbulent interfacial dynamics~\cite{soligo2019breakage,roccon2023phase}, and strong compressibility effects~\cite{huang2023consistent,collis2025thermodynamically,white2025high}. The ability to address such tightly coupled physical processes within a unified modelling framework makes diffuse-interface methods particularly attractive for the simulation of realistic engineering and scientific applications. At the same time, the increasing fidelity and resolution of these simulations often come with a substantial computational cost, motivating the development of efficient surrogate and reduced-order modelling strategies.

In particular, recent developments of the conservative Allen-Cahn equation and its variants have led to numerical methods exhibiting several desirable properties, including improved mass conservation~\cite{soligo2019mass}, strict boundedness~\cite{mirjalili2020conservative,pirozzoli2025efficient,huang2025bound}, and effective interface sharpening~\cite{al2021high,jain2022accurate,tonicello2024high}. These advancements have contributed to the growing popularity of formulations based on the conservative Allen-Cahn equations for high-fidelity two-phase flow simulations. Despite the growing body of literature on high-fidelity simulations based on these approaches, the integration of reduced-order modelling techniques in this context has only recently begun to receive attention~\cite{wiewel2019latent,kani2019reduced,haas2020bubcnn,cutforth2025convolutional}. Developing robust and accurate ROMs for such systems remains particularly challenging due to the strongly nonlinear and transport-dominated nature of the underlying dynamics, such as the case of two-phase flows being characterized by complex interfacial dynamics, sharp gradients, and multi-scale interactions. In this setting, we aim at showing the potentiality of OT-ROM techniques to efficiently capture transport-dominated phenomena.

Indeed, in the context of diffuse-interface models, several properties of the phase-field variable make OT particularly well suited for this class of problems. The phase-field variable distinguishing the two phases is in fact naturally bounded between zero and one. Moreover, in incompressible or low-Mach-number compressible flows, the total volume associated with the phase-field is preserved over time. These properties resemble those of probability distributions within the OT formalism, where quantities are non-negative and globally conserved. This analogy provides a natural connection between diffuse-interface representations of two-phase flows and OT-based reduction strategies, further confirming that OT-based techniques may offer a promising pathway for developing robust  and efficient reduced-order models in this setting.

Finally, in this work we extend the methodological developments of OT-based ROMs and investigate their applicability to two-phase flow simulations. By doing so, we aim to bridge the gap between high-fidelity diffuse-interface models and computationally efficient ROMs. The resulting framework leverages OT to construct reduced representations that remain accurate even in the presence of strongly transport-dominated dynamics and evolving interfaces.

\section{OT-ROMs Preliminaries}
\label{sec:preliminaries}

We consider time-dependent PDE problems defined on a spatial domain $\Omega \subset \R^d$ over a time interval $\mathcal{I}=[0, t_f]$. After spatial discretization using, for instance, FE or FV methods, the state of the system at time $t$ is represented by a vector $\u(t) \in \R^{N_h}$, where $N_h$ denotes the number of degrees of freedom. Given a fixed parametric configuration, the collection of solutions over time traces out a solution manifold $\Mcal = \{\u(t) \mid t \in \mathcal{I}\}$ embedded in the high-dimensional space $\R^{N_h}$. In practice, we work with a discrete set of $N_T$ snapshots collected at times $\{t_k\}_{k=1}^{N_T}$, which we arrange into the snapshot matrix $\mathbf{S} = [\u(t_1) | \u(t_2) | \dots | \u(t_{N_T})] \in \R^{N_h \times N_T}$.

\subsection{Reduced-Order Modeling: The Offline-Online Paradigm}

ROM emerged as class of surrogate strategies to alleviate computational bottlenecks posed by HF methods~\cite{Quarteroni2016,Hesthaven2016,Benner2015}.

In particular, dimensionality reduction techniques aim to identify the underlying low-dimensional structure of the dominant pattern of the solution manifold characterized by the evolution in time and the variation w.r.t.\ input parameters. Toward this goal, the system dynamics is projected onto a reduced set of coordinates, building a surrogate model that is much faster to solve than its HF counterpart, and retains the essential physical characteristics and accuracy.

The construction of a ROM typically involves two main stages. In the \emph{offline stage}, which is computationally intensive and can exploit high-performance computing resources, one generates a set of high-fidelity solutions (snapshots) for several values of the parameter or time instances. These snapshots are processed to extract the main modes that describe the system's behavior, and a low-dimensional basis is constructed to represent the system's dynamics in the reduced space. In the \emph{online stage}, the pre-computed reduced representation is used to rapidly approximate the system's solution for new parameter values or time instances. Since the computations occur in a significantly smaller dimensional space ($N_r \ll N_h$), this stage entails extremely fast evaluation.

The so-called projection-based ROMs approximate the high-dimensional solution $\u(t)$ as the linear combination of a few basis functions:
\begin{equation}
    \u(t) \approx \mathbf{V} \mathbf{a}_r(t),
\end{equation}
where $\mathbf{V} \in \R^{N_h \times N_r}$ is a matrix whose $N_r \ll N_h$ orthonormal columns span a low-dimensional subspace, and $\mathbf{a}_r(t) \in \R^{N_r}$ contains the reduced coordinates. A common choice for $\mathbf{V}$ is the POD basis, obtained by computing the singular value decomposition (SVD) of the snapshot matrix:
\begin{equation}
    \mathbf{S} = \mathbf{U} \bm{\Sigma} \mathbf{W}^T,
\end{equation}
where $\mathbf{U} \in \R^{N_h \times N_T}$ contains the left singular vectors (POD modes), $\bm{\Sigma}$ is a diagonal matrix of singular values $\sigma_1 \geq \sigma_2 \geq \cdots \geq \sigma_{N_T} \geq 0$, and $\mathbf{W}$ contains the right singular vectors. The reduced basis $\mathbf{V}$ is formed by retaining only the first $N_r$ columns of $\mathbf{U}$, corresponding to the largest singular values. The truncation is typically chosen to capture a prescribed fraction of the total energy:
\begin{equation}
    \frac{\sum_{i=1}^{N_r} \sigma_i^2}{\sum_{i=1}^{N_T} \sigma_i^2} \geq 1 - \epsilon,
\end{equation}
where $\epsilon$ is a small tolerance. The rate at which the singular values decay determines how many modes are needed for an accurate approximation, and thus the reliability of the ROM ansatz. As described by the Kolmogorov n-width decay analysis~\cite{Cohen15,OhlbergerReducedBasisMethods2016,PeherstorferBreakingKolmogorovBarrier2022}, for advection-dominated problems with moving features this decay is typically slow, necessitating a large number of modes and thus compromising the computational advantage of the ROM.

\subsection{Optimal Transport and Displacement Interpolation}

\begin{figure}[htbp]
    \centering
    \begin{tikzpicture}[
        source/.style={circle, draw=blue!80, fill=blue!20, thick, minimum size=10pt, inner sep=0pt},
        target/.style={regular polygon, regular polygon sides=4, draw=orange!80, fill=orange!20, thick, minimum size=9pt, inner sep=0pt},
        plan/.style={-Stealth, draw=gray!70, thick, opacity=0.8}
    ]
        \begin{scope}[xshift=-4cm]
            \node[source] (s1) at (0,1.5) {};
            \node[] at (0, 2.05) {Source};
            \node[source] (s2) at (0,0.5) {};
            \node[source] (s3) at (0,-0.5) {};
            \node[source] (s4) at (0,-1.5) {};
            \node[anchor=east] at (-0.5, 0) {$\mu$};
        \end{scope}
        
        \begin{scope}[xshift=-2cm]
            \node[target] (t1) at (0,1.5) {};
            \node[] at (0, 2) {Target};
            \node[target] (t2) at (0,0.4) {};
            \node[target] (t3) at (0,-0.6) {};
            \node[target] (t4) at (0,-1.5) {};
            \node[anchor=west] at (0.5, 0) {$\nu$};
        \end{scope}

        \draw[plan, line width=1pt] (s1) to (t1);
        \draw[plan, line width=1pt] (s2) to (t3);
        \draw[plan, line width=1pt] (s3) to (t2);
        \draw[plan, line width=1pt] (s4) to (t4);
        \node at (-3,-2) {$\bm{\pi}$};
        \node at (-3,-2.8) {Monge problem};

        \begin{scope}[xshift=2cm]
            \node[source] (s1) at (0,1) {};
            \node[] at (0, 2.05) {Source};
            \node[source] (s2) at (0,-1) {};
            \node[anchor=east] at (-0.5, 0) {$\mu$};
        \end{scope}

        \begin{scope}[xshift=4cm]
            \node[target] (t1) at (0,1.5) {};
            \node[] at (0, 2) {Target};
            \node[target] (t2) at (0,0) {};
            \node[target] (t3) at (0,-1.5) {};
            \node[anchor=west] at (0.5, 0) {$\nu$};
        \end{scope}

        \draw[plan, line width=1.5pt] (s1) to (t1);
        \draw[plan, line width=0.8pt] (s1) to (t2);
        \draw[plan, line width=0.8pt] (s2) to (t2);
        \draw[plan, line width=1.5pt] (s2) to (t3);
        \node at (3,-2) {$\bm{\pi}$};
        \node at (3,-2.8) {Kantorovich problem};

    \end{tikzpicture}
    \caption{Illustration of the different OT problem formulations, depicting the mass transport plan $\pi_{ij}$ proportional to the thickness of the arrows.}
    \label{fig:ot_problems}
\end{figure}
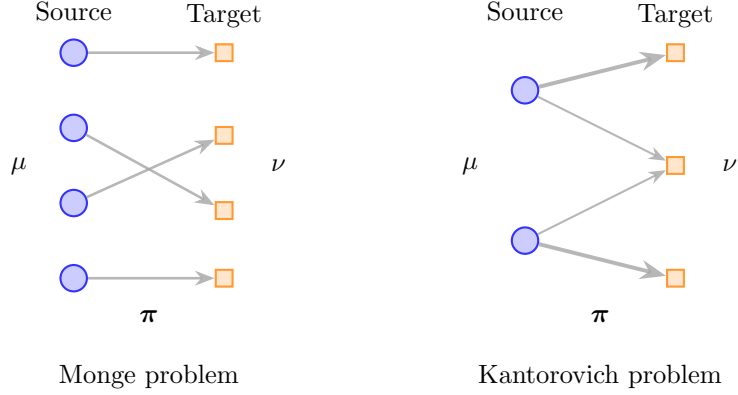

Optimal Transport (OT) theory provides a geometric framework for comparing and morphing probability distributions~\cite{villani2021topics, Peyre2019}. At its core, OT addresses the problem of finding the most efficient way to transform one distribution of mass into another, where efficiency is measured by a total transportation cost. The theory was originally posed by Monge~\cite{Monge1781}, looking for the cheapest way to transport soil to fill a hole, and later relaxed by Kantorovich~\cite{Kantorovich1942} allowing mass from a single source point to be split and distributed among multiple destinations. This latter interpretation led to the modern theory of finding transport maps, for a given cost, which are optimal in matching source $\mu$ and target $\nu$ distributions. Figure~\ref{fig:ot_problems} illustrates the difference between the two formulations in terms of mass splitting and the resulting plans $\bm{\pi}$ to transport the two distributions.

Despite more general formulations are available~\cite{BenamouComputationalFluidMechanics2000,LevyNumericalAlgorithmL22015,KhamlichEfficientNumericalStrategies2026}, here we specialize our discussion directly to the discrete setting, which is of practical importance for numerical applications. Let the source and target measures be discrete probability measures with finite supports, written as weighted sums of Dirac masses $\delta$:
\begin{equation}
    \mu = \sum_{i=1}^n a_i \delta_{\x_i}, \quad \text{with } \boldsymbol{a} \in \mathbb{R}_+^{n}, \quad \sum_{i=1}^n a_i = 1,
\end{equation}
\begin{equation}
    \nu = \sum_{j=1}^m b_j \delta_{\y_j}, \quad \text{with } \boldsymbol{b} \in \mathbb{R}_+^{m}, \quad \sum_{j=1}^m b_j = 1,
\end{equation}
where $\boldsymbol{a}=(a_1, \dots, a_n)^T$ and $\boldsymbol{b}=(b_1, \dots, b_m)^T$ are the weight vectors (or histograms) of the measures. A \emph{coupling} or \emph{transport plan} between $\mu$ and $\nu$ is represented by a matrix $\bm{\pi} \in \R_+^{n \times m}$, where $\pi_{ij}$ specifies how much mass is transported from $\x_i$ to $\y_j$. The set of admissible couplings is the convex polytope of matrices with the prescribed marginals:
\begin{equation}
    \Pi(\mu, \nu) = \left\{ \bm{\pi} \in \mathbb{R}_+^{n \times m} \mid \bm{\pi} \mathbf{1}_m = \boldsymbol{a}, \ \bm{\pi}^T \mathbf{1}_n = \boldsymbol{b} \right\},
\end{equation}
where $\mathbf{1}_k$ denotes a vector of ones of size $k$. The cost function $c(\x,\y)$, representing the cost of moving a unit of mass from $\x$ to $\y$, is encoded in a cost matrix $\mathbf{C} \in \mathbb{R}^{n \times m}$ with $C_{ij} = c(\x_i, \y_j)$. 

The Kantorovich problem of finding the optimal (cheapest) transport plan then becomes the following linear program:
\begin{equation}
\label{eq:ot_kantorovich}
    C(\mu, \nu) \doteq \min_{\bm{\pi} \in \Pi(\mu, \nu)} \sum_{i=1}^n \sum_{j=1}^m C_{ij} \pi_{ij} = \min_{\bm{\pi} \in \Pi(\mu, \nu)} \langle \mathbf{C}, \bm{\pi} \rangle_F,
\end{equation}
where $\langle \cdot, \cdot \rangle_F$ denotes the Frobenius inner product.

When the cost function is the $p$-th power of the Euclidean distance, the $p$-th root of the optimal cost defines the $p$-Wasserstein distance:
\begin{equation}
    W_p(\mu, \nu) \doteq \left( \min_{\bm{\pi} \in \Pi(\mu, \nu)} \sum_{i,j} \|\x_i - \y_j\|^p \, \pi_{ij} \right)^{1/p}.
\end{equation}
The space of probability measures equipped with this metric inherits a rich geometric structure. Unlike classical $L^2$-based metrics, which only measure pointwise differences, the Wasserstein distance is sensitive to the spatial displacement of features. Indeed, two distributions representing a feature and its translated counterpart are close in Wasserstein distance, reflecting the physical intuition that small translations should correspond to small distances, but very far in the standard Euclidean metric. Then, it is clear that the choice of the underlying metric is particularly relevant in the context of physical fields exhibiting sharp contrasts, discontinuous features and small compact supports.

A fundamental property of spaces endowed with the Wasserstein metric is that they admit geodesics, i.e., shortest paths between measures. Given an optimal transport plan $\bm{\pi}^*$ between $\mu_0$ and $\mu_1$, one can construct a geodesic path $\{\mu_t\}_{t \in [0,1]}$ called the \emph{displacement interpolation} or McCann interpolation~\cite{mccan-interpolation,IolloMappingCoherentStructures2022,Khamlich2025JCPDI}:
\begin{equation}
\label{eq:displacement_interp}
    \mu_t = \sum_{i=1}^{n} \sum_{j=1}^{m} \pi^*_{ij} \, \delta_{(1-t)\x_i + t\y_j},
\end{equation}
where each mass $\pi^*_{ij}$ travels along a straight line from source point $\x_i$ to target point $\y_j$, with its position at time $t$ given by the convex combination $(1-t)\x_i + t\y_j$. This path is a constant-speed geodesic satisfying:
\begin{equation}
    W_p(\mu_s, \mu_t) = |t - s| \, W_p(\mu_0, \mu_1), \quad \forall\, s, t \in [0,1].
\end{equation}
The displacement interpolation provides a geometrically natural way to ``morph'' between distributions, and given its nice properties in dealing with moving features, we leverage this as the key tool for generating synthetic snapshots describing the time- and parametric-evolution of physical fields. 

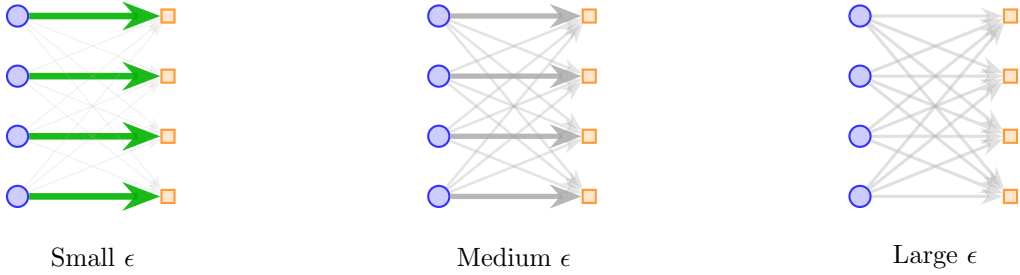
\begin{figure}[t]
    \centering
    \begin{tikzpicture}[scale=0.8,
        source/.style={circle, draw=blue!80, fill=blue!20, thick, minimum size=8pt, inner sep=0pt},
        target/.style={regular polygon, regular polygon sides=4, draw=orange!80, fill=orange!20, thick, minimum size=7pt, inner sep=0pt},
        plan/.style={-Stealth, thick, opacity=0.8}
    ]
    \begin{scope}[xshift=-7cm]
        \foreach \i in {1,...,4}{
          \node[source] (s\i) at (0, \i-0.5) {};
          \node[target] (t\i) at (2.5, \i-0.5) {};
        }
        \draw[plan, green!70!black, line width=2.5pt, opacity=0.9] (s1) to (t1);
        \draw[plan, green!70!black, line width=2.5pt, opacity=0.9] (s2) to (t2);
        \draw[plan, green!70!black, line width=2.5pt, opacity=0.9] (s3) to (t3);
        \draw[plan, green!70!black, line width=2.5pt, opacity=0.9] (s4) to (t4);
        \foreach \i in {1,...,4}{
            \foreach \j in {1,...,4}{
                \ifnum\i=\j\else
                    \draw[plan, gray!50, line width=0.5pt, opacity=0.2] (s\i) to (t\j);
                \fi
            }
        }
        \node at (1.25,-0.5) {Small $\epsilon$};
    \end{scope}

    \begin{scope}
        \foreach \i in {1,...,4}{
          \node[source] (sm\i) at (0, \i-0.5) {};
          \node[target] (tm\i) at (2.5, \i-0.5) {};
        }
        \foreach \i in {1,...,4}{
          \foreach \j in {1,...,4}{
            \pgfmathsetmacro{\op}{ifthenelse(\i==\j, 0.8, 0.3)}
            \pgfmathsetmacro{\lw}{ifthenelse(\i==\j, 2, 1)}
            \draw[plan, gray!70, line width=\lw pt, opacity=\op] (sm\i) to (tm\j);
          }
        }
        \node at (1.25,-0.5) {Medium $\epsilon$};
    \end{scope}

    \begin{scope}[xshift=7cm]
        \foreach \i in {1,...,4}{
          \node[source] (sl\i) at (0, \i-0.5) {};
          \node[target] (tl\i) at (2.5, \i-0.5) {};
        }
        \foreach \i in {1,...,4}{
          \foreach \j in {1,...,4}{
            \draw[plan, gray!60, line width=1.2pt, opacity=0.4] (sl\i) to (tl\j);
          }
        }
        \node at (1.25,-0.5) {Large $\epsilon$};
    \end{scope}
    \end{tikzpicture}
    \caption{Effect of the entropic regularization increasing the parameter $\epsilon$, from left to right the plan becomes smoother and more diffuse.}
    \label{fig:sinkhorn_effect}
\end{figure}

The classical OT problem can be computationally expensive to solve due to its stiffness, particularly for large-scale problems. Entropic regularization~\cite{cuturi13} addresses this by adding a strictly convex penalty term based on the Kullback-Leibler (KL) divergence, ensuring that the problem has a unique minimizer. The regularized OT problem can be formulated as:
\begin{equation}
\label{eq:ot_entropic}
    C_\epsilon(\mu, \nu) = \min_{\bm{\pi} \in \Pi(\mu, \nu)} \left\{ \langle \mathbf{C}, \bm{\pi} \rangle_F + \epsilon \, \text{KL}(\bm{\pi} \| \boldsymbol{a}\boldsymbol{b}^T) \right\},
\end{equation}
where the KL divergence for discrete measures used as penalty term is given by:
\begin{equation}
    \text{KL}(\bm{\pi} \| \boldsymbol{a}\boldsymbol{b}^T) = \sum_{i,j} \pi_{ij} \log\left(\frac{\pi_{ij}}{a_i b_j}\right).
\end{equation}
In Equation \eqref{eq:ot_entropic}, the regularization parameter $\epsilon > 0$ controls the trade-off between fidelity to the original OT problem (small $\epsilon$) and computational efficiency (large $\epsilon$). The effect of $\epsilon$ on the transport plan defined by the unique minimizer $\bm{\pi}_\epsilon$ is illustrated in Figure~\ref{fig:sinkhorn_effect}: for small $\epsilon$, the plan is sharp and concentrated along the optimal paths, while for larger $\epsilon$, it becomes increasingly diffuse, spreading mass more uniformly.

The main advantage of considering the regularized problem is that, exploiting results on the duality of linear minimization problems and the definitions of Kantorovich potentials~\cite{villani2021topics}, it admits a solution for which the optimal plan has the following simple multiplicative structure
\begin{equation}
\label{eq:sinkhorn_plan_form}
    \bm{\pi}_\epsilon = \mathrm{diag}(\mathbf{u}) \, \mathbf{K} \, \mathrm{diag}(\mathbf{v}), \quad \text{or component-wise,} \quad (\bm{\pi}_\epsilon)_{ij} = u_i K_{ij} v_j,
\end{equation}
where $\mathbf{K}$ is the Gibbs kernel defined as:
\begin{equation}
    \mathbf{K} = \exp(-\mathbf{C}/\epsilon), \quad \text{i.e.,} \quad K_{ij} = e^{-C_{ij}/\epsilon},
\end{equation}
and $\mathbf{u} \in \mathbb{R}_+^n$ and $\mathbf{v} \in \mathbb{R}_+^m$ are scaling vectors, that can be computed efficiently via the Sinkhorn algorithm~\cite{Schmitzer2019}. Starting from the initial vector $\mathbf{v}^{(0)} = \mathbf{1}_m$, this iterative fixed-point procedure alternates between row and column normalizations:
\begin{equation}
\label{eq:sinkhorn_iterations}
    \mathbf{u}^{(\ell+1)} = \boldsymbol{a} \oslash (\mathbf{K} \mathbf{v}^{(\ell)}), \qquad \mathbf{v}^{(\ell+1)} = \boldsymbol{b} \oslash (\mathbf{K}^T \mathbf{u}^{(\ell+1)}),
\end{equation}
where $\oslash$ denotes element-wise division, and converges to the unique scaling vectors satisfying the marginal constraints. We remark that the algorithm consists entirely of matrix-vector products, making it highly parallelizable and well-suited for GPU acceleration.

Displacement interpolation extends naturally to the entropy-regularized setting by replacing the optimal plan $\bm{\pi}^*$ with its regularized counterpart $\bm{\pi}_\epsilon$ in \eqref{eq:displacement_interp}:
\begin{equation}
\label{eq:discrete_disp_int_entropic}
    \mu_t^{(\epsilon)} = \sum_{i=1}^{n} \sum_{j=1}^{m} (\pi_\epsilon)_{ij} \, \delta_{(1-t)\x_i + t\y_j}.
\end{equation}

A critical requirement for applying OT is that the fields to be interpolated must be non-negative, as they are interpreted as probability mass distributions. Solution fields that are inherently non-negative, such as densities or concentrations, can be directly normalized to have unit mass. For fields containing both positive and negative values, such as general physical fields or residuals, a simple decomposition into positive and negative parts can be applied, with each part interpolated separately.

\subsection{OT-Based Displacement Interpolation for ROM}
\label{sec:ot_rom}

Finally, we summarize the core OT-ROM methodology introduced in our previous work~\cite{Khamlich2025JCPDI}, which serves as the foundation for the extensions presented in this paper. The central idea is to use displacement interpolation to generate synthetic snapshots between sparse high-fidelity checkpoints, addressing data scarcity in ROM construction.

Consider a time-dependent problem whose solution $\u(t)$ is available at a sparse set of $N_c$ checkpoint times $\{t_k\}_{k=1}^{N_c}$, with $N_c \ll N_T$. The goal is to construct a continuous-in-time approximation $\u_{\text{synth}}(t)$ for any $t \in [t_1, t_{N_c}]$ via the displacement interpolation framework. The procedure, allowing for the efficient evaluation of time-trajectory snapshots and data-augmentation tasks, is structured as follows:

\begin{enumerate}
    \item \textbf{Offline Stage: computation of the optimal transport plans}
    \begin{enumerate}
        \item For each checkpoint, normalize the field $\u(t_k)$ by its mass $m_k = \|\u(t_k)\|_{L^1(\Omega)}$ to obtain a probability measure.
        \item For each consecutive pair of checkpoints $(t_k, t_{k+1})$, compute the entropic OT plan $\bm{\pi}_\varepsilon^{k}$ between the normalized fields using the Sinkhorn algorithm.
    \end{enumerate}

    \item \textbf{Online Stage: prediction at time $t^*$}
    \begin{enumerate}
        \item For a query time $t^* \in [t_k, t_{k+1}]$, compute the local interpolation parameter playing the role of "virtual time" defined by $\alpha = (t^* - t_k) / (t_{k+1} - t_k)$.
        \item Interpolate the mass via $m(\alpha) = (1-\alpha) m_k + \alpha m_{k+1}$.
        \item Apply displacement interpolation using the pre-computed plan $\bm{\pi}_\varepsilon^{k}$:
        \begin{equation}
        \label{eq:synth_formula_base}
            \u_{\text{synth}}(\alpha) = m(\alpha) \sum_{\ell,m} (\bm{\pi}_\varepsilon^{k})_{\ell m} \, \delta_{\text{proj}_{\mathcal{X}}((1-\alpha)x_\ell + \alpha x_m)},
        \end{equation}
        where $\mathcal{X} = \{x_\ell\}_{\ell=1}^{N_h}$ denotes the discretization points and $\text{proj}_{\mathcal{X}}(\cdot)$ maps each point to its nearest grid location.
    \end{enumerate}
\end{enumerate}

This approach exploits the geometric structure of Wasserstein space to transport features between their positions at time $t_k$ and $t_{k+1}$, rather than simply blending them pointwise as in linear interpolation, which would be completely ineffective in the case of moving features with compact support. Moreover, the resulting synthetic snapshots preserve mass and follow geodesic paths in the space of probability measures, and thus, despite being a data-driven approach, this promotes  physical and geometric consistency with the underlying PDE.

\begin{rem}
Instead of choosing the standard linear time mapping for the virtual time $\alpha$, one could also consider the more accurate nonlinear time mapping introduced in~\cite{Khamlich2025JCPDI}, which minimizes the $L^2$ error between synthetic and available true snapshots. In general, accounting for multi-scale behavior in time could prove crucial for the accuracy of the online reconstruction, but for the problems considered here the linear mapping has been sufficient to provide accurate results.
\end{rem}

\begin{rem}
We remark that the accuracy of the OT-based interpolation can be further enhanced by a POD-based correction step~\cite{Khamlich2025JCPDI}. This involves constructing a POD basis for the dominant error modes and using Gaussian Process Regression to predict correction coefficients. To be able to assess the effectiveness of the plain displacement interpolation approach in the context of two-phase flows, we do not employ the correction in the present work.
\end{rem}

\section{Multi-Fidelity Residual Interpolation}
\label{sec:mf_ot_rom}

In many applications, one has access to models of varying fidelity, which could originate from different assumptions or discretization. A high-fidelity model provides accurate solutions $\uHF(t)$ but at substantial computational cost, while a low-fidelity solution $\uLF(t)$ can be evaluated cheaply but with reduced accuracy. In particular, the LF model might arise from a coarser spatial discretization (the setting considered in this work), a larger time step, or simplified physics. The goal of multi-fidelity modeling is to combine the strengths of both: exploiting the efficiency of LF models to obtain many simulation data to complete the missing information obtained from a limited number of HF simulations~\cite{Peherstorfer2018,ContiMultifidelityReducedorderSurrogate2023,KastNonintrusiveMultifidelityMethod2020,MorrisonGFNGraphFeedforward2024a}. In practice, this strategy allows for an efficient exploration of the solution manifold, both in the offline and in the online phases, exploiting a cheaper model to generate many snapshots and building a non-intrusive real-time surrogate, respectively.

A straightforward strategy is to model the residual or discrepancy field $\r(t) = \uHF(t) - \uLF(t)$, which captures the difference between the two models. The HF solution can then be approximated as $\uHF(t) \approx \uLF(t) + \r_{\text{approx}}(t)$, where $\r_{\text{approx}}(t)$ is a reduced-order and efficient-to-evaluate representation of the residual, that, for example, can be constructed via a POD-regression strategy on the residual snapshots. However, for problems with moving features, the residual field itself exhibits advective dynamics that are poorly captured by a linear combination in the POD coefficient space, motivating our OT-based approach.

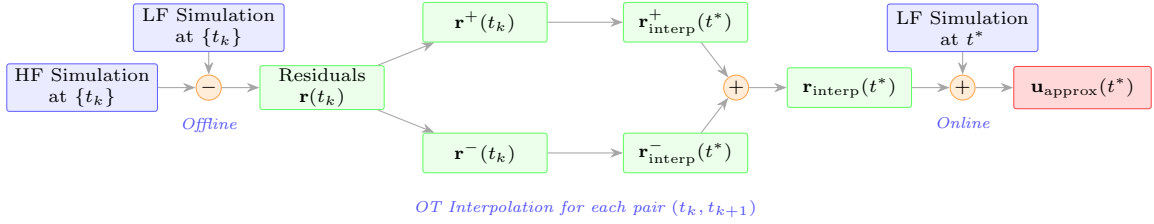
\begin{figure}[t!]
    \centering
    \begin{tikzpicture}[
        node distance=0.4cm and 0.7cm,
        >=Stealth,
        block/.style={
            rectangle,
            draw=blue!60,
            fill=blue!8,
            text width=5.5em,
            text centered,
            rounded corners=1pt,
            minimum height=1.8em,
            font=\scriptsize,  
            inner sep=1pt
        },
        op/.style={
            circle,
            draw=orange!70,
            fill=orange!15,
            text centered,
            minimum size=0.6em,
            font=\scriptsize\bfseries,  
            inner sep=0.5pt
        },
        data/.style={
            rectangle,
            draw=green!60,
            fill=green!8,
            text centered,
            rounded corners=1pt,
            minimum height=1.5em,
            font=\scriptsize,  
            inner sep=1pt,
            text width=4.5em
        },
        final/.style={
            rectangle,
            draw=red!70,
            fill=red!15,
            text centered,
            rounded corners=1pt,
            minimum height=1.6em,
            font=\scriptsize\bfseries,  
            inner sep=1pt,
            text width=5em
        },
        arrow/.style={
            ->,
            color=gray!70
        }
    ]

        \node[block] (hf_sim) {HF Simulation\\at $\{t_k\}$};
        
                \node[op, right=0.5cm of hf_sim] (subtract) {$-$};

        \node[block, above=0.3cm of subtract] (lf_sim) {LF Simulation\\at $\{t_k\}$};

        \node[data, right=0.5cm of subtract] (residual) {Residuals\\$\mathbf{r}(t_k)$};

        \node[data, above right=0.3cm and 0.5cm of residual] (res_pos) {$\mathbf{r}^+(t_k)$};
        \node[data, below right=0.3cm and 0.5cm of residual] (res_neg) {$\mathbf{r}^-(t_k)$};

        \node[data, right=1cm of res_pos] (interp_pos) {$\mathbf{r}^+_{\text{interp}}(t^*)$};
        \node[data, right=1cm of res_neg] (interp_neg) {$\mathbf{r}^-_{\text{interp}}(t^*)$};

        \node[below=0.5cm, font=\tiny\itshape, color=blue!70] at ($(res_neg)!0.5!(interp_neg)$) {OT Interpolation for each pair $(t_k, t_{k+1})$};

        \node[op, right=0.5cm of $(interp_pos)!0.5!(interp_neg)$] (recombine) {$+$};

        \node[data, right=0.5cm of recombine] (interp_res) {$\mathbf{r}_{\text{interp}}(t^*)$};

        \node[op, right=0.5cm of interp_res] (add) {$+$};

        \node[block, above=0.3cm of add] (lf_online) {LF Simulation\\at $t^*$};

        \node[final, right=0.5cm of add] (final_sol) {$\mathbf{u}_{\normalfont{\text{approx}}}(t^*)$};

        \draw[arrow] (hf_sim) -- (subtract);
        \draw[arrow] (lf_sim) -- (subtract);
        \draw[arrow] (subtract) -- (residual);
        \draw[arrow] (residual) -- (res_pos);
        \draw[arrow] (residual) -- (res_neg);
        \draw[arrow] (res_pos) -- (interp_pos);
        \draw[arrow] (res_neg) -- (interp_neg);
        \draw[arrow] (interp_pos) -- (recombine);
        \draw[arrow] (interp_neg) -- (recombine);
        \draw[arrow] (recombine) -- (interp_res);
        \draw[arrow] (interp_res) -- (add);
        \draw[arrow] (lf_online) -- (add);
        \draw[arrow] (add) -- (final_sol);

        \node[below=0.1cm of subtract, font=\tiny\itshape, color=blue!70] {Offline};
        \node[below=0.1cm of add, font=\tiny\itshape, color=blue!70] {Online};

    \end{tikzpicture}
    \caption{Schematic of the multi-fidelity OT-ROM framework, using displacement interpolation on the  residual between HF and LF models.}
    \label{fig:mf_ot_rom_scheme}
\end{figure}

Thus, we propose to use displacement interpolation to model the temporal evolution of the residual field, developing the Multi-Fidelity OT-ROM (MF-OT-ROM) strategy illustrated in Figure~\ref{fig:mf_ot_rom_scheme}. 

Given HF and LF snapshots at checkpoint times $\{t_k\}_{k=1}^{N_c}$, we compute residual checkpoints $\r_k = \uHF(t_k) - \uLF(t_k)$. In practice, for resolution-based multi-fidelity approaches, the LF and HF models are typically defined on different spatial grids. In all the numerical experiments considered here, we employ Cartesian grids with refinement factors of $2$ or $4$, therefore, the LF solution is interpolated onto the HF grid so that the residual can be evaluated on a common spatial discretization. We then apply OT-based interpolation to generate $\r_{\text{interp}}(t^*)$ for any query time $t^*$. Since residual fields typically contain both positive and negative values, we apply a simple decomposition into positive and negative parts, $\r = \r^+ - \r^-$, with each part interpolated separately before recombination, so that the final multi-fidelity approximation is given by 
\begin{equation}
    \u_{\text{approx}}(t^*) = \uLF(t^*) + \r_{\text{interp}}(t^*).
\end{equation}
This multi-fidelity approach offers clear advantages over simply applying standard interpolation (\eg, linear or spline) to residual coefficients in a POD basis, especially in this context. Indeed, the OT-based interpolation inherently respects the geometric structure and movement of features within the residual field, and thus, for diffused-interface models, where the residual is often concentrated at the moving interface, the method can effectively ``transport'' the correction from its position at $t_k$ to its new position at $t_{k+1}$, rather than just fading out the old correction and fading in the new one.

\section{Parametric Displacement Interpolation}
\label{sec:pmf_ot_rom}

In the previous sections, we have specifically considered time-dependent problems as the natural framework where to exploit our strategy. Nevertheless, being a data-driven strategy, both temporal and parametric nature of the moving features can be treated in a similar way, thus we are interested in extending the strategy proposed in~\cite{Khamlich2025JCPDI} to perform displacement interpolation for problems with parametric dependencies. More specifically, in the formulation presented so far, time $t$ serves as the interpolation coordinate in the solution manifold, with OT plans computed between consecutive temporal checkpoints, while we aim at including additional parameters $\mu \in \Pcal \subset \R^p$, modeling physical characteristics of the systems, such as material properties and boundary conditions, or geometric properties.

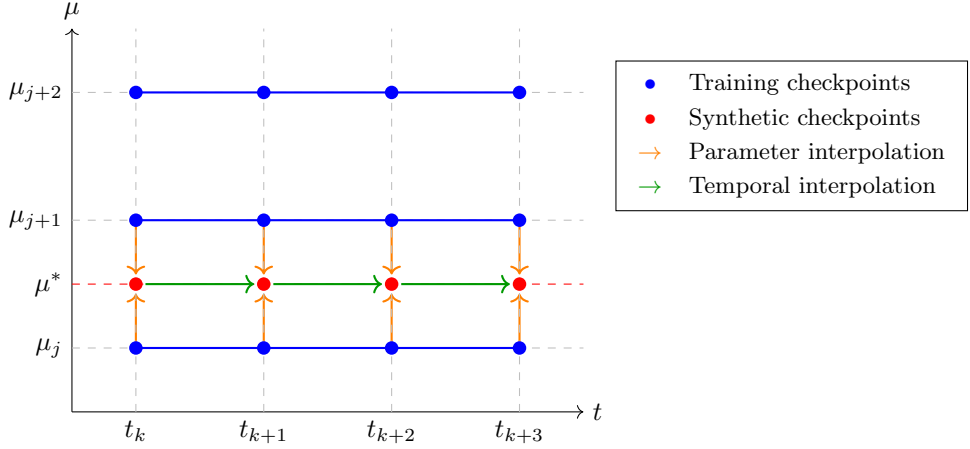
\begin{figure}[btp]
    \centering
    \begin{tikzpicture}[scale=0.85]
        \draw[->] (0,0) -- (8,0) node[right] {$t$};
        \draw[->] (0,0) -- (0,6) node[above] {$\mu$};

        \foreach \y in {1,3,5} {
            \draw[dashed, gray!50] (0,\y) -- (8,\y);
        }
        \node[left] at (0,1) {$\mu_j$};
        \node[left] at (0,3) {$\mu_{j+1}$};
        \node[left] at (0,5) {$\mu_{j+2}$};

        \foreach \x in {1,3,5,7} {
            \draw[dashed, gray!50] (\x,0) -- (\x,6);
        }
        \node[below] at (1,0) {$t_{k}$};
        \node[below] at (3,0) {$t_{k+1}$};
        \node[below] at (5,0) {$t_{k+2}$};
        \node[below] at (7,0) {$t_{k+3}$};

        \foreach \y in {1,3,5} {
            \foreach \x in {1,3,5,7} {
                \fill[blue] (\x,\y) circle (3pt);
            }
            \draw[blue, thick] (1,\y) -- (3,\y) -- (5,\y) -- (7,\y);
        }

        \draw[red, dashed] (0,2) -- (8,2);
        \node[left] at (0,2) {$\mu^*$};

        \begin{scope}[on background layer]
        \foreach \x in {1,3,5,7} {
            \draw[->, thick, orange] (\x,1) -- (\x,1.85);
            \draw[->, thick, orange] (\x,3) -- (\x,2.15);
        }
        \end{scope}

        \foreach \x in {1,3,5,7} {
            \fill[red] (\x,2) circle (3pt);
        }

        \draw[->, thick, green!60!black] (1.15,2) -- (2.85,2);
        \draw[->, thick, green!60!black] (3.15,2) -- (4.85,2);
        \draw[->, thick, green!60!black] (5.15,2) -- (6.85,2);

        \node[draw, fill=white, anchor=north west, font=\small] at (8.5,5.5) {
            \begin{tabular}{cl}
                {\color{blue}$\bullet$} & Training checkpoints\\[2pt]
                {\color{red}$\bullet$} & Synthetic checkpoints\\[2pt]
                {\color{orange}$\rightarrow$} & Parameter interpolation\\[2pt]
                {\color{green!60!black}$\rightarrow$} & Temporal interpolation
            \end{tabular}
        };
    \end{tikzpicture}
    \caption{Parametric and temporal displacement interpolation on the $(t, \mu)$ grid.} 
    \label{fig:parametric_strategy}
\end{figure}

In Figure~\ref{fig:parametric_strategy}, we present the procedure to exploit a dataset of dynamical systems corresponding to different parameter values performing a two-stage displacement interpolation strategy to obtain synthetic checkpoints for newly given parameters and recover their temporal evolution. Practically speaking, given solution snapshots at scalar ($p=1$) training parameter values $\{\mu_j\}_{j=1}^{N_p}$ and checkpoint times $\{t_k\}_{k=1}^{N_c}$, we can interpolate along both coordinate directions: temporally at fixed parameter, or parametrically at fixed time. Following this latter option, for a test parameter $\mu^* \in (\mu_j, \mu_{j+1})$, OT plans are computed between snapshots at adjacent training parameters for each checkpoint, and synthetic snaphshots are generated via parametric displacement interpolation with ``virtual parameter'' $\alpha_\mu = (\mu^* - \mu_j)/(\mu_{j+1} - \mu_j)$. Next, temporal displacement interpolation is exploited to connect these generated checkpoints, and obtain efficient and reliable predictions for any query time. We remark that the extension to higher-dimensional parameter spaces can be done in several ways, exploiting additional knowledge on the parametric dependency, considering tensor product of grids. Since all these extensions could potentially require a dramatic increase of the computational cost when the number of parameters increases, we postpone the investigation of an efficient strategy to following works.

By combining the multi-fidelity framework from Section~\ref{sec:mf_ot_rom} with the parametric setting, here we also aim at introducing the Parametric Multi-Fidelity OT-ROM (PMF-OT-ROM). We start by noticing that two strategies are possible, depending on the field being considered and the characteristics of the problem at hand. The first approach entails applying parametric displacement interpolation directly to the HF training snapshots data to generate synthetic HF checkpoints $\uHF^{\text{synth}}(t_k; \mu^*)$. A cheap LF simulation is run at $\mu^*$, and the synthetic residuals $\r^{\text{synth}}(t_k; \mu^*) = \uHF^{\text{synth}}(t_k; \mu^*) - \uLF(t_k; \mu^*)$ are interpolated in time to produce the final approximation $\u_{\text{approx}}(t; \mu^*) = \uLF(t; \mu^*) + \r_{\text{interp}}(t; \mu^*)$. 
On the other hand, one could compute residuals at each training parameter value $\r^{\text{synth}}(t_k; \mu_j) = \uHF^{\text{synth}}(t_k; \mu_j) - \uLF(t_k; \mu_j)$, apply parametric displacement interpolation to these residual fields to generate synthetic residual checkpoints $\r^{\text{synth}}(t_k; \mu^*)$, and finally perform temporal displacement interpolation to have a continuous corrected model in time. 

Both strategies avoid constructing a monolithic ROM over the joint $(t, \mu)$ space, rather they explore the displacement interpolation along the two coordinates shifting the focus from the physical to the residual fields. The former approach is conceptually simpler, especially for applications where the high-fidelity solution has relevant physical meaning, while residual interpolation could be less intuitive, but preserving the same characteristics of moving features, it may be advantageous when its behavior with respect to $\mu$ varies more smoothly than the full solution. 

In terms of computational cost, it is important to highlight two advantages of the former approach (i.e., interpolating HF solutions first and then computing the residual). First, the DI step needs to be performed only three times instead of four, since the phase-field is strictly non-negative and therefore does not require a prior decomposition into positive and negative parts. Second, this approach only requires HF solutions along the neighboring parametric trajectories, while for the unseen parameter $\mu^{*}$ only the LF solution is needed. Instead, interpolating residuals in both time and parameter space requires LF solutions also for the neighboring trajectories. These additional computations ar thus essential to reconstruct the residual associated with the trajectory corresponding to $\mu^{*}$. For these reasons, in the numerical experiments considered in this work we adopt this approach within the parametric multi-fidelity strategy.

A key assumption underlying the proposed approach is that the solution (or the residual field) varies smoothly with respect to the parameter $\mu$. In problems where the parametric dependence exhibits bifurcations~\cite{PichiArtificialNeuralNetwork2023,HessLocalizedReducedorderModeling2019} or sharp transitions in parameter space~\cite{HesthavenNonlinearModelReduction2026,RimModelReductionTransportDominated2023}, the interpolation-based strategy may become less effective. In such situations, more sophisticated techniques such as local model order reduction~\cite{AmsallemNonlinearModelOrder2012,GeelenLocalizedNonintrusiveReducedorder2022} or adaptive sampling strategies~\cite{PeherstorferModelReductionTransportDominated2020,JinAdaptiveHybridReduced2025} may be required.

\section{Problem Definition: Five equation model}
\label{sec:problem_def}
The proposed methodological extensions of OT-based ROMs are applied in this work to the numerical simulation of compressible two-phase flows using diffuse-interface methods. Specifically, the full-order model relies on the five-equation formulation originally introduced by Allaire et al.~\cite{allaire2002five}.
From the Baer-Nunziato seven-equation model~\cite{baer1986two}, assuming mechanical and thermodynamic equilibrium between the two phases, the following five-equation model can be derived:
\begin{align}
 \frac{\partial \psi}{\partial t} + \textbf{v} \cdot \nabla \psi &= 0, \label{eq:5eq_finalPC2}\\
 \frac{\partial \phi_{1}}{\partial t} + \textbf{v} \cdot \nabla \phi_{1} &= \nabla \cdot \textbf{a}_{1} \label{eq:AC}, \\
 \frac{\partial (\rho_{l} \phi_{l})}{\partial t} + \nabla \cdot (\rho_{l} \phi_{l} \textbf{v}) &= \nabla \cdot \textbf{R}_{l},  \quad l=\{1,2\}, \label{eq:5eq_mass}  \\
 \frac{\partial (\rho \textbf{v})}{\partial t} + \nabla \cdot (\rho \textbf{v}  \otimes \textbf{v} + p \mathbb{I}) &= \nabla \cdot (\textbf{f} \otimes \textbf{v}) + \nabla \cdot \boldsymbol{\tau} + \sigma \kappa \nabla \phi_{1}, \label{eq:5eq_mom} \\
 \frac{\partial (\rho E)}{\partial t} + \nabla \cdot [(\rho E + p) \textbf{v}] &= \nabla \cdot (\textbf{f} g) + \sum_{l=1}^{2} \nabla \cdot ( \rho_{l} H_{l} \textbf{a}_{l}) + \nabla \cdot (\boldsymbol{\tau} \cdot \textbf{v}) + \sigma \kappa \nabla \phi_{1} \cdot \textbf{v}, \label{eq:5eq_energy}
\end{align}
consisting of one evolution equation for the phase-field (\ref{eq:AC}), where $\phi_{1}$ denotes the phase-field referred to the first phase, two separate mass conservation laws for the individual phases (\ref{eq:5eq_mass}), one balance equation for the mixture total momentum (\ref{eq:5eq_mom}) and one balance equation for the total energy of the mixture (\ref{eq:5eq_energy}) with $H_{l}$ the specific enthalpy of the $l^{\mathrm{th}}$ phase. Compared to the original model, the present approach incorporates an additional transport equation for the level-set function $\psi$ (\ref{eq:5eq_finalPC2}), which is used to evaluate the interface normal vectors, following a strategy similar to that proposed by Al-Salami et al.~\cite{al2021high}. The variable $\rho$ is the mixture density, $\rho \textbf{v}$ the total momentum, and $p$ the pressure. The total energy is defined as $\rho E = \rho e + \frac{1}{2}\rho \|\textbf{v}\|^{2}$, i.e., the sum of internal and kinetic energy. The parameter $\sigma$ denotes the surface tension coefficient, $\kappa = -\nabla \cdot \widehat{\textbf{n}}_{1}$ the interface curvature, and the following auxiliary quantities are commonly introduced:
\begin{align}
& \textbf{a}_{l} = \Gamma\left[ \varepsilon \nabla \phi_{l} - \phi_{l} (1- \phi_{l})\widehat{\textbf{n}}_{l}\right], 
\qquad \widehat{\textbf{n}}_{l} = \frac{\nabla \psi}{\|\nabla \psi\|}, \\
& \textbf{R}_{l} = \rho_{l}^{(0)} \textbf{a}_{l}, 
\qquad \textbf{f} = \sum_{l=1}^{2} \textbf{R}_{l}, 
\qquad g=\frac{1}{2} \|\textbf{v}\|^{2}.
\end{align}
The regularization term $\mathbf{a}_{1}$ in Equation \eqref{eq:AC} is related to the conservative Allen-Cahn formulation proposed by~\cite{chiu2011conservative}, and is governed by two main parameters, $\varepsilon$ and $\Gamma$. The former controls the interface thickness and is typically chosen to be proportional to the grid size, while the latter regulates the magnitude of the regularization relative to advection and is commonly scaled with the maximum velocity magnitude in the computational domain. All the additional terms involving contributions that depend directly on $\mathbf{a}_{l}$ are classically added to the remaining conservation laws in order to ensure consistency and preserve kinetic energy and entropy in the overall coupled system~\cite{jain2020conservative}.

The viscous stress tensor is defined as
\[
\boldsymbol{\tau}= 2 \mu \left(\mathbb{S} - \frac{1}{3}(\nabla \cdot \textbf{v}) \mathbb{I}\right),
\]
where $\mu$ is the mixture dynamic viscosity, evaluated as $\mu = \mu_{1} \phi_{1} + \mu_{2} \phi_{2}$, and $\mathbb{S} = [\nabla \textbf{v} + (\nabla \textbf{v})^{\intercal}]/2$ denotes the strain-rate tensor.

The system is closed by relating the internal energy to the pressure through a stiffened-gas equation of state (EOS)~\cite{harlow1971fluid}:
\begin{equation}\label{eq:EOS}
p = \frac{\rho e - \left ( \frac{ \gamma_{1} p^{\infty}_{1}}{\gamma_{1}-1} \phi_{1} + \frac{ \gamma_{2} p^{\infty}_{2}}{\gamma_{2}-1} \phi_{2} \right)}{ \left ( \frac{\phi_{1}}{\gamma_{1}-1} +   \frac{\phi_{2}}{\gamma_{2}-1} \right)},
\end{equation}
where $\gamma_{l}$ and $p^{\infty}_{l}$ are the EOS parameters of the  $l^{\mathrm{th}}$ phase. The corresponding speed of sound and specific enthalpy for each phase are given by
\begin{equation}
c_{l} = \sqrt{\gamma_{l}\left( \frac{p + p^{\infty}_{l}}{\rho_{l}} \right)} 
\qquad \text{and} \qquad 
H_{l} = \frac{(p + p^{\infty}_{l})\gamma_{l}}{\rho_{l} (\gamma_{l} -1)}, 
\quad l =\{1,2\}.
\end{equation}
For completeness, the following mixture relations hold:
\begin{align}
\phi_{2} &= 1- \phi_{1},   \\
\rho &= \rho_{1}\phi_{1} + \rho_{2} \phi_{2},  \\
\frac{1}{\gamma-1} &= \phi_{1} \frac{1}{\gamma_{1}-1} + \phi_{2} \frac{1}{\gamma_{2}-1}, \\
p^{\infty} \frac{\gamma}{\gamma-1} &= \phi_{1} \frac{\gamma_{1} p_{1}^{\infty}}{\gamma_{1}-1} + \phi_{2} \frac{\gamma_{2} p_{2}^{\infty}}{\gamma_{2}-1}.
\end{align}

The spatial discretization of the five-equation model, in both kinetic and fully coupled configurations, is performed using a high-order Spectral Difference scheme (see~\cite{tonicello2025extension} for further details regarding the numerical discretization of the five-equation model).

In kinetic tests, the phase-field evolves under a prescribed analytical velocity field, whereas in the fully coupled configuration the complete system of governing equations is solved. As representative kinetic benchmark, we consider the Rider-Kothe vortex~\cite{rider1998reconstructing}, while the Rayleigh-Taylor instability~\cite{tryggvason1988numerical} is adopted as the fully coupled test case. Both problems are inherently unsteady and involve explicit parametric dependencies in their initial conditions.

In all cases, the analysis focuses on the accurate prediction of interfacial dynamics. The primary quantity of interest is consequently represented by the phase-field variable $\phi_{1}$. Although the same framework can be extended to other flow variables with different physical interpretations, $\phi_{1}$ is selected due to its central role in the modeling of two-phase compressible flows. Furthermore, for the sake of simplicity, we will refer to $\phi_{1}$ simply as $\phi$ for the remainder of the paper.

\section{Numerical Results}
\label{sec:numerical_results}
In this section, we assess the performance of the proposed MF-OT-ROM and PMF-OT-ROM frameworks on benchmark problems based on the diffuse-interface models described in Section~\ref{sec:problem_def}. In order to present a comprehensive evaluation of the performance of both frameworks, for the MF-OT-ROM, we demonstrate the accuracy of the \textit{corrected} residual interpolation and the improvement over the standalone \textit{uncorrected} LF solution for single trajectories involving multiple fidelities, while for the PMF-OT-ROM, we evaluate the model's ability to improve LF solutions for completely unseen parametric trajectories.
\subsection{Rider-Kothe vortex}
A classical test case for interface capturing techniques consists in the deformation of a two-dimensional liquid column by a shear flow~\cite{rider1998reconstructing}. In particular, in the classical test, a droplet of radius $R=0.15$ centered at $(0.5,0.75)$ in a $[0,1]^{2}$ square periodic domain is advected by the following divergence-free velocity field:
\begin{equation}
\textbf{v} = \begin{bmatrix} v_{1}\\[0.1cm] v_{2} \end{bmatrix} = \begin{bmatrix}
    \sin^{2} (\pi x_{1}) \sin (2 \pi x_{2}) \cos \left( \frac{\pi t}{T}\right)\\[0.1cm]
    -\sin ( 2 \pi x_{1}) \sin^{2} (\pi x_{2}) \cos \left( \frac{\pi t}{T}\right)
\end{bmatrix},
\end{equation}
where $T$ denotes the characteristic period of the shear flow, which we have chosen as the classical value $T=4$ for this particular test. Under the action of such a velocity field, during the first half of the period the initial circular bubble is strongly deformed into a thin filament. After the velocity reversal, the filament is advected back toward its initial configuration. 

This benchmark already highlights several challenges from the full-order perspective. At the full-order level, the strong deformation of the interface represents a well-known difficulty for interface-capturing techniques, which should ideally preserve the interface as sharp as possible while satisfying boundedness requirements and maintaining numerical stability. 

Challenges also emerge from a model order reduction standpoint, as the presence of a sharp interface moving rapidly across the computational domain is a typical feature of advection-dominated flows. As a result, the need for highly refined grids to resolve thin interface structures, combined with the pronounced variability between snapshots at different time instances, can significantly affect the performance of conventional ROM approaches. By contrast, displacement interpolation is anticipated to deliver consistent and accurate predictions, producing synthetic yet physically meaningful approximations at a negligible computational cost.

For this problem, we consider as the Full Order Model (FOM) a $3^{\mathrm{rd}}$-order polynomial spatial discretization on a series of increasingly finer meshes from $128^2$ to $512^2$ degrees of freedom, while we performed time integration using a classical $4^{\mathrm{th}}$ order Runge-Kutta scheme. For this first test case, we consider a dataset consisting of $401$ snapshots, corresponding to a time window of $\Delta \tau = 0.01$. The time step sizes used in the simulations are approximately $\Delta t_{128} \approx 10^{-3}$, $\Delta t_{256} \approx 5 \cdot 10^{-4}$, and $\Delta t_{512} \approx 2.5 \cdot 10^{-4}$ for the three spatial resolutions, to satisfy CFL constraints.

\begin{figure}[h!]
\centering
\subfigure[$R=0.1$]{\includegraphics[width=0.32\textwidth]{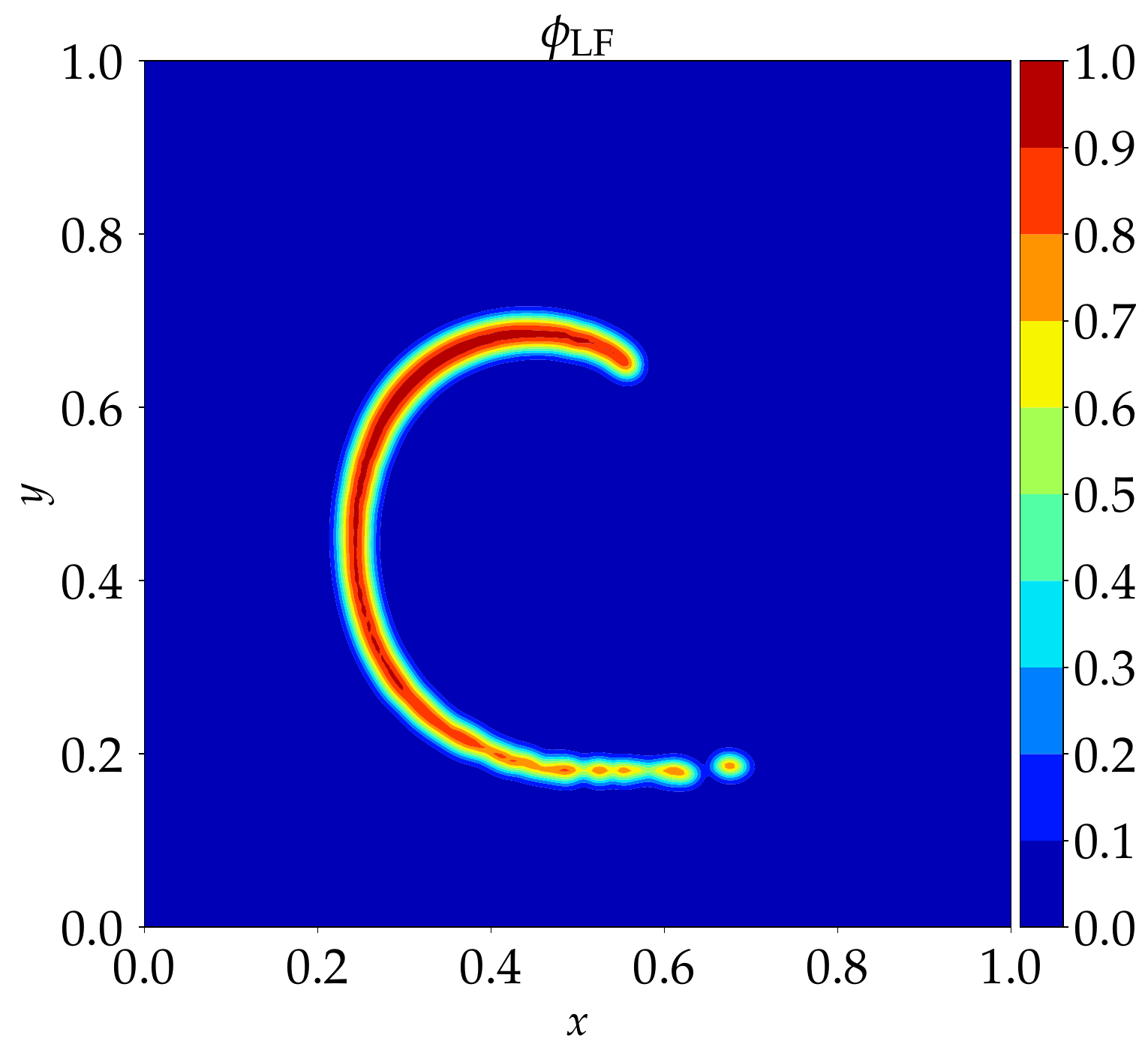}}
\subfigure[$R=0.15$]{\includegraphics[width=0.32\textwidth]{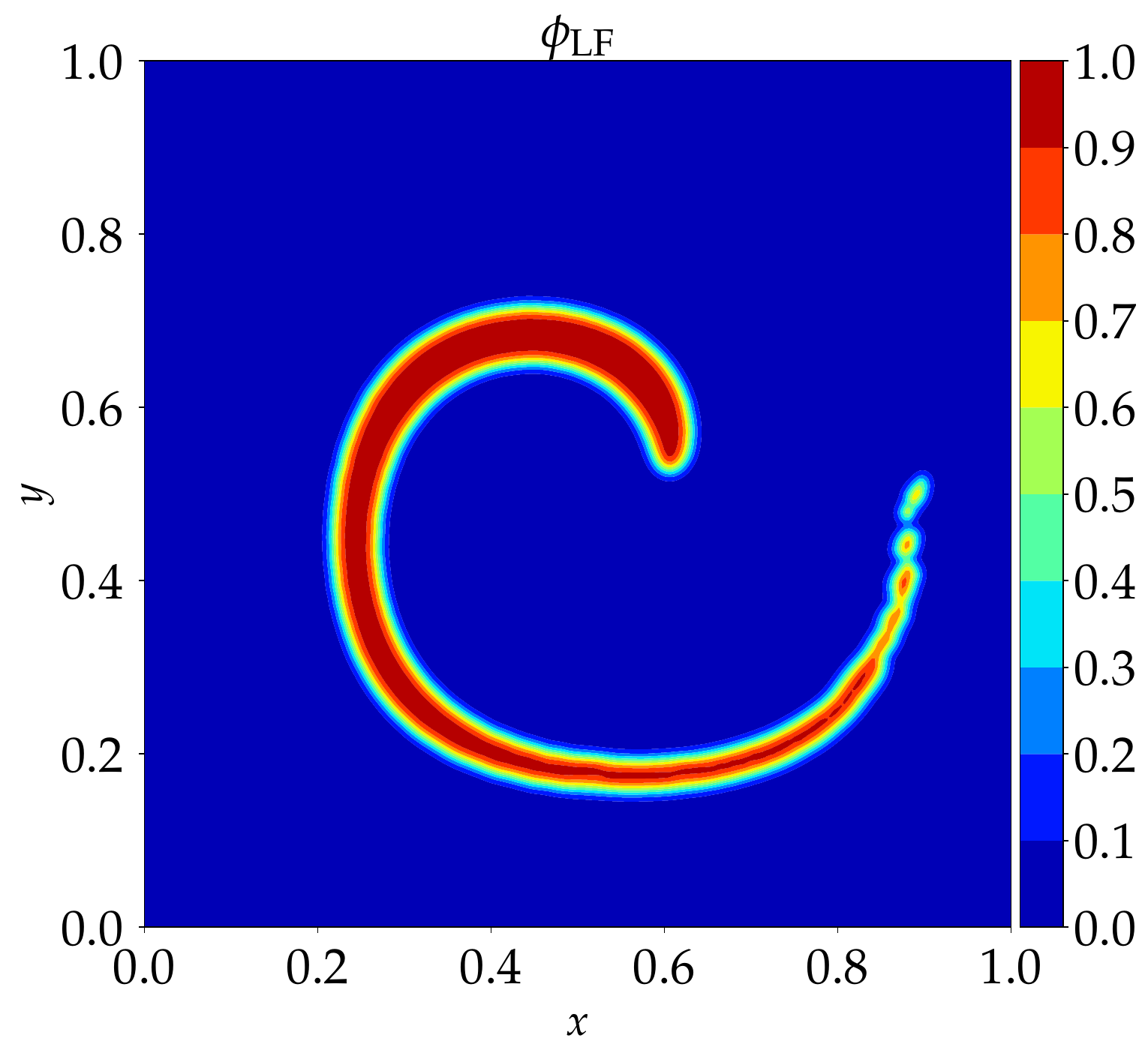}}
\subfigure[$R=0.2$]{\includegraphics[width=0.32\textwidth]{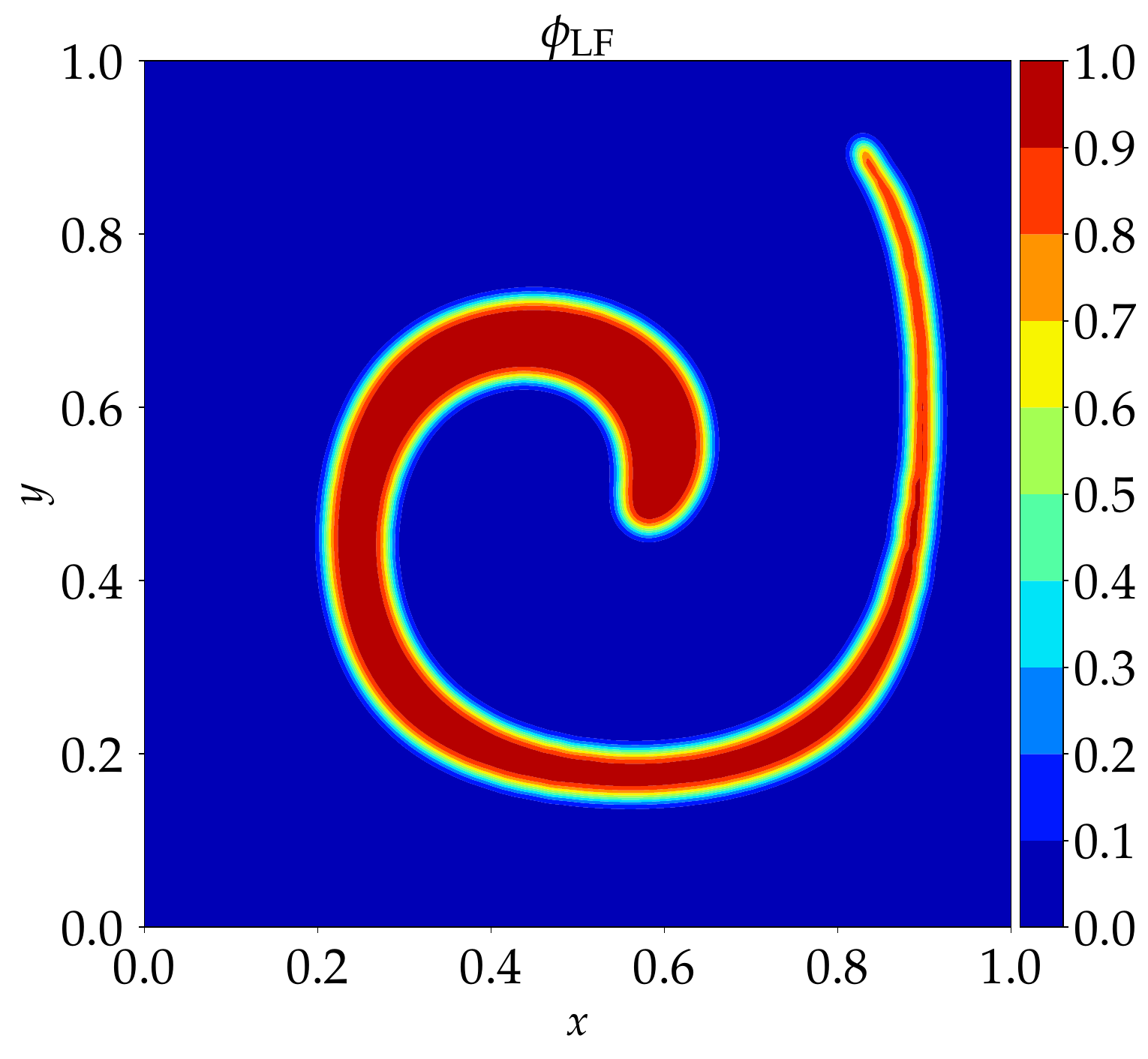}}
\subfigure[$R=0.1$]{\includegraphics[width=0.32\textwidth]{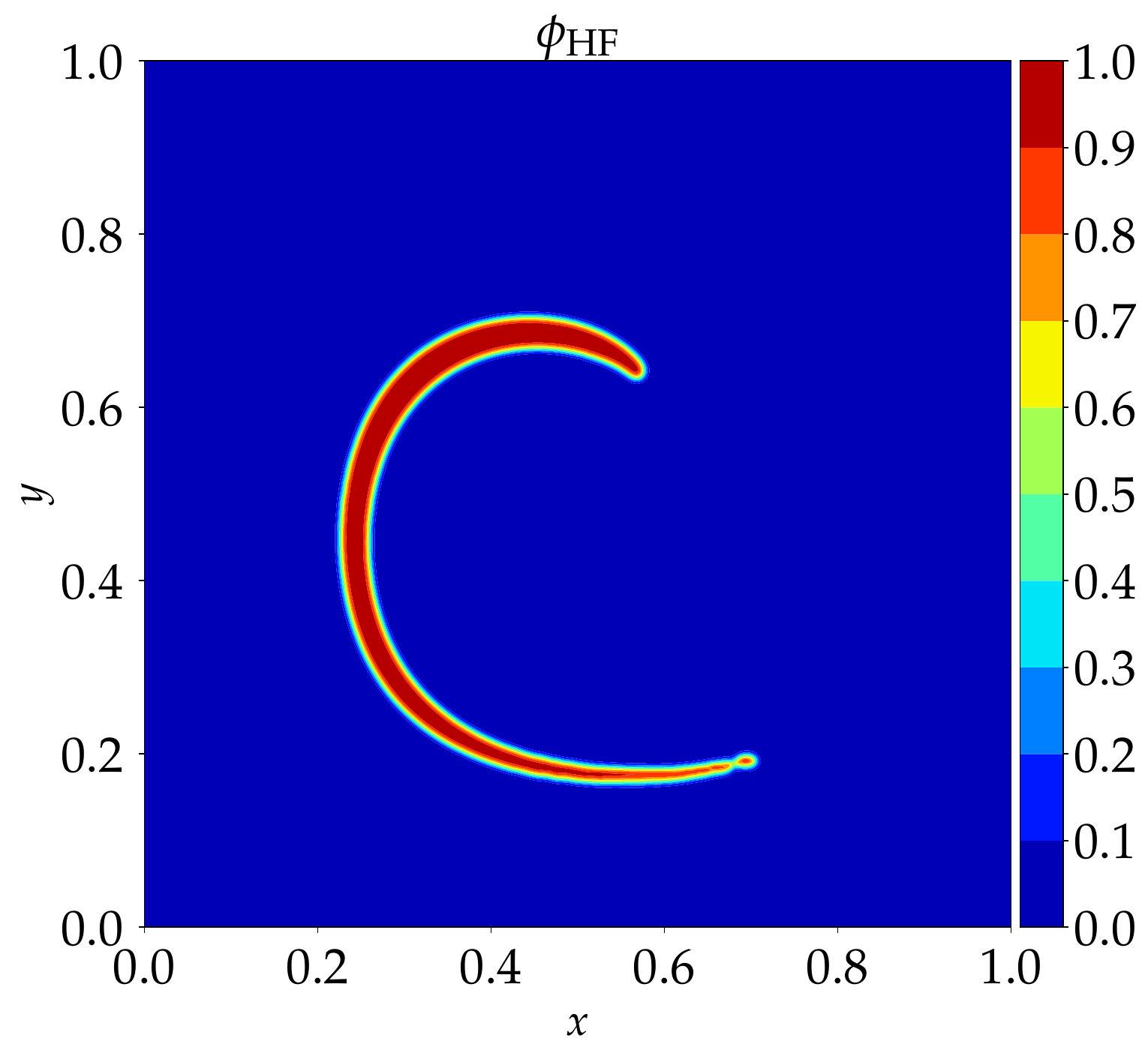}}
\subfigure[$R=0.15$]{\includegraphics[width=0.32\textwidth]{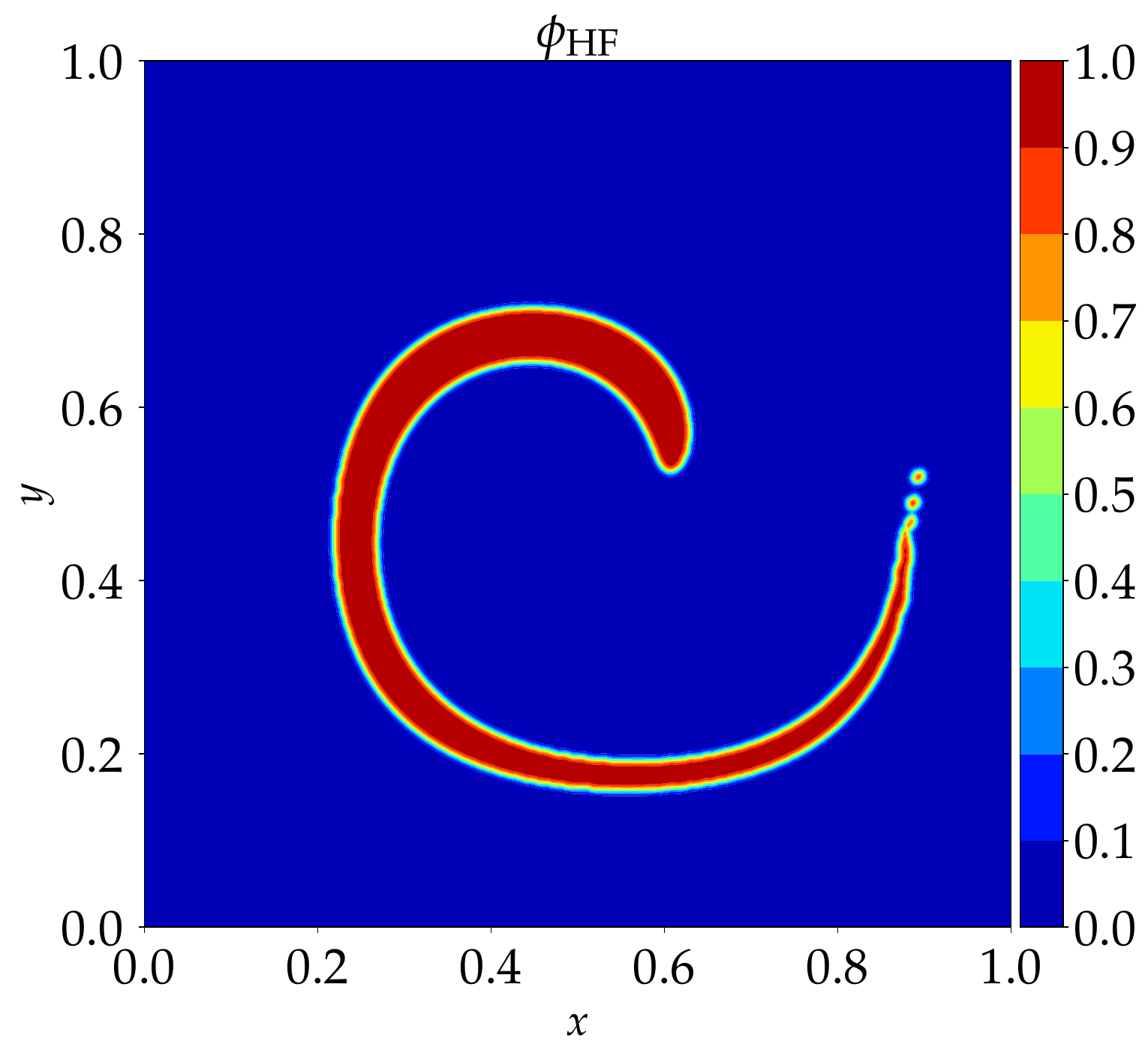}}
\subfigure[$R=0.2$]{\includegraphics[width=0.32\textwidth]{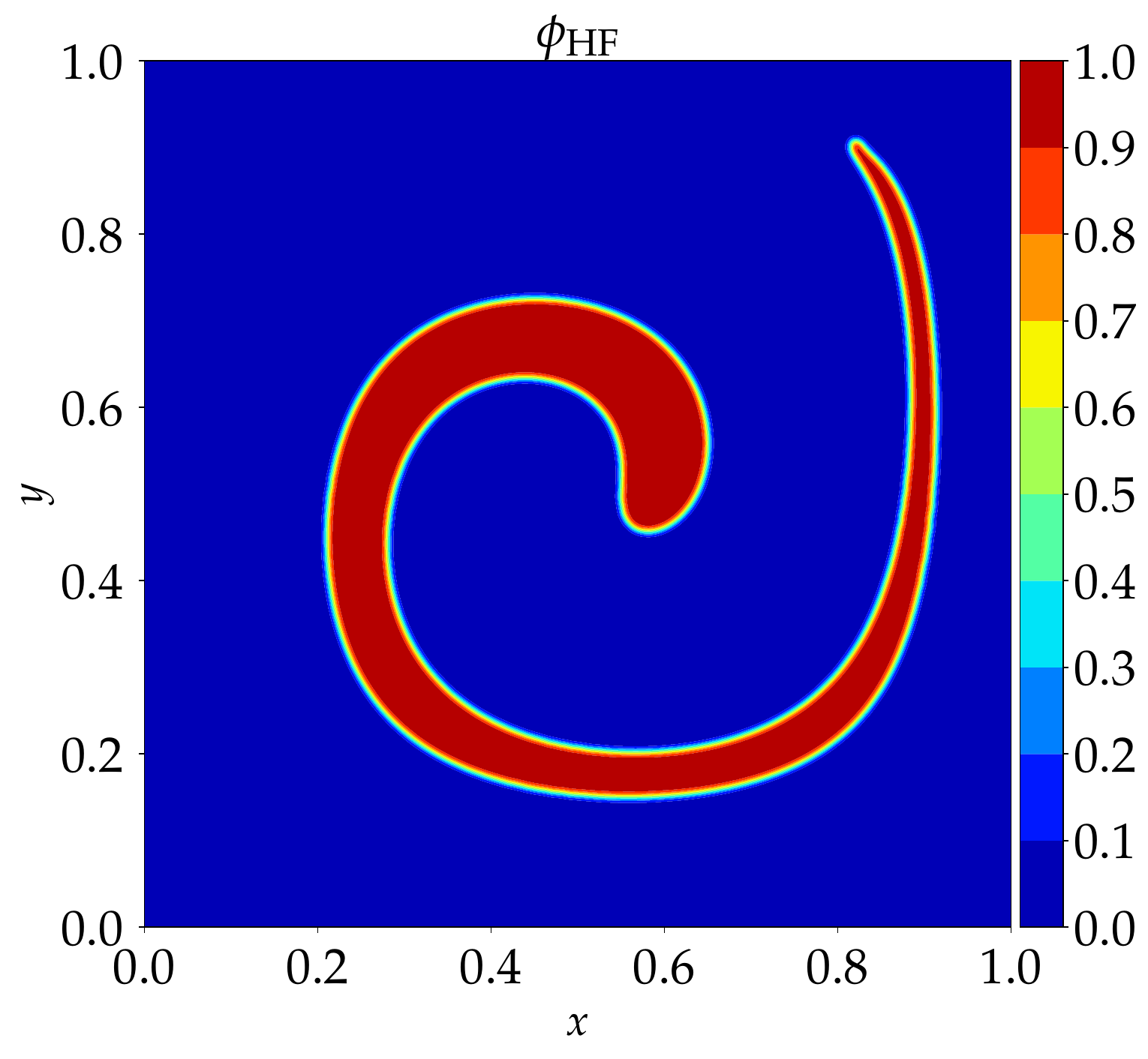}}
\caption{LF (top) and HF (bottom) phase-field isocontours for different initial conditions of the Rider-Kothe vortex at $t^{*}=T/2$.}
\label{fig:RK_overview}
\end{figure}

Moreover, we propose the following parametric extension of the test case by changing the initial radius of the bubble as:
\begin{equation}
    R_{k}=0.1 + k \Delta R \quad \mathrm{with} \quad \Delta R = 0.005 \quad \mathrm{and} \quad k=0,\dots,20.
\end{equation}
An overview of the problem's dynamics corresponding to different values for the initial radius $R$ and two different levels of refinement (LF with $256^{2}$ and HF with $512^{2}$) is shown in Figure~\ref{fig:RK_overview}. If not specified otherwise, we will always consider the case $R=0.11$ to test our framework for the Rider-Kothe vortex in the non-parametric scenario.
%
\subsubsection{Proper Orthogonal Decomposition and Displacement Interpolation}
Before presenting the performance of the novel strategies proposed in this work, we want to further motivate and numerically show the main reasons to exploit displacement interpolation instead of a classical POD-based reduced-order models for advection-dominated flows. Here, we aim to highlight the limitations of standard linear reduction strategies in reconstructing the temporal evolution of the phase-field in two-phase flows. Towards this goal, we consider the aforementioned Rider-Kothe vortex problem and compare the results obtained with POD projection and displacement interpolation of the HF solutions.

We remark that we have introduced the nomenclature \emph{checkpoints} for the set of snapshots used to perform the pairwise displacement interpolation. Thus, to provide a fair comparison with POD, the number of checkpoints used in the displacement interpolation will always match the number of snapshots available for the POD projection. Indeed, the checkpoints are the only \emph{data} nominally available to both OT- and POD-based methodologies.

\begin{figure}[t!]
\centering
\includegraphics[width=0.96\textwidth]{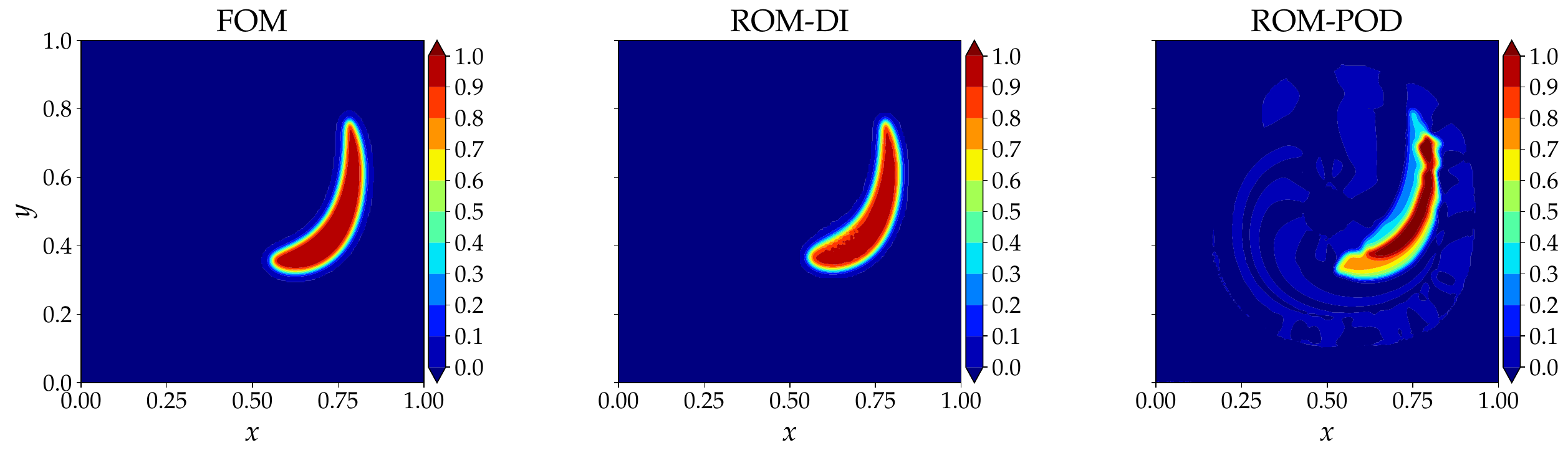}
\caption{Phase-field at $t\approx 0.2$ for the Rider-Kother vortex at resolution $256^2$ with FOM (left), DI (center) and POD (right) for $N_c=30$.}
\label{fig:err_POD_RK}
\end{figure}

Figure~\ref{fig:err_POD_RK} shows an example of the predicted phase-field obtained by exploiting the two strategies (POD and DI) with $N_c=30$ checkpoints. It is evident that performing displacement interpolation between HF data produces more accurate results than using the POD basis. This expected improvement is indeed related to the spatial differences in the support regions between snapshots: the POD projection essentially computes a linear combination of the available modes leading to a noticeable smearing of the interface, large oscillations in the union of the supports, and intermediate unphysical values for the reconstructed field. In contrast, displacement interpolation correctly predicts sharply varying values between zero and one in the whole spatial domain, detecting with high accuracy the evolution of the interface.

\begin{figure}[t!]
\centering
\includegraphics[width=0.85\textwidth]{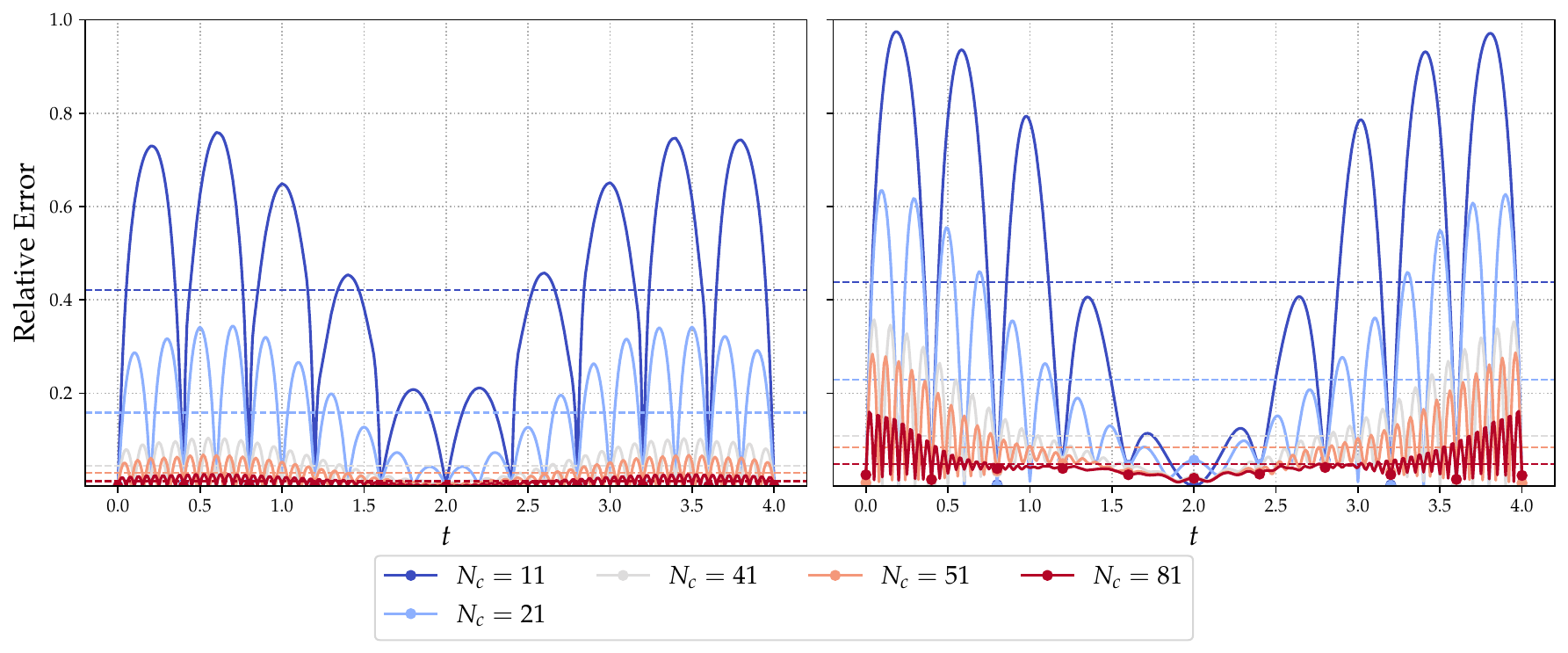}
\caption{Relative error behavior w.r.t.\ time for the Rider-Kother vortex at resolution $256^2$ with DI (left) and POD projection (right).}
\label{fig:err_POD_RK_time}
\end{figure}

\begin{figure}[t!]
\centering
\includegraphics[width=0.75\textwidth]{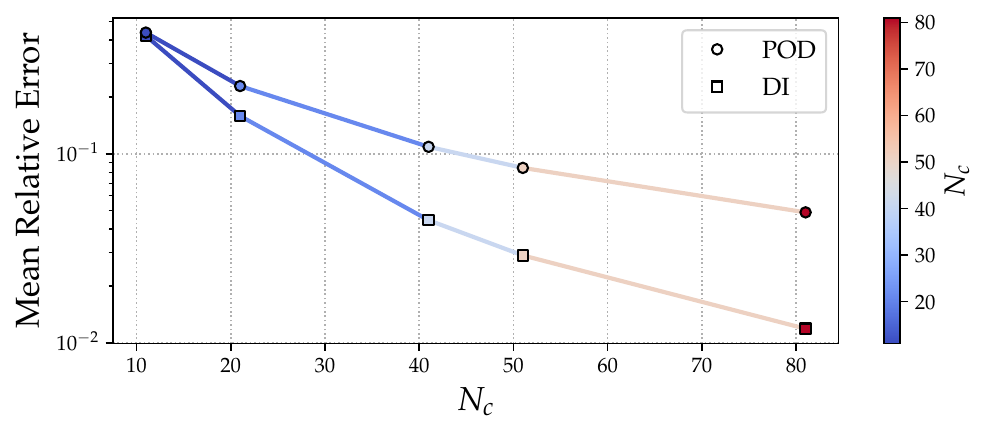}
\caption{Mean relative error w.r.t.\ the number of checkpoints $N_c$ for the Rider-Kother vortex at resolution $256^2$ with DI and POD projection.}
\label{fig:err_POD_RK_relative}
\end{figure}

In Figures~\ref{fig:err_POD_RK_time} and~\ref{fig:err_POD_RK_relative}, we present a more quantitative evaluation of the relative errors of displacement interpolation and POD projection varying the total number of checkpoints. In particular, for this analysis we select the number of checkpoints so that they are uniformly spaced across the available snapshots, enabling a fair comparison between DI, POD, and the FOM solutions. Since the dataset contains $401$ snapshots in total, this condition is satisfied for $N_{c} = 11, 21, 41, 51, 81$.
We notice that for a small number of checkpoints, the relative errors show only minor differences between the two approaches. However, as the number of checkpoints increases, displacement interpolation between consecutive data points provides significantly more accurate results compared to POD projection. In particular, note the clearly intermittent behavior at the initial and late stages of the dynamics (see Figure~\ref{fig:err_POD_RK_time}). In these regions, where the droplet moves at its highest speed, the POD projection struggles to accurately capture the dynamics of the problem, and provides seemingly better results only in reconstructing the states that have been used to construct the POD basis. From Figure~\ref{fig:err_POD_RK_relative}, depicting the analysis of the decay of the error w.r.t.\ the number of checkpoints $N_c$, it may seem that displacement interpolation yields only modest improvements compared to POD projection in terms of relative error, but we emphasize that this observation is largely influenced by the specific error metric adopted, which is biased towards pointwise discrepancy methodologies. Indeed, such a metric is partially misleading, as it does not fully capture the key features of interest in the solution, namely the interface morphology and the accuracy of the reconstructed phase-field values.


\begin{figure}[t!]
\centering
\subfigure{\includegraphics[width=0.45\textwidth]{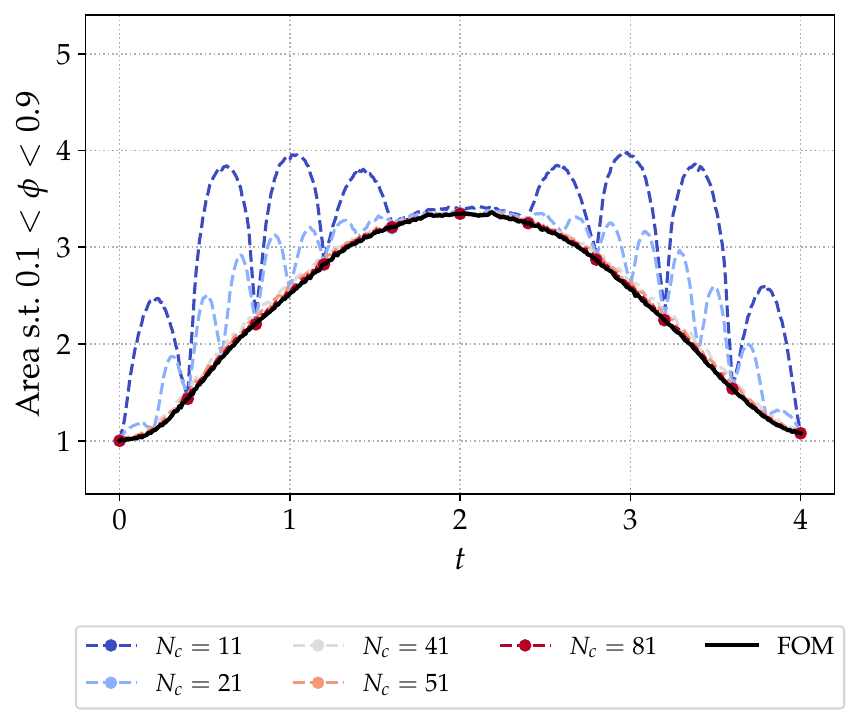}}
\subfigure{\includegraphics[width=0.45\textwidth]{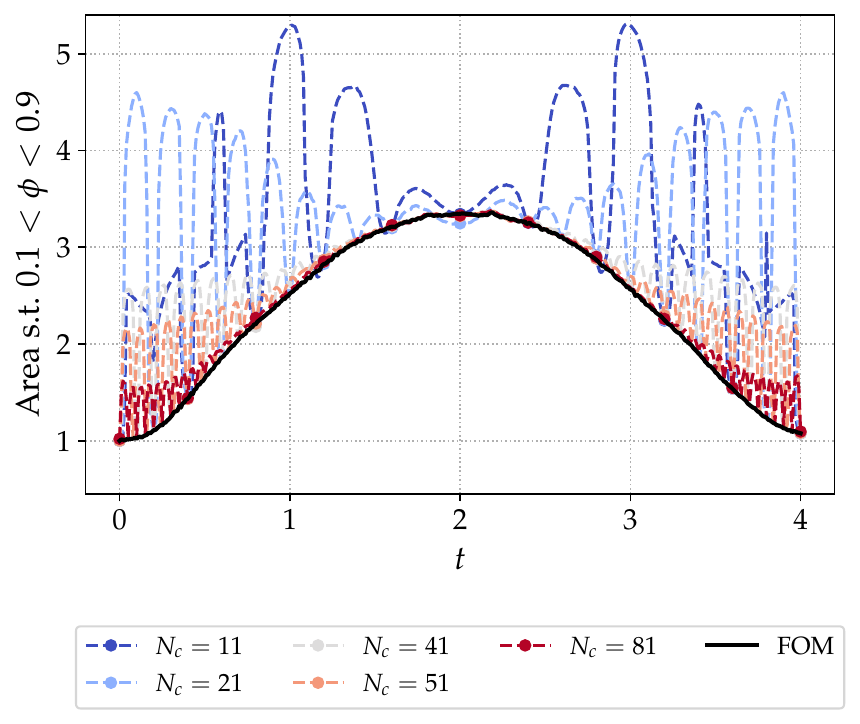}}
\caption{Normalized area of the spatial domain such that $0.1 < \phi < 0.9$ using DI (left) and POD projection (right) for different $N_c$ values.}
\label{fig:surf}
\end{figure}

\begin{figure}[t!]
\centering
\subfigure{\includegraphics[width=0.45\textwidth]{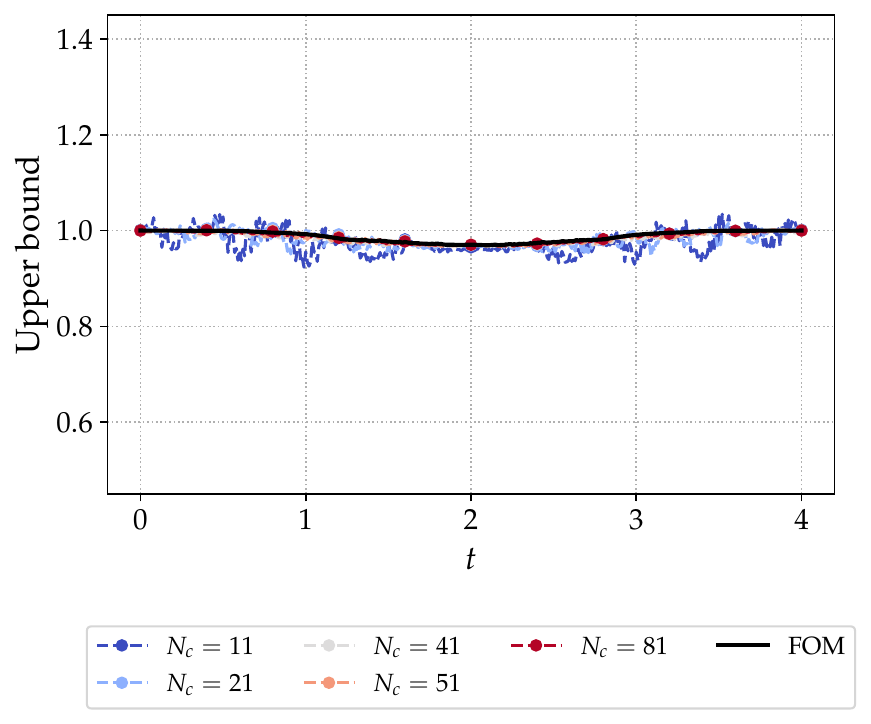}}
\subfigure{\includegraphics[width=0.45\textwidth]{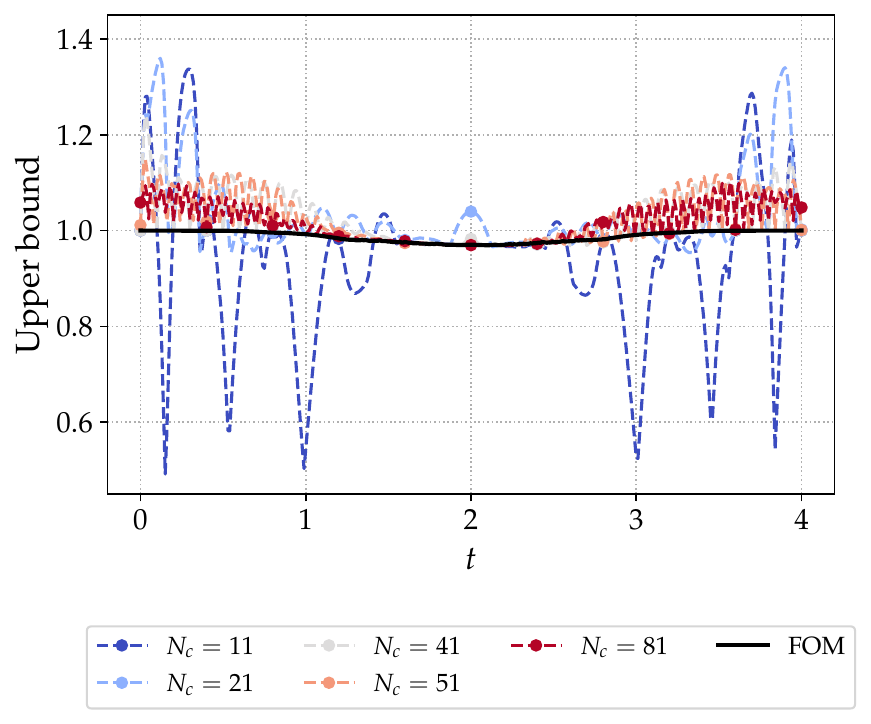}}
\caption{Maximum value of the predicted phase-field as a function of time using DI (left) and POD projection (right) for different $N_c$ values.}
\label{fig:bounds}
\end{figure}

For these reasons, we also analyze the problem from a more modeling-oriented perspective by examining the total area where $0.1 < \phi < 0.9$ that, in the case of bounded values, provides an approximation of the droplet's interfacial area.
Specifically, Figure~\ref{fig:surf} shows that the general trend of the FOM behaves as expected for this problem: the droplet is initially circular, then stretches into a thin ligament, maximizing the total surface area, and finally comes back to its initial shape due to the periodicity of the advection field. Displacement interpolation captures this behavior accurately, whereas, in contrast, POD projection deviates significantly from the expected trend, particularly at the initial and final stages of the dynamics, where the prescribed advection velocity reaches its maximum magnitude.
Furthermore, even in regions where common support exists, the POD projection essentially acts as a high-order interpolation, and even increasing the number of basis spanning the POD space, the benefits in accuracy are counteracted by the introduction of spurious oscillations near sharp features. This behavior is clearly visible in Figure~\ref{fig:bounds}, where we show a one-to-one comparison of the upper bounds of the phase-field predicted by displacement interpolation and POD projection. It is evident that DI consistently remains very close to the FOM simulations, whereas the POD projection exhibits strong oscillations around the unitary upper bound.
%
\subsubsection{MF-OT-ROM}
Having highlighted how POD-based ROMs perform poorly when applied to two-phase flows, especially when evaluated on more case-specific metrics such as boundedness and interfacial area, we are ready to present the performance of the proposed extensions of OT-based surrogate modeling. We begin with the MF-OT-ROM approach introduced previously, where displacement interpolation is used to interpolate the residuals between high- and low-fidelity solutions. The interpolated residual is then employed to \emph{correct} the LF solutions in order to recover, as accurately as possible, an approximation of the HF solutions.  
For this test case, we consider three levels of fidelity: $128^2$, $256^2$ and $512^2$, computing the pairwise residuals between the different grid resolutions.

\begin{figure}[t!]
\centering
\includegraphics[width=0.75\textwidth]{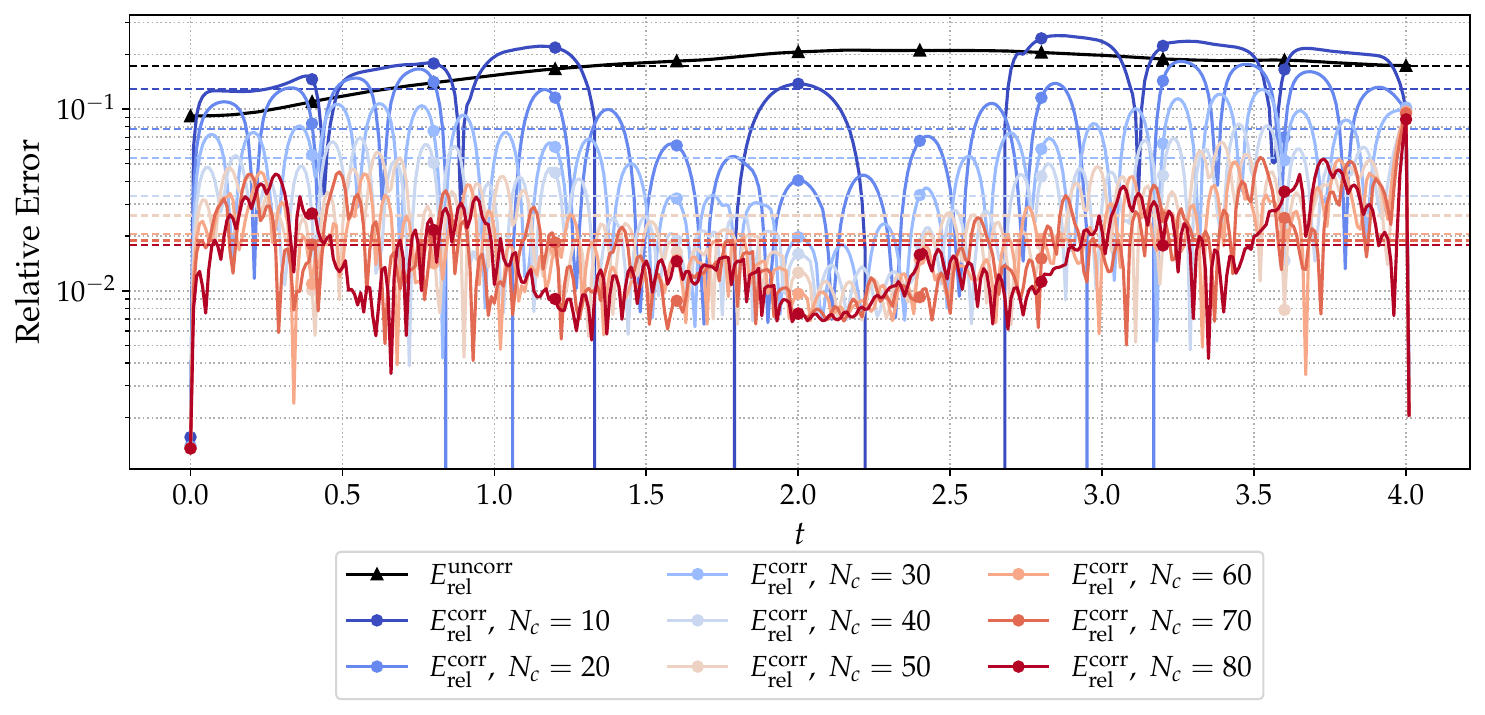}
\caption{Relative error between low- and high-fidelity solutions ($256^2$ and $512^2$) w.r.t.\ different number of checkpoints for Rider-Kothe vortex with $R=0.11$. Horizontal dashed lines represent the mean error for each temporal trajectory.}
\label{fig:RK_err}
\end{figure}

\begin{figure}[t!]
\centering
\includegraphics[width=0.75\textwidth]{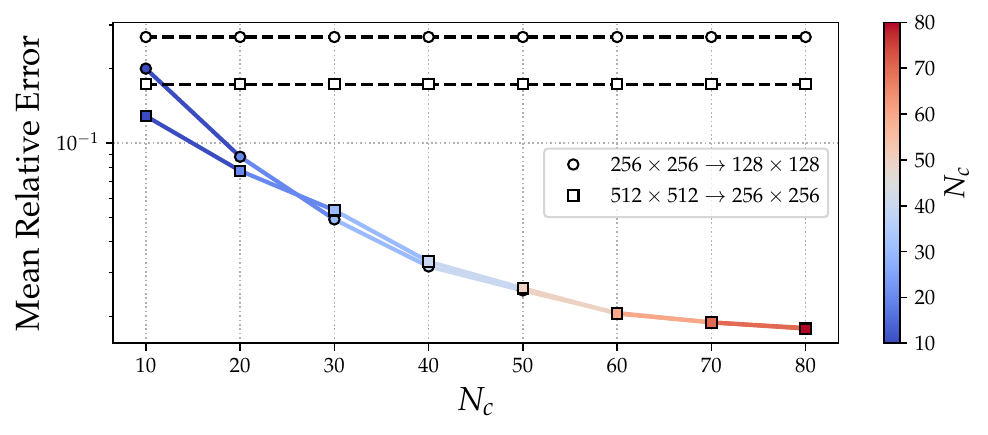}
\caption{Mean relative errors against number of checkpoints.}
\label{fig:RK_mean}
\end{figure}

As in the previous section, to perform a complete accuracy study of the strategy with respect to the available data, we vary the number of checkpoints used for displacement interpolation. 
Figures~\ref{fig:RK_err} and~\ref{fig:RK_mean} quantify the relative errors produced by the proposed approach, measured as the difference between the LF solution (corrected or not via OT) and the HF solution. In Figure~\ref{fig:RK_err}, we show the relative error between low- and high-fidelity solutions ($256^2$ and $512^2$) as a function of time, varying the number of checkpoints. The uncorrected LF solution is shown in black, while a color gradient from blue to red indicates an increasing number of checkpoints. With relatively few checkpoints (e.g., $N_c = 10$), the corrected LF solution can locally increase the error relative to the HF solution, although the mean error decreases (see Figure~\ref{fig:RK_mean}). As the number of checkpoints increases, the relative error decreases throughout the whole trajectory, demonstrating progressively better agreement with the HF solution.
In Figure~\ref{fig:RK_mean}, we consider both settings with $128^2$ and $256^2$ exploited as LF, and we investigate the decay of the mean error. Notice that the trend is almost identical between the two cases, but always over-performing the uncorrected approximations, suggesting that the correction is effective for both resolutions and that in both cases we are accurately capturing all the features from the dynamics. Indeed, by reporting the relative errors, the improvement obtained by correcting the $128^2$ solution to match the $256^2$ solution is, in relative terms, almost the same as that obtained by correcting the $256^2$ solution to match the $512^2$ solution.

\begin{figure}[t!]
\centering
\includegraphics[width=0.75\textwidth]{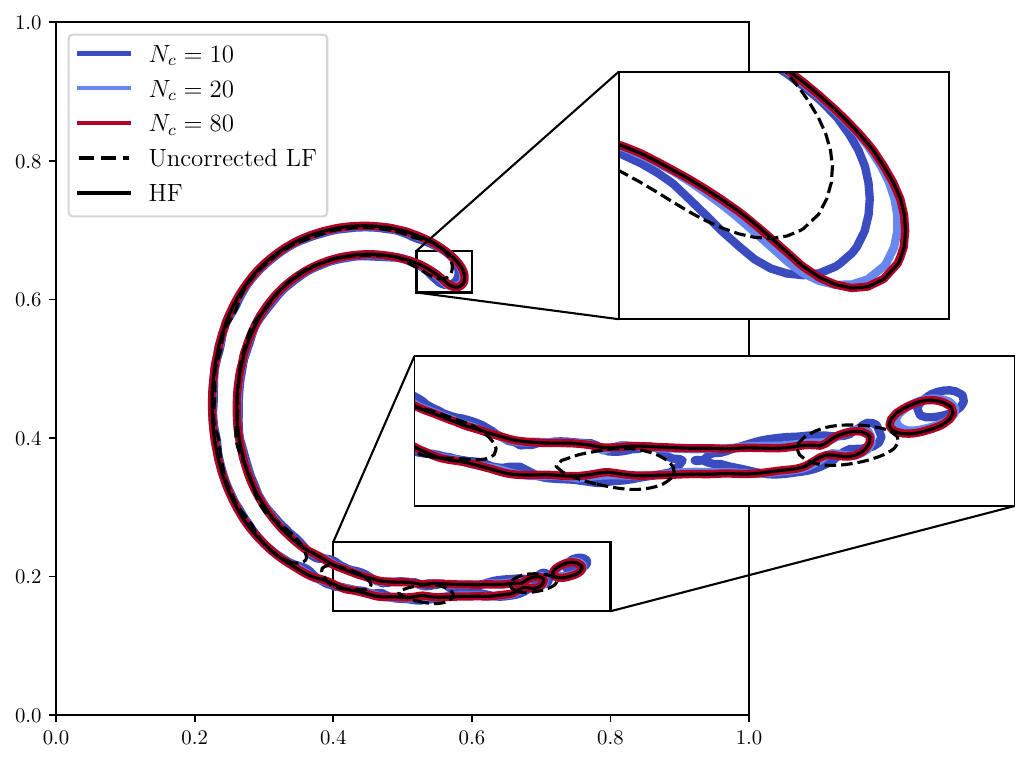}
\caption{Isocontour of $\phi=0.5$ for $256^2$ HF (dash-dotted black), $128^2$ LF (dashed black) and corrected LFs (color gradient) at $t=1.95$.}
\label{fig:iso_MF1}
\end{figure}

\begin{figure}[t!]
\centering
\includegraphics[width=0.75\textwidth]{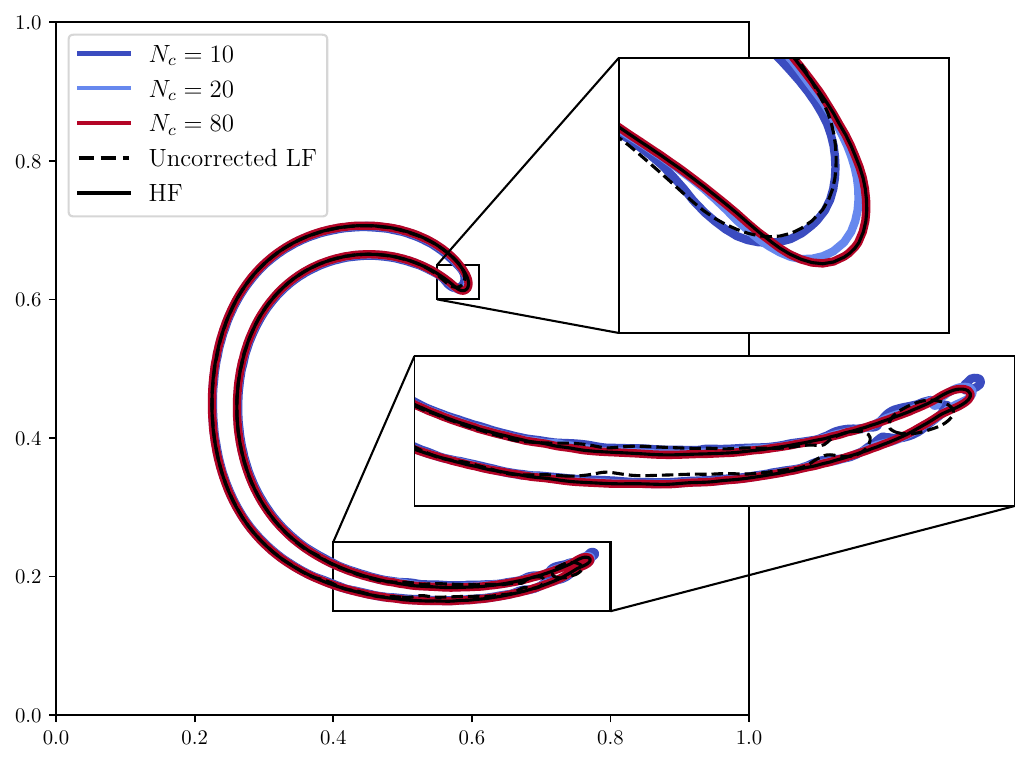}
\caption{Isocontour of $\phi=0.5$ for $512^2$ HF (dash-dotted black), $256^2$ LF (dashed black) and corrected LFs (color gradient) at $t=1.95$.}
\label{fig:iso_MF2}
\end{figure}

From a qualitative perspective, Figures~\ref{fig:iso_MF1} and~\ref{fig:iso_MF2} show the isolines $\phi=0.5$, representing the interface location, predicted by the corrected LF solutions, along with the uncorrected LF (dashed black) and HF (dash-dotted black) references. In Figure~\ref{fig:iso_MF1}, the $128^2$ solution is considered as LF, whereas in Figure~\ref{fig:iso_MF2}, the $256^2$ solution serves as the LF reference, and increasing the number of checkpoints significantly improves in both cases the predicted interface location, correctly recovering tail and bubble structures. This analysis confirms the importance of exploiting refined meshes for the simulation of two-phase flows, and thus further motivates the need for multi-fidelity corrections to obtain efficient simulations.

\begin{figure}[t!]
\centering
\includegraphics[width=0.95\textwidth]{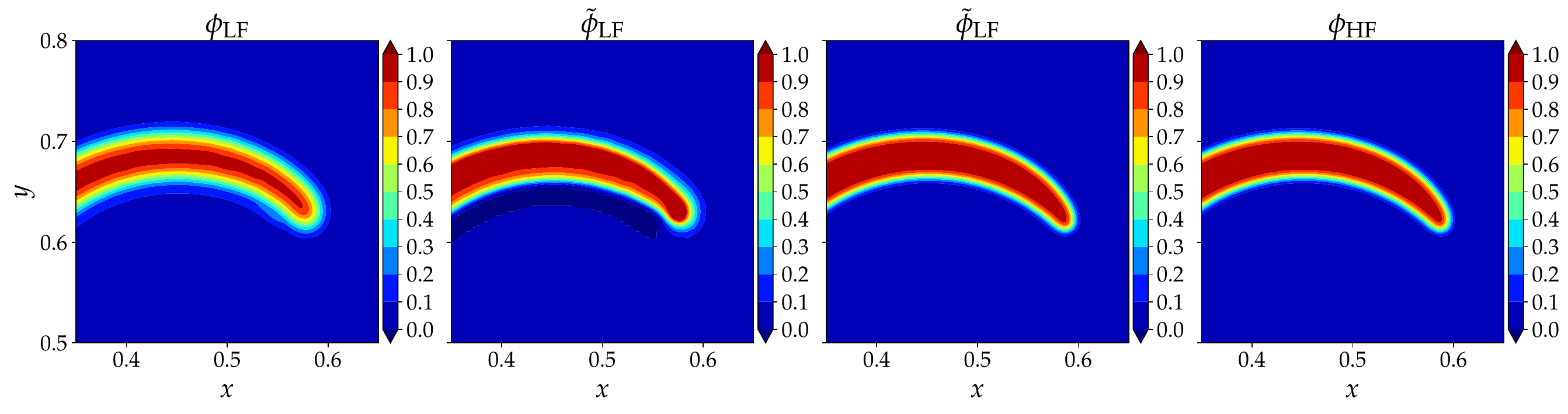}
\caption{Phase-field isocontours at $t=1.95$. From left to right: uncorrected LF ($256^2$), corrected LF with $N_{c}=10$, corrected LF with $N_{c}=80$, and reference HF ($512^2$).}
\label{fig:fullfield}
\end{figure}

In Figure~\ref{fig:fullfield}, we show the predicted field for $\phi$ near the top of the droplet at $t=1.95$. From left to right, the amount of OT correction is increased: the left panel shows the uncorrected LF solution, the two center panels show the corrected LF solutions using $N_c = 10$ and $N_c = 80$ checkpoints, respectively, while the right panel show the HF solution. The nice performance of the methodology is already evident exploiting only $10$ checkpoints, with the solution showing a significant improvement in terms of interface thickness, while using $N_c = 80$ produces a corrected solution that is nearly indistinguishable from the HF simulation.
%
\subsubsection{PMF-OT-ROM}
In this section, we increase the complexity of the Rider-Kothe vortex by considering its parametric extension. Indeed, HF solutions may not be directly available for all parametric trajectories, making the problem more challenging but also more relevant for applications where new parametric dynamics must be computed efficiently. In this context, we propose to apply displacement interpolation to HF data from neighboring trajectories in order to generate a new synthetic HF trajectory for previously unseen parameter values (see the discussion in Section \ref{sec:pmf_ot_rom} concerning the computational cost). Such synthetic checkpoints can then be used to compute the residual and correct the corresponding LF solution. 

\begin{figure}[t!]
\centering
\includegraphics[width=0.75\textwidth]{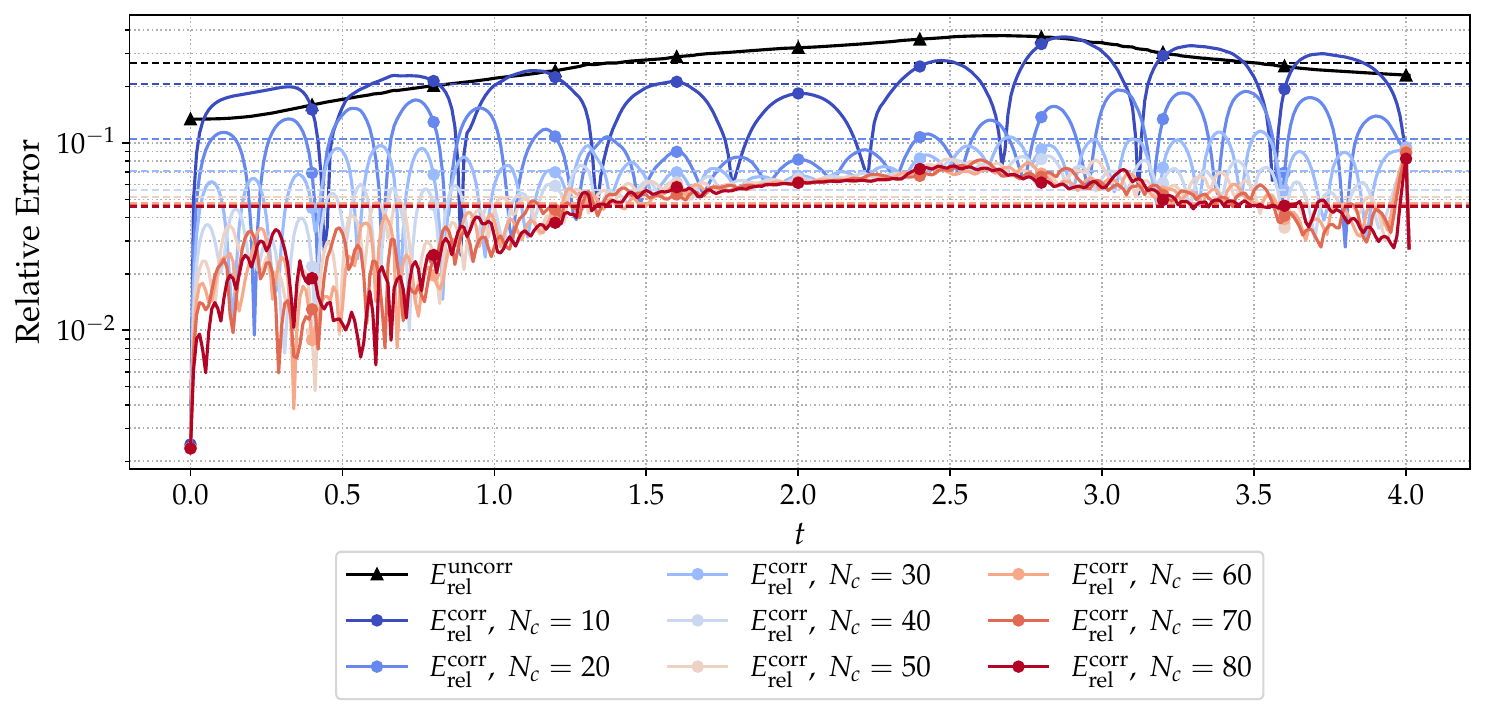}
\caption{Relative error between low- and high-fidelity ($128^2$ and $256^2$) w.r.t.\ different number of checkpoints for the Rider-Kothe vortex with unseen parameter $R=0.11$. Horizontal dashed lines represent the mean error for each temporal trajectory.}
\label{fig:RK_err_par}
\end{figure}

\begin{figure}[t!]
\centering
\includegraphics[width=0.75\textwidth]{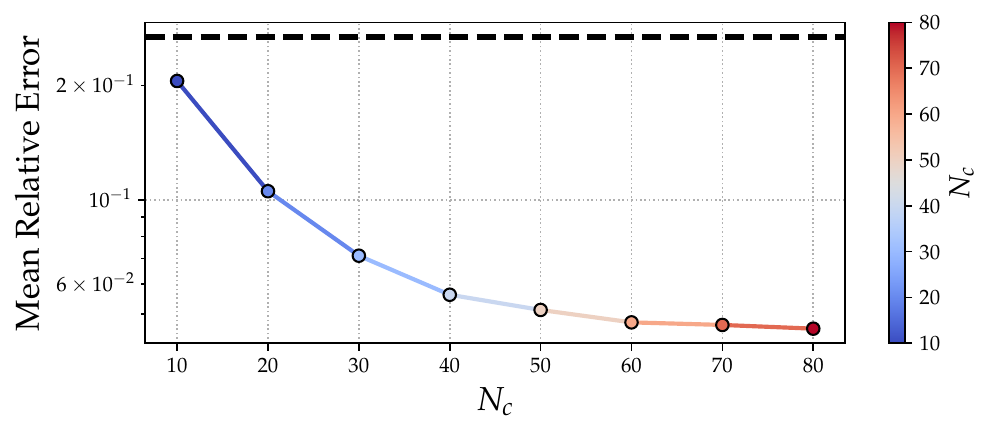}
\caption{Mean relative errors against number of checkpoints for the unseen parameter $R=0.11$ with LF $128^2$ and HF $256^2$.}
\label{fig:RK_mean_par}
\end{figure}

In Figure~\ref{fig:RK_err_par} we report the relative error between the corrected LF solution and the validation HF solution for an unseen parametric trajectory. We can observe that increasing the number of checkpoints significantly decreases the errors throughout the whole dynamics. In Figure~\ref{fig:RK_mean_par} we show the mean error over the entire trajectory, which further highlights the decrease in error as the number of checkpoints increases.

\begin{figure}[t!]
\centering
\subfigure[Upper bound]{\includegraphics[width=0.45\textwidth]{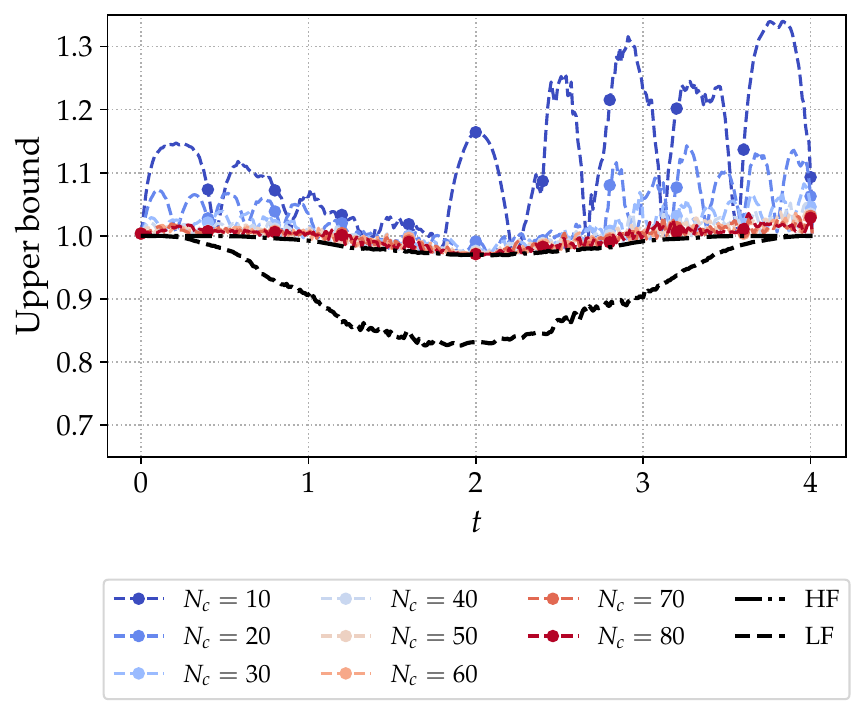}}
\subfigure[Interfacial area]{\includegraphics[width=0.45\textwidth]{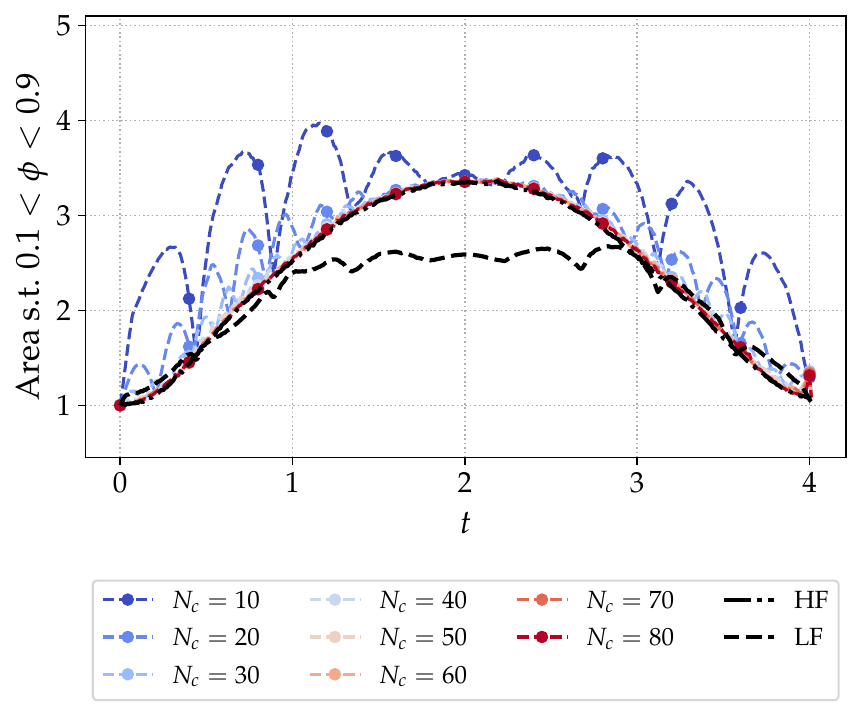}}
\caption{Maximum value of the predicted phase-field as a function of time (left) and interfacial area (right). Dash-dotted line, HF $256^2$; dashed, uncorrected LF $128^2$; color gradient, corrected LFs.}
\label{fig:bounds_inter_par}
\end{figure}

Finally, in Figure~\ref{fig:bounds_inter_par}, as a metric which is more specific for this problem, we investigate the upper bound and total interfacial area by comparing HF, uncorrected LF and a series of corrected LFs with the proposed approach with different number of checkpoints. Notice that as the number of checkpoints increases we obtain an almost perfect match in both quantities.
%
\subsection{Rayleigh-Taylor instability}
We now consider an even more complex benchmark represented by the Rayleigh-Taylor instability, which occurs when an interface between two fluids with different densities experiences a pressure gradient opposing the density gradient~\cite{tryggvason1988numerical}. The computational domain is $[d \times 4d]$ with $d=1$, and the interface is initially defined as the curve $y(x) = 2d + 0.1 d \cos{(2\pi x/d)}$. The Rayleigh-Taylor instability is characterized by the Reynolds number $\mathrm{Re} = (\rho_{1} d^{3/2} \|\textbf{g}\|^{1/2})/\mu$ and the Atwood number $\mathrm{At} = (\rho_{1} - \rho_{2})/(\rho_{1} + \rho_{2})$ which are classically set to $3000$ and $0.5$. Time is adimensionalized as $t^{*} = t /\sqrt{d/(\|\textbf{g}\| \mathrm{At})}$, while for the EOS relation in Equation \eqref{eq:EOS} we choose the following parameters $\gamma_{1}=\gamma_{2}=1.4$, $P^{\infty}_{1}=P^{\infty}_{2}=6000$. 

The Reynolds number characterizes the relative importance of inertial to viscous effects in the coupled system and is governed by gravity and the mixture viscosity, which is assumed constant across the two phases for this case. The Atwood number, on the other hand, quantifies the density ratio between the phases. When $\mathrm{At}=0$, the two phases have identical densities, whereas increasing values of the Atwood number correspond to progressively larger density ratios. As the Atwood number increases, one therefore expects the flow to develop increasingly complex interfacial patterns.

The top boundary is treated as a Riemann-invariant boundary condition with zero velocity and constant pressure, the bottom boundary is a no-slip wall whereas slip-wall boundary conditions are prescribed on the other lateral sides of the domain. For this study, two different meshes with quadrilateral elements are considered, involving $40$k and $160$k degrees of freedom, LF and HF, respectively. A $4^{\mathrm{th}}$ order Runge-Kutta scheme was used for time integration.  For this second test case, we consider a dataset consisting of $301$ snapshots, corresponding to a time window of $\Delta \tau = 0.01$. The time step sizes used in the simulations are approximately $\Delta t_{40\mathrm{k}} \approx 10^{-5}$ and $\Delta t_{160\mathrm{k}} \approx 5 \cdot 10^{-6}$ for the two spatial resolutions considered, to satisfy CFL constraints.

\begin{figure}[t!]
\centering
\subfigure[$\mathrm{At}=0.25$]{%
\includegraphics[width=0.18\textwidth]{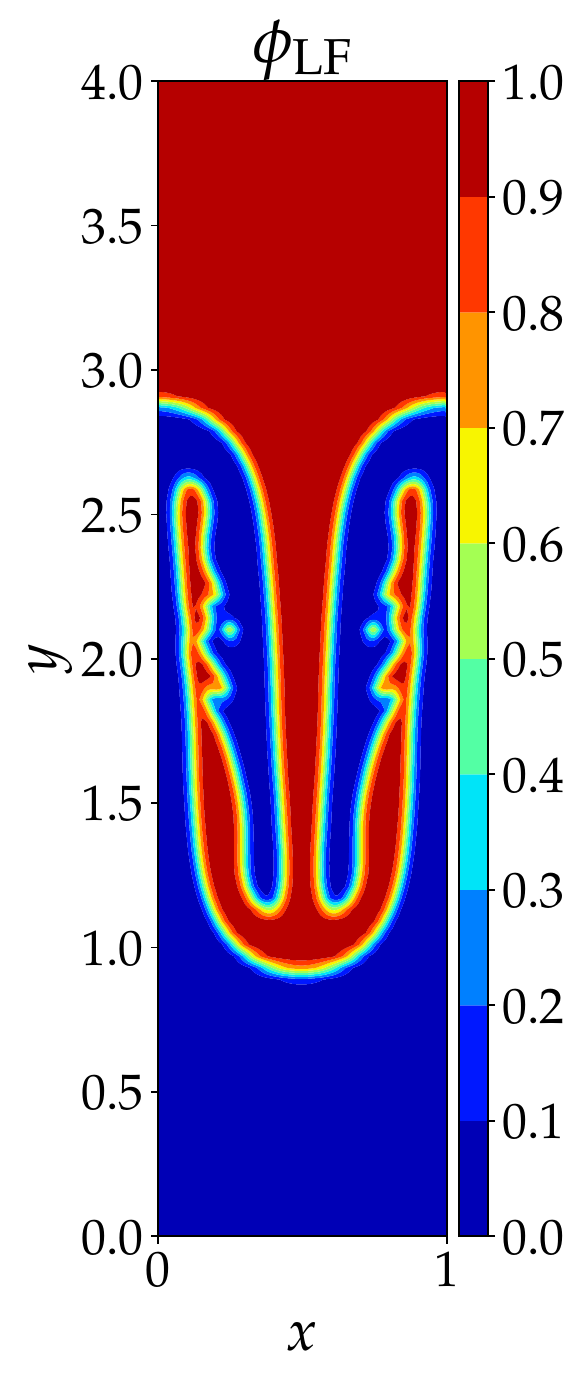}}%
\hspace{0.01\textwidth}
\subfigure[$\mathrm{At}=0.5$]{%
\includegraphics[width=0.18\textwidth]{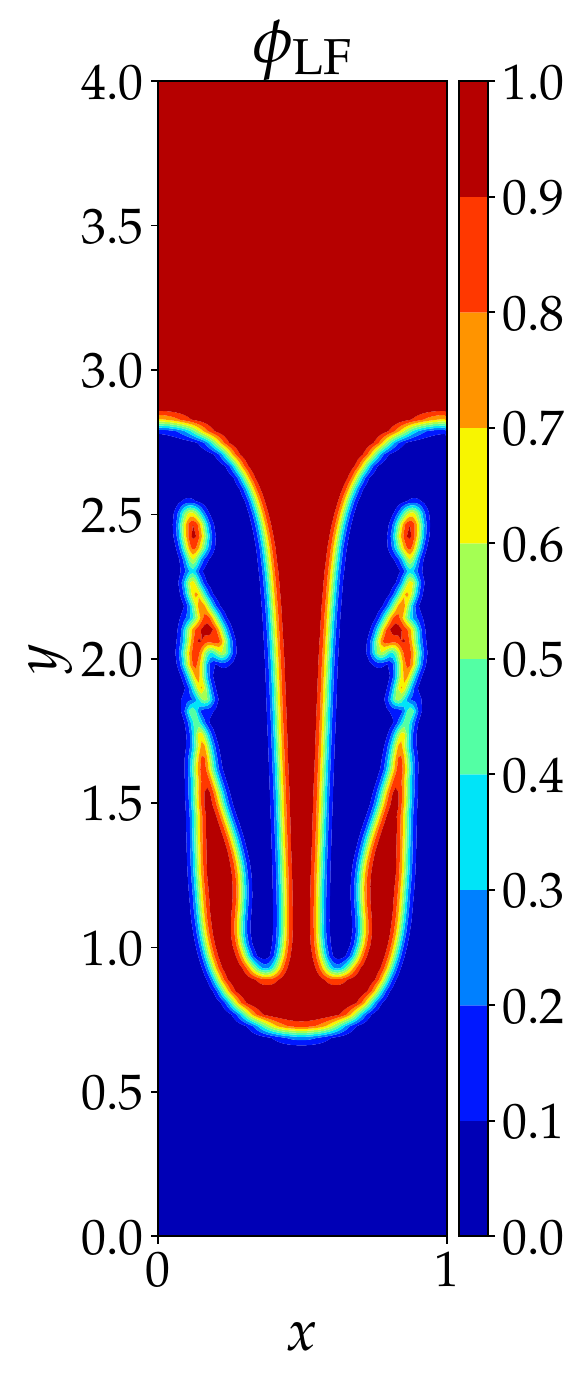}}%
\hspace{0.01\textwidth}
\subfigure[$\mathrm{At}=0.75$]{%
\includegraphics[width=0.18\textwidth]{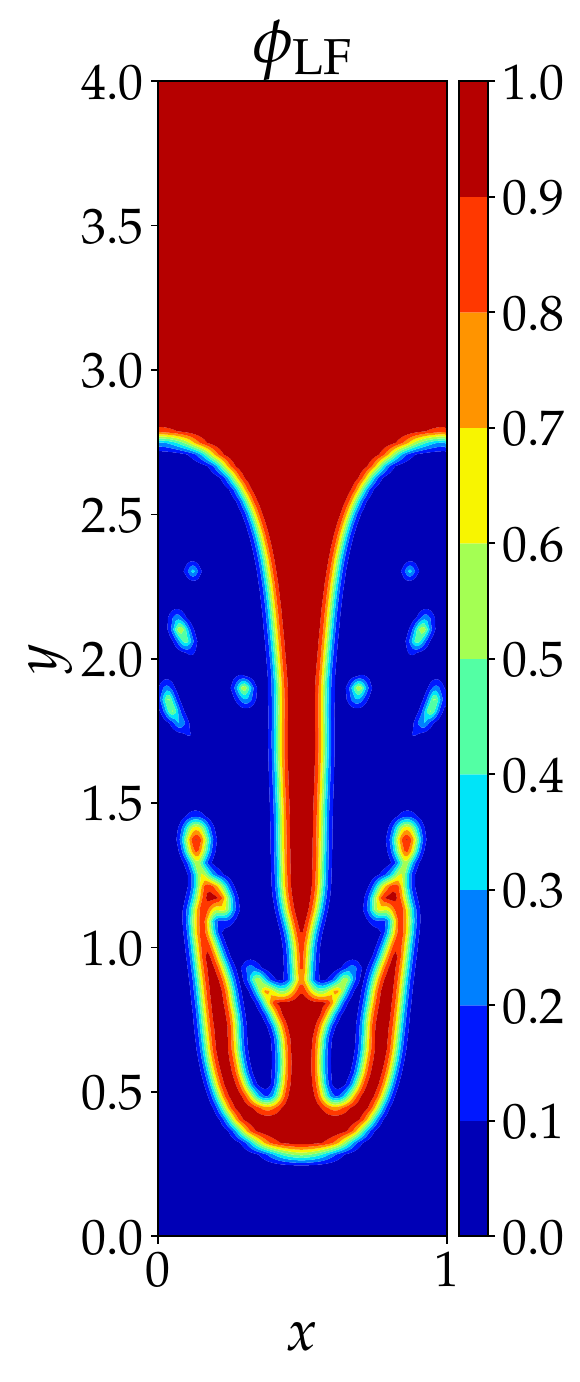}}

\vspace{0.3em}

\subfigure[$\mathrm{At}=0.25$]{%
\includegraphics[width=0.18\textwidth]{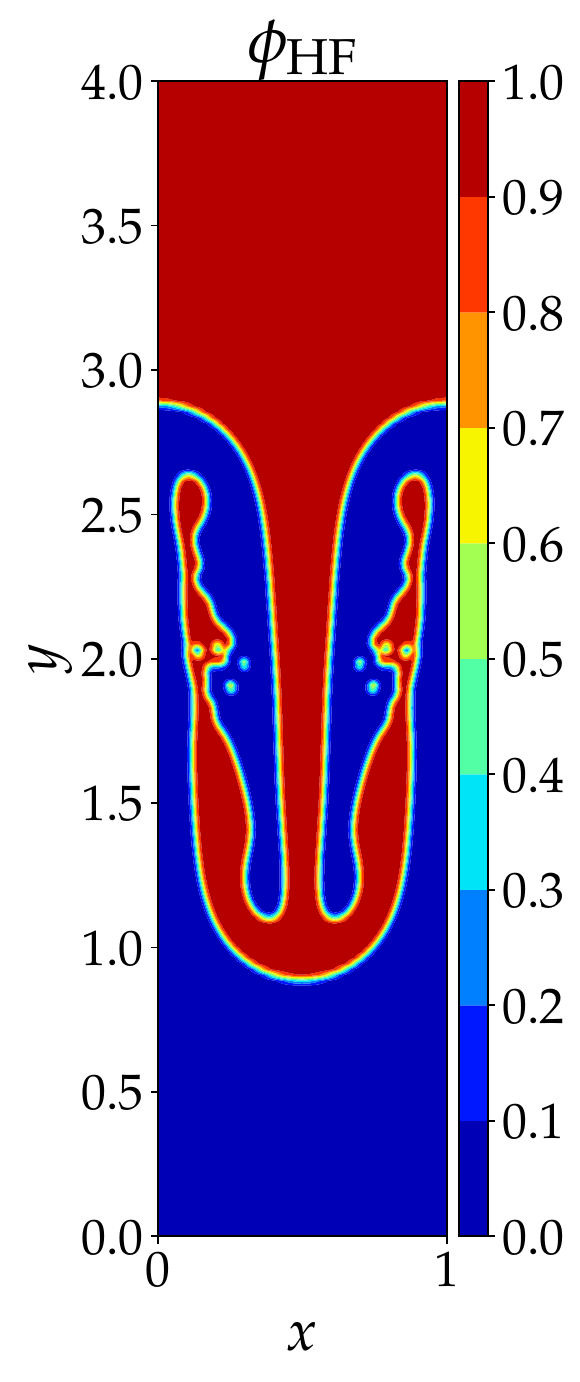}}%
\hspace{0.01\textwidth}
\subfigure[$\mathrm{At}=0.5$]{%
\includegraphics[width=0.18\textwidth]{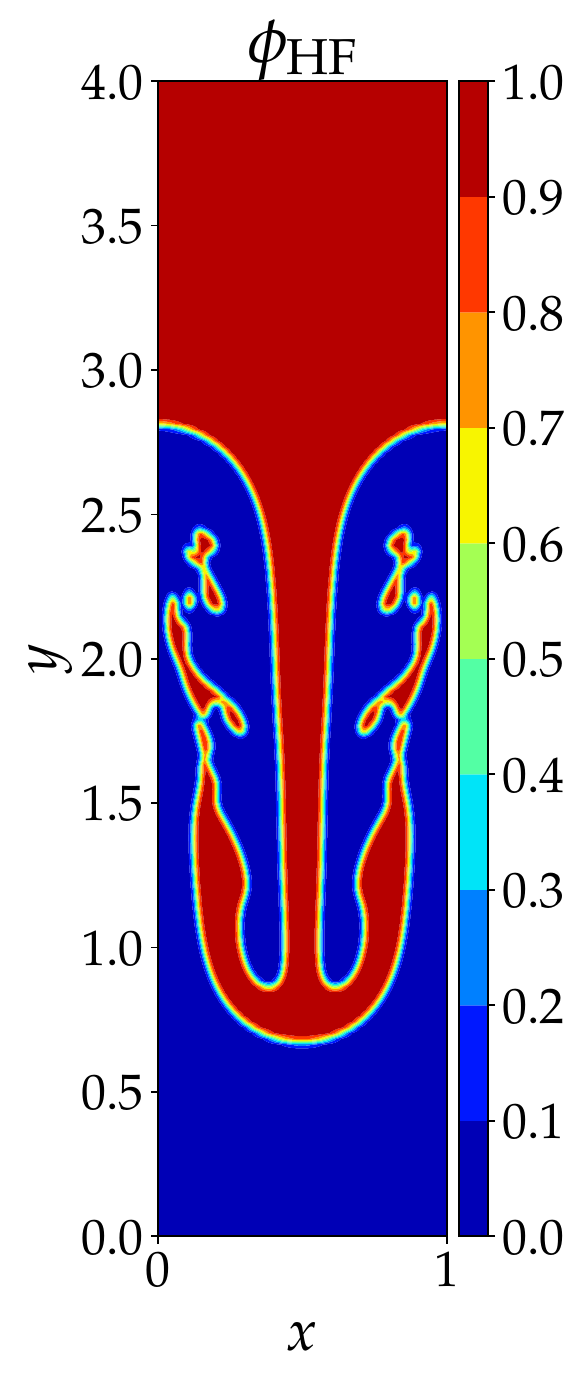}}%
\hspace{0.01\textwidth}
\subfigure[$\mathrm{At}=0.75$]{%
\includegraphics[width=0.18\textwidth]{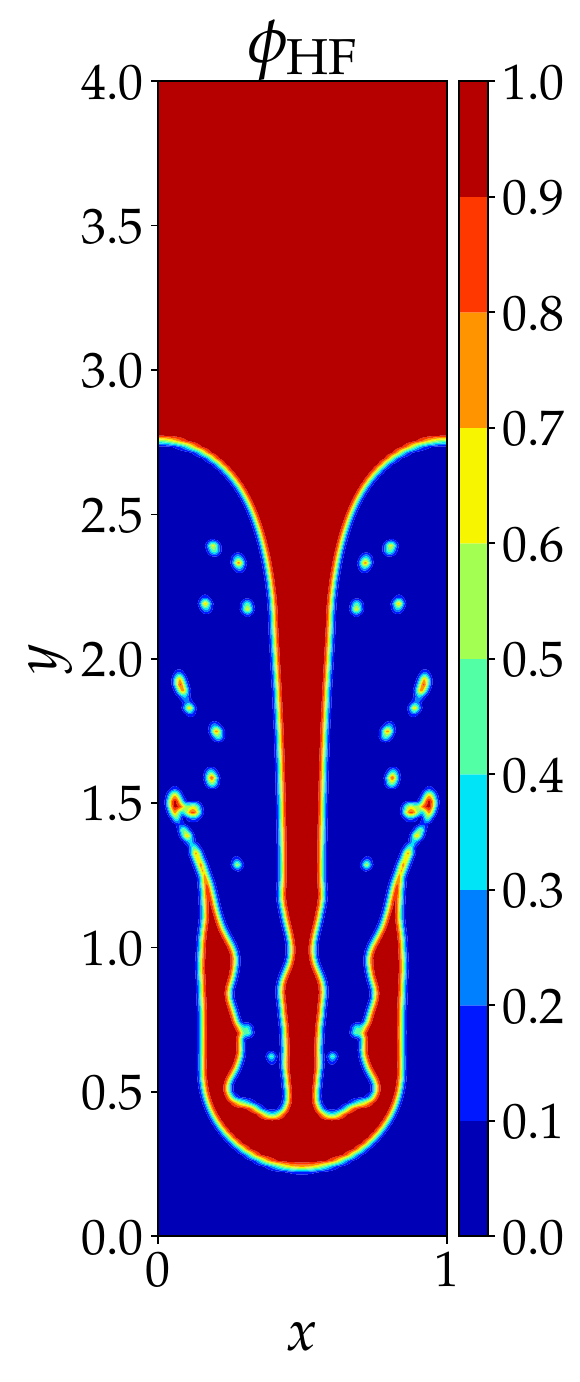}}
\caption{LF (top) and HF (bottom) phase-field isocontours at different Atwood numbers for the Rayleigh-Taylor instability at $t^{*}=3.0$.}
\label{fig:RT_overview}
\end{figure}

Concerning its parametric extension, we consider a variation of Atwood numbers in the following range:
\begin{equation}
    \mathrm{At}_{k}=0.25 + k \Delta \mathrm{At} \quad \mathrm{with} \quad \Delta \mathrm{At} = 0.025 \quad \mathrm{and} \quad k=0,\dots,20.
\end{equation}
An overview of the parametric problem's dynamics is shown in Figure~\ref{fig:RT_overview}. If not specified otherwise, we will always consider the $\mathrm{At}=0.35$ case to test our framework for the Rayleigh-Taylor instability.
\subsubsection{MF-OT-ROM}
\begin{figure}[t!]
\centering
\includegraphics[width=0.75\textwidth]{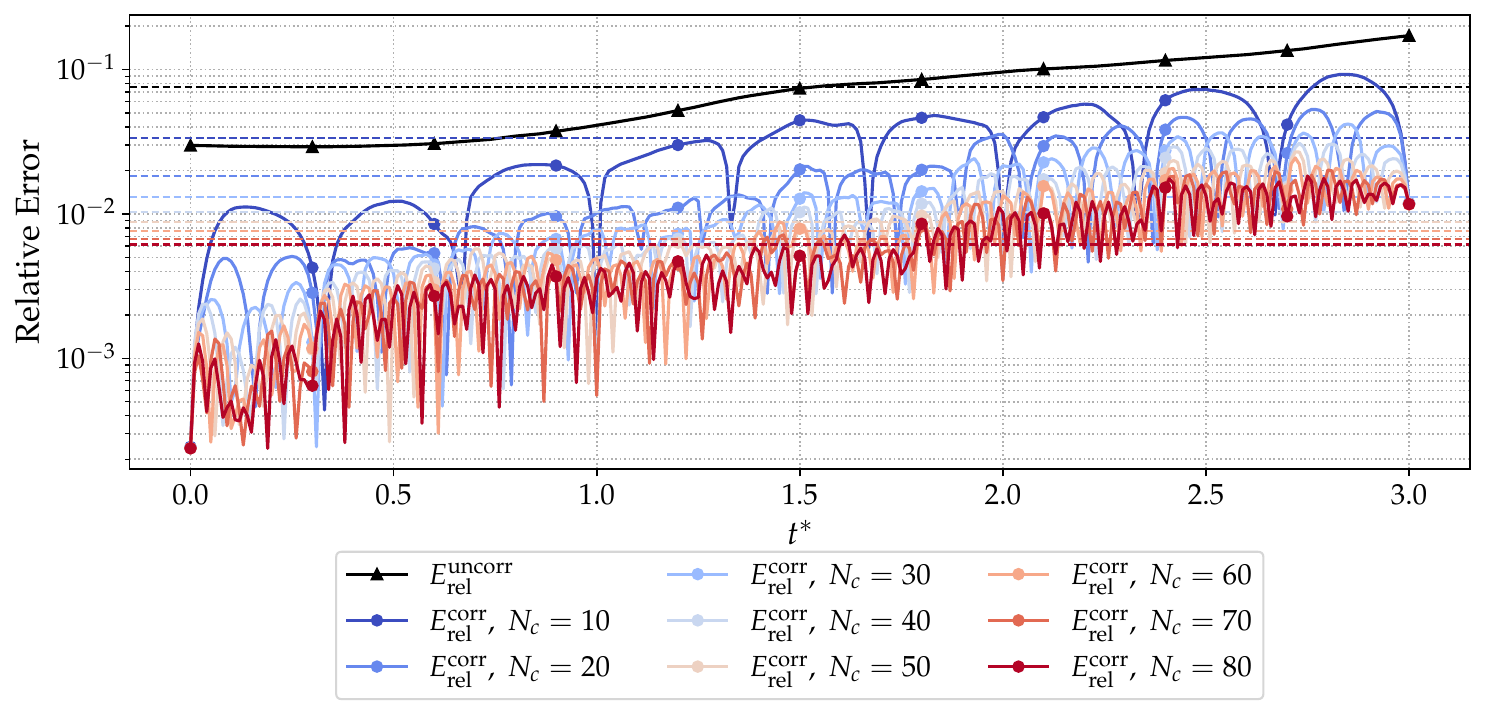}
\caption{Relative error between low- and high-fidelity for different number of checkpoints for Rayleigh-Taylor instability with $\mathrm{At}=0.35$. Horizontal dashed lines represent the mean error for each temporal trajectory.}
\label{fig:RT_err}
\end{figure}

Similarly to the previous test case, we aim to quantify the relative errors between the low- and high-fidelity solutions, both with and without DI correction. In Figure~\ref{fig:RT_err}, we report the relative error as a function of time for different numbers of checkpoints. 

A markedly different trend emerges compared to the previous test case: here, the dynamics \emph{accelerates} over time. Initially, the interface evolves relatively slowly, but it subsequently accelerates due to gravitational effects, leading to complex breakup patterns. This behavior is also reflected in the MF-OT-ROM procedure, as the errors tend to grow as time progresses. This indicates that a more involved regression of the virtual time $\alpha$ may be needed to better capture the accelerating dynamics. Despite this, the current approximation already provides a significant improvement over the uncorrected LF solution, even with a limited amount of data, and increasing the number of checkpoints significantly improves the overall accuracy of the method, consistently with the behavior observed in the previous test case.

\begin{figure}[t!]
\centering
\includegraphics[width=0.65\textwidth]{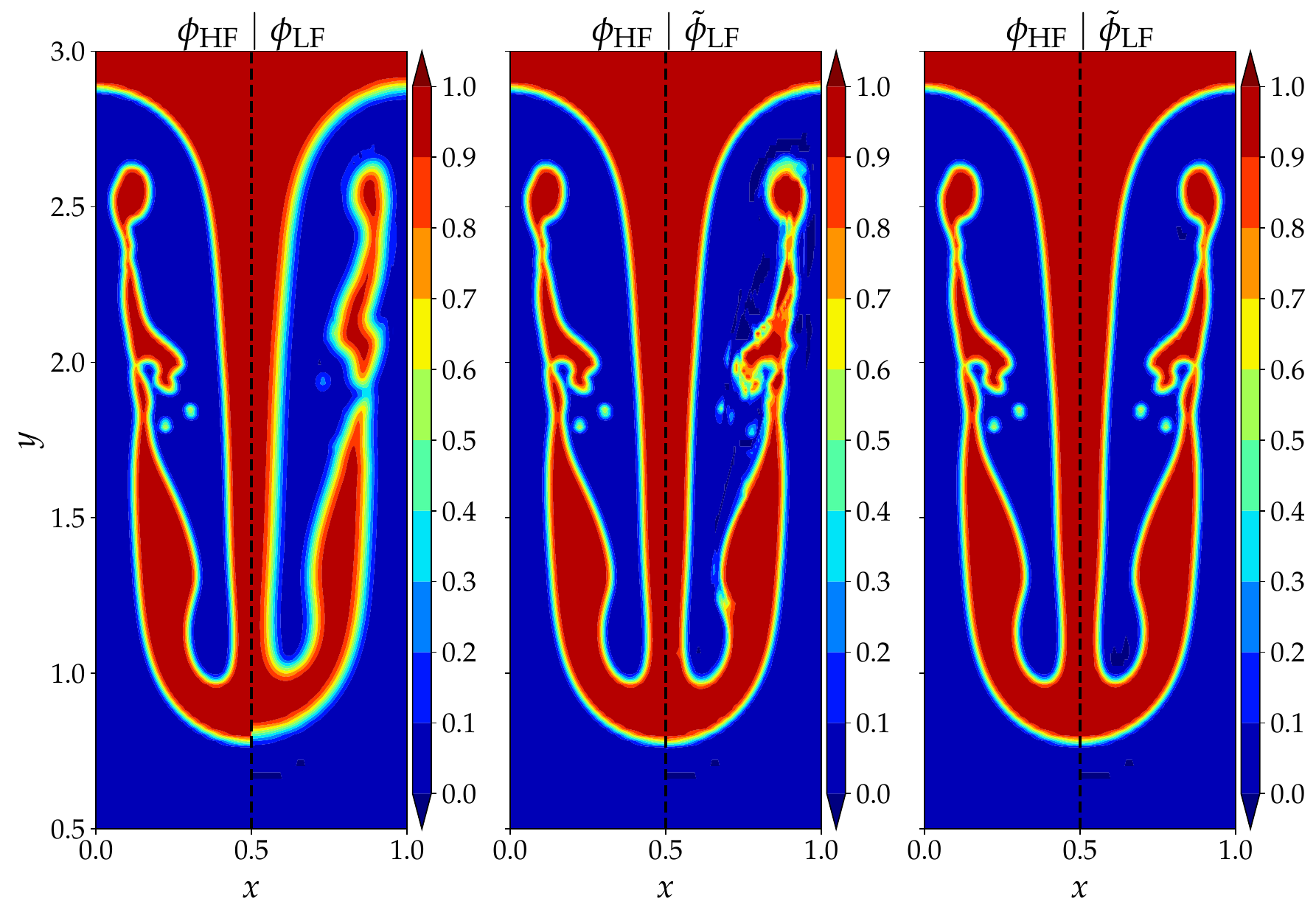}
\caption{Phase-field isocontours for Rayleigh-Taylor instability at $t^{*}=2.95$ and $\mathrm{At}=0.35$. In each subfigure, the left field is given by the HF solution, while right field shows respectively: LF solution, corrected LF solution with $N_{c}=10$ and corrected LF solution with $N_{c}=80$.}
\label{fig:RT_fullfield}
\end{figure}

In Figure~\ref{fig:RT_fullfield}, we present the predicted field of $\phi$ for the Rayleigh-Taylor instability problem. Each panel is split into two halves, where the HF solution is mirrored against the corresponding LF prediction (either corrected or uncorrected) to facilitate visual comparison. 
From left to right, the level of OT correction increases. The left panel displays the uncorrected LF solution, while the central and right panels show the corrected LF solutions obtained with $N_c = 10$ and $N_c = 80$ checkpoints, respectively.
We observe that even for the smallest value of $N_{c}$ considered here, the OT-based correction substantially improves the overall prediction. Increasing the number of checkpoints to $N_{c}=80$ further enhances the accuracy, to the point that the corrected LF solution becomes almost indistinguishable from the HF one.

\begin{figure}[t!]
\centering
\includegraphics[width=0.75\textwidth]{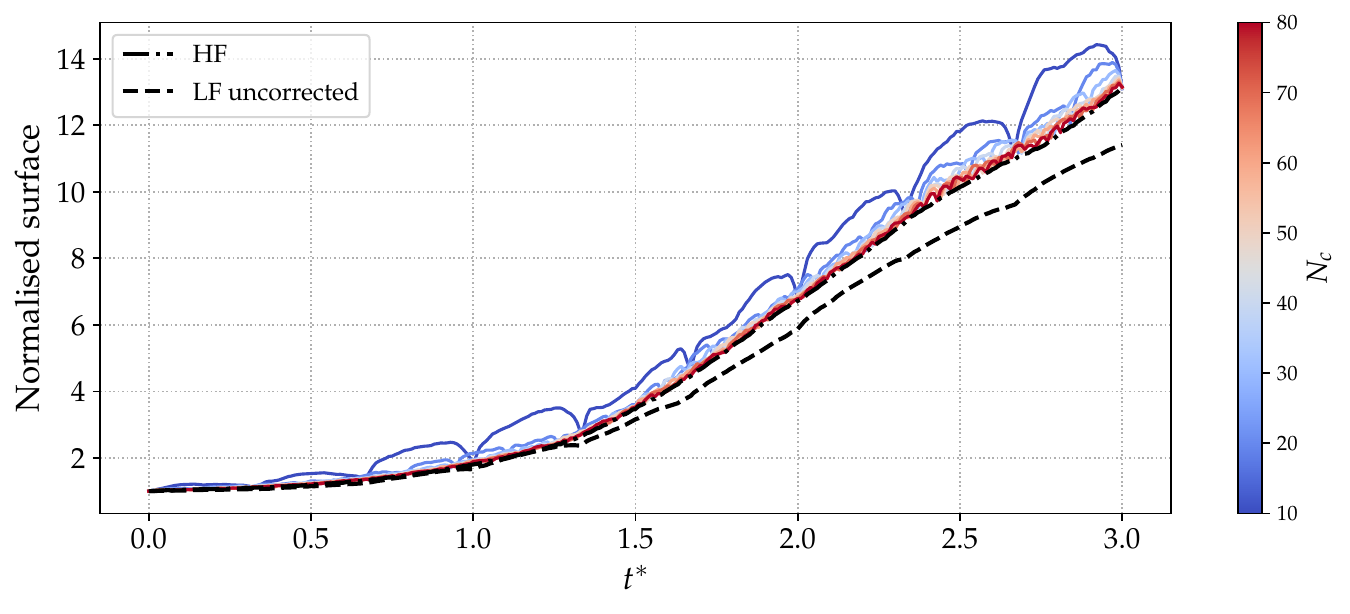}
\caption{Normalized interfacial surface as a function of time for the Rayleigh-Taylor instability problem. Dash-dotted line, HF; dashed, uncorrected LF; color gradient, corrected LFs.
}
\label{fig:RT_surf}
\end{figure}

More importantly, for this specific problem, where breakup phenomena are widespread, we assess the capability of the MF-OT-ROM approach to accurately predict the total interfacial surface from a quantitative perspective. In Figure~\ref{fig:RT_surf}, we report the normalized interface surface as a function of time for both fidelities, considering the uncorrected and OT-corrected LF solutions.
First, the temporal evolution clearly reflects the breakup dynamics. At early times, the interface between the two phases is relatively smooth, and the total surface is limited. As the instability develops, gravitational effects enhance mixing, generating complex structures and significantly increasing the contact area between the phases.
Second, a classical resolution effect can be observed: coarse simulations are generally unable to capture small-scale features, such as vortices and thin interfacial ligaments. As a consequence, the LF solution systematically under-predicts the total interfacial area.

Finally, as the number of checkpoints $N_{c}$ increases, the MF-OT-ROM prediction progressively converges to the reference HF solution, confirming the effectiveness of the OT-based correction in recovering fine-scale interfacial dynamics. Interestingly, the MF-OT-ROM curve approaches the HF reference from above. This behavior may be attributed to the regularization employed within the OT solver. Indeed, as observed in the middle panel of Figure~\ref{fig:RT_fullfield}, a limited number of checkpoints produces a solution that appears less smooth. This residual lack of regularity can artificially increase the small-scale oscillations at the interface, ultimately leading to a slight over-prediction of the total interfacial area.
\subsubsection{PMF-OT-ROM}
Finally, we test the parametric multi-fidelity approach for the Rayleigh-Taylor instability problem, for which we have observed significant variations in the dynamics for different parameters in Figure~\ref{fig:RT_overview}.

\begin{figure}[t!]
\centering
\includegraphics[width=0.75\textwidth]{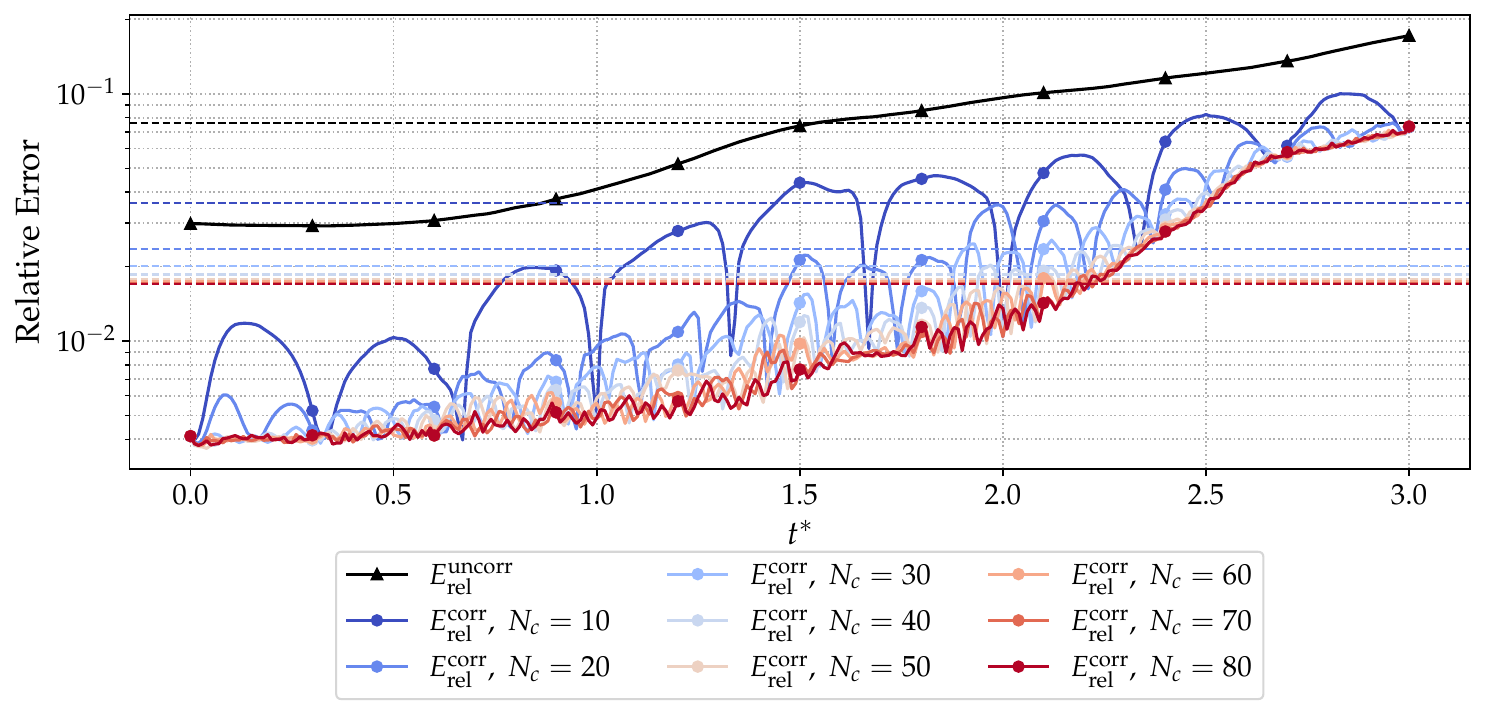}
\caption{Relative error between LF and HF for different number of checkpoints. Horizontal dashed lines represent the mean error for each temporal trajectory.}
\label{fig:RT_Nc}
\end{figure}

\begin{figure}[t!]
\centering
\includegraphics[width=0.75\textwidth]{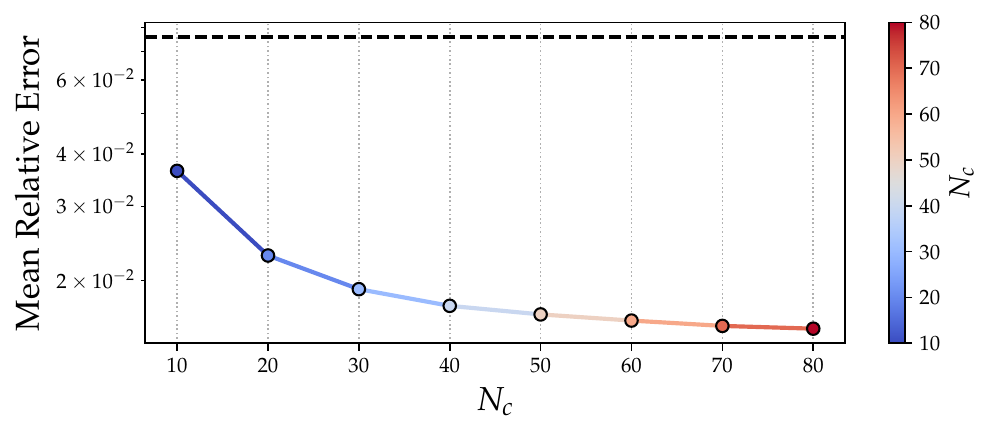}
\caption{Mean relative errors against number of checkpoints.}
\label{fig:RT_mean_par}
\end{figure}

In Figure~\ref{fig:RT_Nc} we compare how the relative error behaves for different number of checkpoints. Normally, we would expect that for larger numbers of checkpoints the error would decrease. This is generally the case but in a different manner depending on the problem's dynamics and its ``velocity''. At early stages, the improvement given by considering more checkpoints is significant. In fact, the initial dynamics is very similar across the different parameters as the initial condition for the phase-field is exactly the same across the whole dataset. However, as the dynamics of the problem unfolds, the type of breakups obtained at later stages is significantly different. At those stages, considering more checkpoints in time does not significantly improve the prediction of the OT solver. This clearly indicates that, as already previously observed, the benefit of using more checkpoints to reconstruct the fast-dynamical features of the problem is limited by the naive choice of the linear interpolator for the virtual time $\alpha$.
The same concept in terms of mean error can be seen in Figure~\ref{fig:RT_mean_par}, where we can still observe that the error decrease by increasing the number of checkpoints used in the OT strategy, but at a slower rate.

\begin{figure}[t!]
\centering
\includegraphics[width=0.95\textwidth]{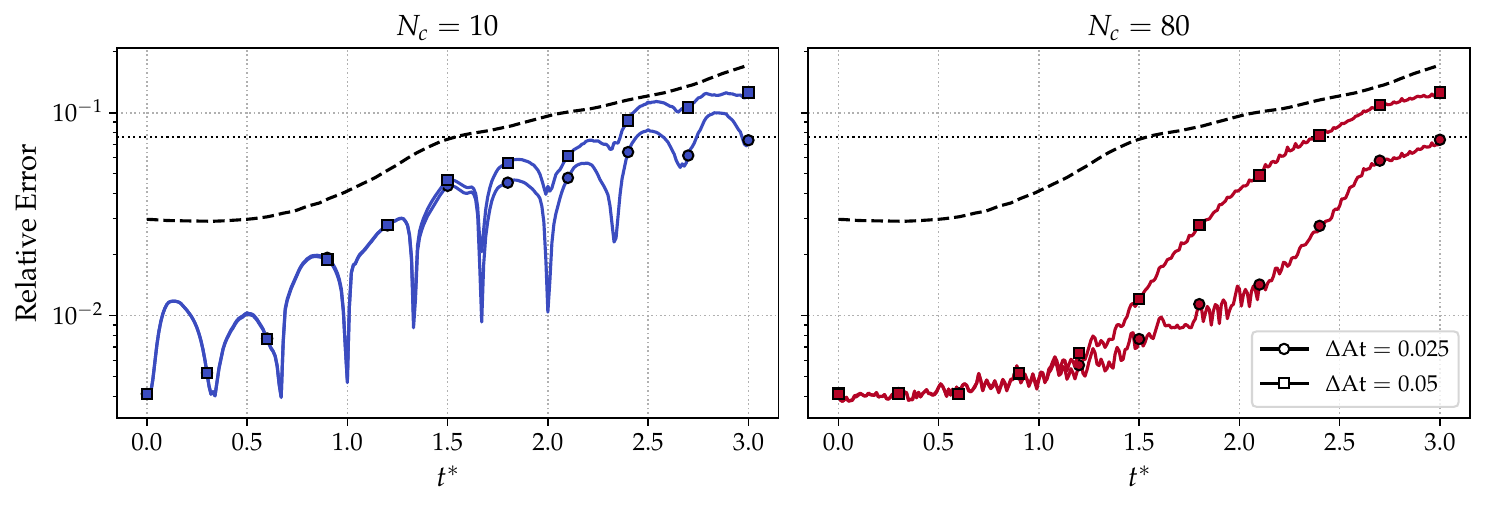}
\caption{Relative error between low- and high-fidelity using the full dataset (circles) and half of it (squares), for $N_{c}=10$ and $N_{c}=80$, left and right respectively. Dashed black line, uncorrected LF solution; dotted line, mean relative error for uncorrected LF solution.}
\label{fig:RT_Nc10_80}
\end{figure}

Within this framework we also modified the total number of snapshots at disposal so that we can observe the same trends in terms of parametric refinements. 
In Figure~\ref{fig:RT_Nc10_80} we compare how the relative error behaves for different number of checkpoints and different number of parametric trajectories considered. Normally, we would expect that given more parametric trajectories the error would decrease. In this scenario, at late stages of the problem's dynamics, considering more data improves the prediction of the overall procedure, even better when a large number of temporal checkpoints is considered (left panel of Figure~\ref{fig:RT_Nc10_80}).

\begin{figure}[t!]
\centering
\includegraphics[width=0.75\textwidth]{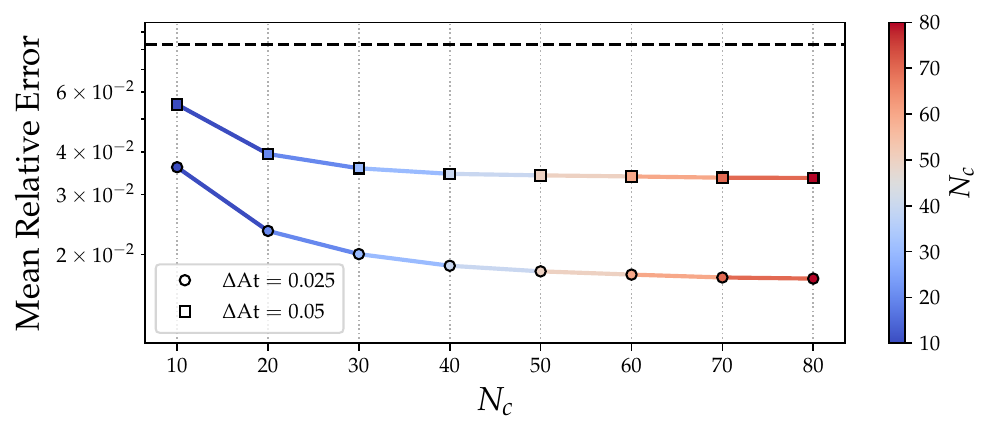}
\caption{Mean relative errors against number of checkpoints using the full dataset (circles) and half of it (squares).}
\label{fig:RT_mean_par2}
\end{figure}

The same concept in terms of mean error is shown in Figure~\ref{fig:RT_mean_par2} where we can observe that indeed the error decreases both by enlarging the initial dataset or alternatively by increasing the number of checkpoints used in the OT solver. In general, it is then interesting to notice that a certain trade-off between temporal and parametric refinements is necessary in order to achieve an optimal result.

\begin{figure}[t!]
\centering
\includegraphics[width=0.95\textwidth]{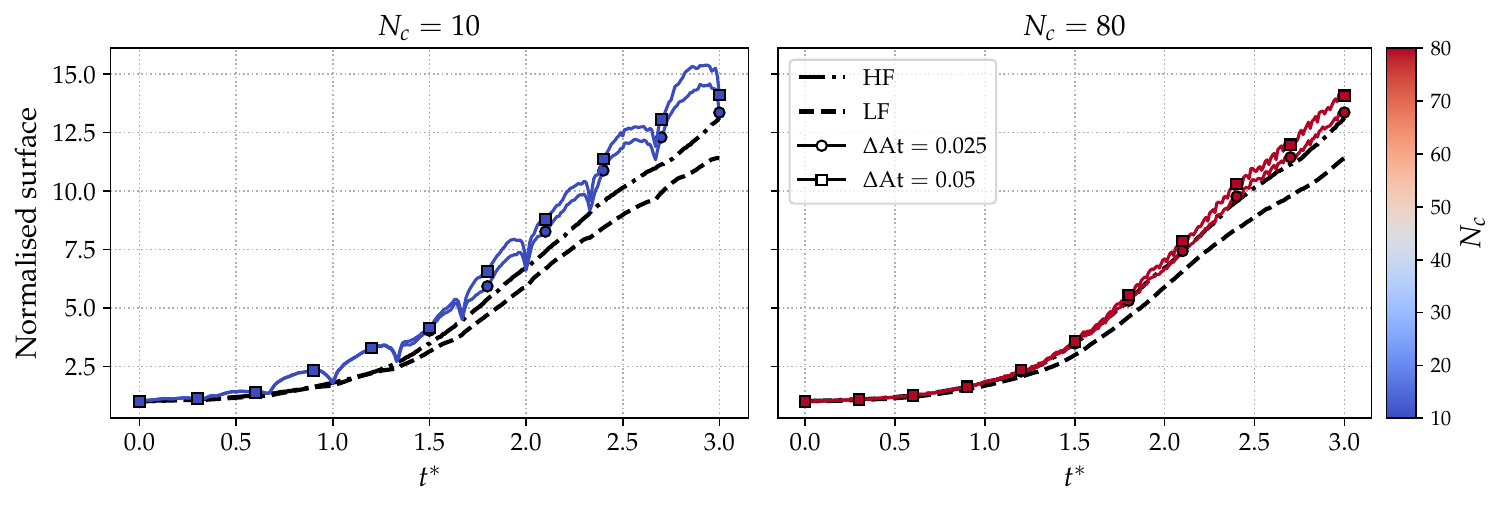}
\caption{Normalized interface area predicted using the full dataset (circles) and half of it (squares) for $N_{c}=10$ and $N_{c}=80$, left and right respectively.}
\label{fig:RT_surf_Nc10_80}
\end{figure}

As a more informative metric for this test case, we examine the total interfacial area as a function of time. This quantity is reported in Figure~\ref{fig:RT_surf_Nc10_80}, where we compare the predictions obtained with different numbers of checkpoints and different levels of parametric refinement.

In the early stages of the evolution, even a relatively small number of checkpoints ($N_c = 10$) is sufficient to accurately capture the interfacial dynamics. However, as time progresses and increasingly complex breakup events occur, discrepancies between the OT prediction and the HF solution become more pronounced. This deviation reflects the growing complexity of the interfacial topology, which is more challenging to reconstruct with a limited number of trajectories.

The agreement can be improved by increasing the number of trajectories in the parametric space, which enhances the resolution of the transported structures. Nevertheless, despite this refinement, we are not able to achieve the same level of agreement observed with the MF-OT-ROM approach applied exclusively in time (compare with Figure~\ref{fig:RT_surf}).

\begin{figure}[t!]
\centering
\includegraphics[width=0.65\textwidth]{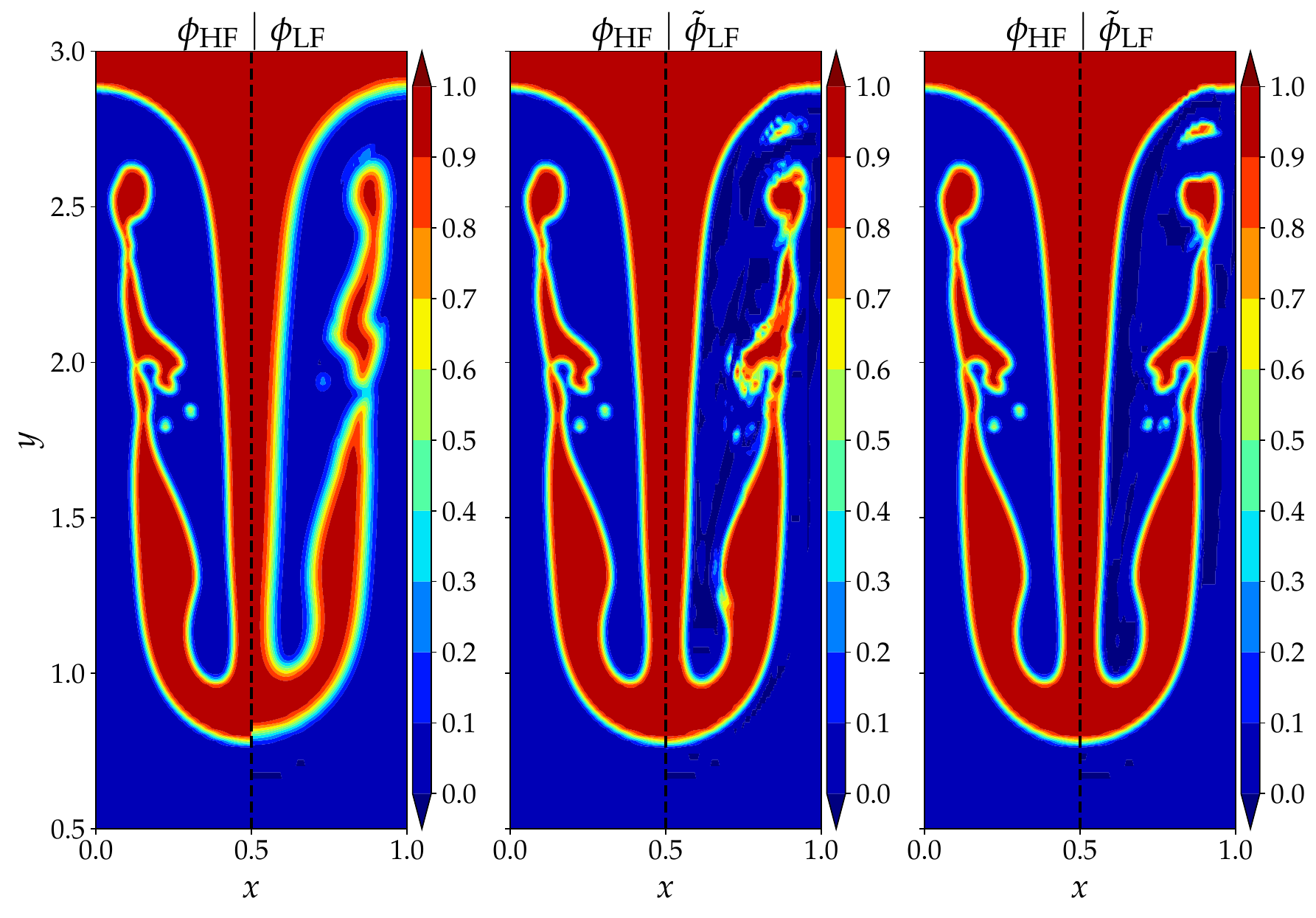}
\caption{Phase-field isocontours for Rayleigh-Taylor instability at $t^{*}=2.95$ and $\mathrm{At}=0.35$. In each subfigure, the left field is given by the HF solution, while right field shows respectively: LF solution, corrected LF solution with $N_{c}=10$ and corrected LF solution with $N_{c}=80$.}
\label{fig:RT_fullfield_par}
\end{figure}

In Figure~\ref{fig:RT_fullfield_par} we compare the fields in the parametric setting as well. We can notice that for $N_{c}=10$ the differences are quite noticeable but if we increase the number of checkpoints to $N_{c}=80$ we obtain a quite good agreement between the HF solution and the OT-corrected LF one. Similarly to the predicted total interfacial area, we are obtaining a good prediction in the parametric setting too although not as accurate as the single parametric multi-fidelity approach presented earlier.

In the previous analysis, we tested OT-based approaches on a single test parameter to validate the proposed methodologies. Thus, we conclude the investigation by examining the performance of the most general approach, namely the PMF-OT-ROM method, in the most challenging scenario considered in this work, i.e.\ the Rayleigh-Taylor instability, across the entire range of parametric configurations.

\begin{figure}[t!]
\centering
\includegraphics[width=0.95\textwidth]{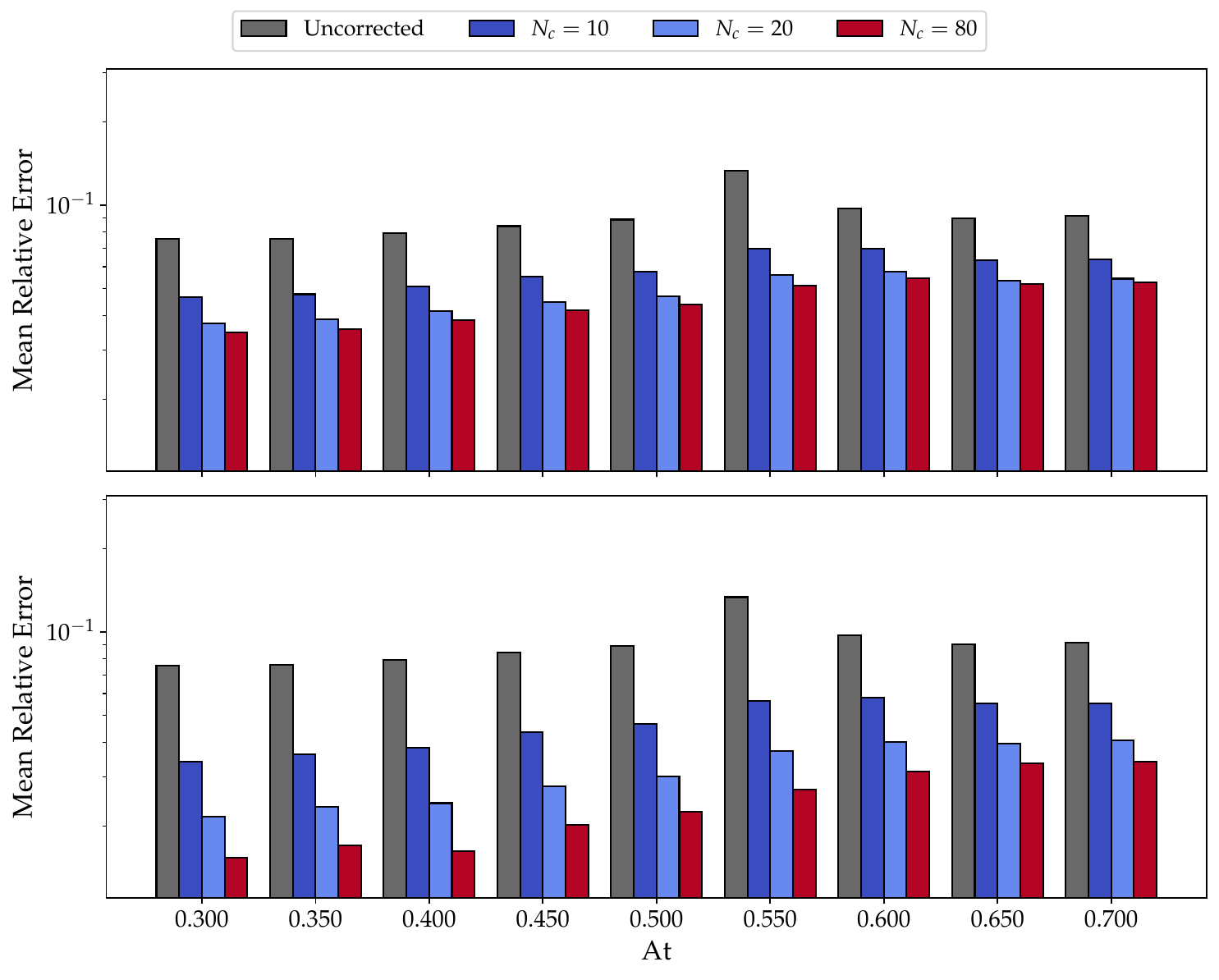}
\caption{Mean relative error analysis across the full parametric range with $\Delta \mathrm{At}=0.05$ (top) and $\Delta \mathrm{At}=0.025$ (bottom).}
\label{fig:RT_err_glob}
\end{figure}

Figure~\ref{fig:RT_err_glob} reports the error analysis performed over the entire parametric range in such setting, from which a few relevant trends can be observed. In general, the error tends to increase as the Atwood number grows, up to approximately $\mathrm{At}=0.5$. This behavior is expected, since for relatively small Atwood numbers the overall dynamics of the system remains more regular and less affected by strong interfacial breakups. Conversely, larger Atwood numbers lead to increasingly complex interfacial dynamics, characterized by stronger nonlinear effects and fragmentation phenomena (see Figure~\ref{fig:RT_overview}).

Concerning the PMF-OT-ROM approach, we observe that, even with a relatively small number of checkpoints ($N_c=10$), the corrected solution consistently improves the baseline uncorrected LF one for all the parametric trajectories. Furthermore, in agreement with previous observations for the reference parametric configuration, increasing the number of checkpoints significantly improves the accuracy of the corrected solution. This improvement appears to saturate relatively quickly when considering the coarser parametric spacing ($\Delta \mathrm{At}=0.05$). In contrast, for the finer parametric resolution ($\Delta \mathrm{At}=0.025$), a noticeable improvement in accuracy is observed when increasing the number of checkpoints from $N_c=20$ to $N_c=80$. This gain in accuracy is particularly evident for lower Atwood numbers, where the flow dynamics remain comparatively less complex.

\section{Conclusions and Perspectives}
\label{sec:conclusions}
In this paper, we have introduced two extensions to the Optimal Transport-based ROM framework. The multi-fidelity residual interpolation provides a method for correcting low-fidelity models by interpolating the HF-LF discrepancy field in time using displacement interpolation (MF-OT-ROM). The parametric displacement interpolation extends the framework to handle additional physical parameters beyond time, enabling prediction at unseen parameter values through a combination of parametric and temporal interpolation (PMF-OT-ROM). For multi-fidelity parametric problems, the HF solution is first interpolated across neighboring trajectories to produce synthetic data. These data are then used to compute residuals, which are subsequently employed to correct the LF solutions for previously unseen parametric configurations.

The mathematical formulation of these frameworks is general and can be applied to a wide range of problems governed by PDEs. By framing the discussion in the context of challenging diffuse-interface methods, we highlight the potential of these methods to handle complex, nonlinear dynamics involving moving fronts and interfaces. We have successfully applied the proposed extensions to two standard test cases for two-phase flow modeling, namely the Rider-Kothe vortex and the Rayleigh-Taylor instability. Both problems are characterized by unsteady dynamics with sharp structures that evolve and propagate in space. The Rayleigh-Taylor instability, in particular, involves complex interfacial dynamics with frequent breakup events, which become even more pronounced at large Atwood numbers. In both test cases, the proposed extensions significantly improve the baseline predictions, highlighting the potential of OT-based strategies for surrogate modeling in two-phase flows.

Several directions for future work can be identified.
First, the selection of training parameters and time checkpoints could be made adaptive, using error estimators to place computational effort where it is most needed.
Second, the linear interpolation of the OT map in parameter space could be extended to barycentric interpolation for multiple parameters.
Third, machine learning methods such as Gaussian Process Regression could be used to learn the mapping from physical parameters to the OT plans or interpolation weights.
Finally, OT methodologies could be further leveraged beyond surrogate modeling, for instance in the modeling of turbulent two-phase flows, such as predicting breakup dynamics.

\section*{Acknowledgments}
This work was partially supported by the European Union - NextGenerationEU, in the framework of the iNEST - Interconnected Nord-Est Innovation Ecosystem (iNEST ECS00000043 - CUP G93C22000610007).
The authors also acknowledge the support of the INdAM-GNCS.

\bibliographystyle{plainnat}
\bibliography{main}

@ARTICLE{agueh2011barycenters,
  author = {Agueh, Martial and Carlier, Guillaume},
  journal = {SIAM Journal on Mathematical Analysis},
  title = {Barycenters in the {W}asserstein Space},
  year = {2011},
  number = {2},
  pages = {904--924},
  volume = {43},
  doi = {10.1137/100805741},
}

@ARTICLE{Benner2015,
  author = {Benner, Peter and Gugercin, Serkan and Willcox, Karen},
  journal = {SIAM Review},
  title = {A Survey of Projection-Based Model Reduction Methods for Parametric Dynamical Systems},
  year = {2015},
  number = {4},
  pages = {483--531},
  volume = {57},
  doi = {10.1137/130932715},
}

@ARTICLE{Blickhan23,
  author = {Blickhan, Tobias},
  journal = {SIAM Journal on Scientific Computing},
  title = {A Registration Method for Reduced Basis Problems Using Linear Optimal Transport},
  year = {2024},
  number = {5},
  pages = {A3177--A3204},
  volume = {46},
  doi = {10.1137/23m1570715},
}

@ARTICLE{Chatterjee2000,
  author = {Chatterjee, Anindya},
  journal = {Current Science},
  title = {An Introduction to the Proper Orthogonal Decomposition},
  year = {2000},
  number = {7},
  pages = {808--817},
  volume = {78},
}

@ARTICLE{Cohen15,
  author = {Cohen, Albert and DeVore, Ronald},
  journal = {IMA Journal of Numerical Analysis},
  title = {Kolmogorov Widths Under Holomorphic Mappings},
  year = {2015},
  pages = {dru066},
  doi = {10.1093/imanum/dru066},
}

@INPROCEEDINGS{cuturi13,
  author = {Cuturi, Marco},
  booktitle = {Advances in Neural Information Processing Systems},
  title = {Sinkhorn Distances: Lightspeed Computation of Optimal Transport},
  year = {2013},
  editor = {Burges, C.J. and Bottou, L. and Welling, M. and Ghahramani, Z. and Weinberger, K.Q.},
  publisher = {Curran Associates, Inc.},
  volume = {26},
}

@ARTICLE{fresca2021,
  author = {Fresca, Stefania and Manzoni, Andrea and Ded{\`e}, Luca},
  journal = {Journal of Scientific Computing},
  title = {A Comprehensive Deep Learning-Based Approach to Reduced Order Modeling of Nonlinear Time-Dependent Parametrized {PDE}s},
  year = {2021},
  number = {2},
  pages = {1--42},
  volume = {87},
  doi={10.1007/s10915-021-01462-7}
}

@BOOK{Hesthaven2016,
  author = {Hesthaven, Jan S. and Rozza, Gianluigi and Stamm, Benjamin},
  publisher = {Springer International Publishing},
  title = {Certified Reduced Basis Methods for Parametrized Partial Differential Equations},
  year = {2016},
  series = {SpringerBriefs in Mathematics},
  doi = {10.1007/978-3-319-22470-1},
}

@ARTICLE{Khamlich2025JCPDI,
  author = {Khamlich, Moaad and Pichi, Federico and Girfoglio, Michele and Quaini, Annalisa and Rozza, Gianluigi},
  journal = {Journal of Computational Physics},
  title = {Optimal Transport-Based Displacement Interpolation with Data Augmentation for Reduced Order Modeling of Nonlinear Dynamical Systems},
  year = {2025},
  pages = {113938},
  volume = {531},
  doi = {10.1016/j.jcp.2025.113938},
}

@ARTICLE{LeeModelReductionDynamical2020,
  author = {Lee, Kookjin and Carlberg, Kevin T.},
  journal = {Journal of Computational Physics},
  title = {Model Reduction of Dynamical Systems on Nonlinear Manifolds Using Deep Convolutional Autoencoders},
  year = {2020},
  pages = {108973},
  volume = {404},
  doi = {10.1016/j.jcp.2019.108973},
}

@ARTICLE{mccan-interpolation,
  author = {McCann, Robert J.},
  journal = {Advances in Mathematics},
  title = {A Convexity Principle for Interacting Gases},
  year = {1997},
  number = {1},
  pages = {153--179},
  volume = {128},
  doi = {10.1006/aima.1997.1634},
}

@ARTICLE{Peherstorfer2018,
  author = {Peherstorfer, Benjamin and Willcox, Karen and Gunzburger, Max},
  journal = {SIAM Review},
  title = {Survey of Multifidelity Methods in Uncertainty Propagation, Inference, and Optimization},
  year = {2018},
  number = {3},
  pages = {550--591},
  volume = {60},
  doi = {10.1137/16M1082469},
}

@ARTICLE{Peyre2019,
  author = {Peyr{\'e}, Gabriel and Cuturi, Marco},
  journal = {Foundations and Trends in Machine Learning},
  title = {Computational Optimal Transport: With Applications to Data Science},
  year = {2019},
  number = {5--6},
  pages = {355--607},
  volume = {11},
  doi = {10.1561/2200000073},
}

@BOOK{Quarteroni2016,
  author = {Quarteroni, Alfio and Manzoni, Andrea and Negri, Federico},
  publisher = {Springer International Publishing},
  title = {Reduced Basis Methods for Partial Differential Equations},
  year = {2016},
  series = {UNITEXT},
  doi = {10.1007/978-3-319-15431-2},
}

@ARTICLE{Schmitzer2019,
  author = {Schmitzer, Bernhard},
  journal = {SIAM Journal on Scientific Computing},
  title = {Stabilized Sparse Scaling Algorithms for Entropy Regularized Transport Problems},
  year = {2019},
  number = {3},
  pages = {A1443--A1481},
  volume = {41},
  doi = {10.1137/16m1106018},
}

@ARTICLE{Taddei2020,
  author = {Taddei, Tommaso},
  journal = {SIAM Journal on Scientific Computing},
  title = {A Registration Method for Model Order Reduction: Data Compression and Geometry Reduction},
  year = {2020},
  number = {2},
  pages = {A997--A1027},
  volume = {42},
  doi = {10.1137/19M1271270},
}

@article{allaire2002five,
  title={A five-equation model for the simulation of interfaces between compressible fluids},
  author={Allaire, Gr{\'e}goire and Clerc, S{\'e}bastien and Kokh, Samuel},
  journal={Journal of Computational Physics},
  volume={181},
  pages={577--616},
  year={2002},
  doi={10.1006/jcph.2002.7143},
  publisher={Elsevier}
}

@article{al2021high,
  title={A high order flux reconstruction interface capturing method with a phase field preconditioning procedure},
  author={Al-Salami, Jabir and Kamra, Mohamed M and Hu, Changhong},
  journal={Journal of Computational Physics},
  volume={438},
  pages={110376},
  year={2021},
  doi={10.1016/j.jcp.2021.110376},
  publisher={Elsevier}
}

@article{tonicello2025extension,
  title={Extension of a spectral difference method for the diffused-interface five-equation model},
  author={Tonicello, Niccol{\`o} and Lodato, Guido and Ihme, Matthias},
  journal={Computers \& Fluids},
  pages={106880},
  year={2025},
  doi={10.1016/j.compfluid.2025.106880},
  publisher={Elsevier}
}

@article{rider1998reconstructing,
  title={Reconstructing volume tracking},
  author={Rider, William J and Kothe, Douglas B},
  journal={Journal of Computational Physics},
  volume={141},
  pages={112--152},
  year={1998},
  publisher={Elsevier}
}

@article{chiu2011conservative,
  title={A conservative phase field method for solving incompressible two-phase flows},
  author={Chiu, Pao-Hsiung and Lin, Yan-Ting},
  journal={Journal of Computational Physics},
  volume={230},
  pages={185--204},
  year={2011},
  doi={10.1016/j.jcp.2010.09.021},
  publisher={Elsevier}
}

@book{QuarteroniNumericalApproximationPartial1994,
  title = {Numerical {{Approximation}} of {{Partial Differential Equations}}},
  author = {Quarteroni, Alfio and Valli, Alberto},
  year = 1994,
  publisher = {Springer-Verlag},
  googlebooks = {MB23keacoWMC},
}

@book{ManzoniOptimalControlPartial2021,
  title = {Optimal {{Control}} of {{Partial Differential Equations}}: {{Analysis}}, {{Approximation}}, and {{Applications}}},
  shorttitle = {Optimal {{Control}} of {{Partial Differential Equations}}},
  author = {Manzoni, Andrea and Quarteroni, Alfio and Salsa, Sandro},
  year = 2021,
  series = {Applied {{Mathematical Sciences}}},
  volume = {207},
  publisher = {Springer International Publishing},
  address = {Cham},
  doi = {10.1007/978-3-030-77226-0},
}

@article{PichiDrivingBifurcatingParametrized2022,
  title = {Driving Bifurcating Parametrized Nonlinear {{PDEs}} by Optimal Control Strategies: Application to {{Navier}}--{{Stokes}} Equations with Model Order Reduction},
  shorttitle = {Driving Bifurcating Parametrized Nonlinear {{PDEs}} by Optimal Control Strategies},
  author = {Pichi, Federico and Strazzullo, Maria and Ballarin, Francesco and Rozza, Gianluigi},
  year = 2022,
  journal = {ESAIM: Mathematical Modelling and Numerical Analysis},
  volume = {56},
  number = {4},
  pages = {1361--1400},
  publisher = {EDP Sciences},
  doi = {10.1051/m2an/2022044}
}

@article{ChenReducedBasisMethods2017,
  title = {Reduced {{Basis Methods}} for {{Uncertainty Quantification}}},
  author = {Chen, Peng and Quarteroni, Alfio and Rozza, Gianluigi},
  year = 2017,
  journal = {SIAM/ASA Journal on Uncertainty Quantification},
  volume = {5},
  number = {1},
  pages = {813--869},
  publisher = {{Society for Industrial and Applied Mathematics}},
  doi = {10.1137/151004550}
}

@incollection{MulaInverseProblemsDeterministic2023,
  title = {Inverse {{Problems}}: {{A Deterministic Approach Using Physics-Based Reduced Models}}},
  shorttitle = {Inverse {{Problems}}},
  booktitle = {Model {{Order Reduction}} and {{Applications}}: {{Cetraro}}, {{Italy}} 2021},
  author = {Mula, Olga},
  editor = {Hinze, Michael and Kutz, J. Nathan and Mula, Olga and Urban, Karsten and Falcone, Maurizio and Rozza, Gianluigi},
  year = 2023,
  pages = {73--124},
  publisher = {Springer Nature Switzerland},
  address = {Cham},
  doi = {10.1007/978-3-031-29563-8_2}
}

@article{GiacominiSurrogateModelTopology2026,
  title = {A Surrogate Model for Topology Optimisation of Elastic Structures via Parametric Autoencoders},
  author = {Giacomini, Matteo and Huerta, Antonio},
  year = 2026,
  journal = {Computer Methods in Applied Mechanics and Engineering},
  volume = {448},
  pages = {118503},
  doi = {10.1016/j.cma.2025.118503}
}

@article{BabuskaStochasticCollocationMethod2007,
  title = {A Stochastic Collocation Method for Elliptic Partial Differential Equations with Random Input Data},
  author = {Babuska, Ivo and Nobile, Fabio and Tempone, Ra{\'u}l},
  year = 2007,
  journal = {SIAM Journal on Numerical Analysis},
  volume = {45},
  number = {3},
  eprint = {https://doi.org/10.1137/050645142},
  pages = {1005--1034},
  doi = {10.1137/050645142}
}

@article{BoulleControlBifurcationStructures2022,
  title = {Control of {{Bifurcation Structures}} Using {{Shape Optimization}}},
  author = {Boull{\'e}, Nicolas and Farrell, Patrick E. and Paganini, Alberto},
  year = 2022,
  journal = {SIAM Journal on Scientific Computing},
  volume = {44},
  number = {1},
  pages = {A57-A76},
  publisher = {{Society for Industrial and Applied Mathematics}},
  doi = {10.1137/21M1418708}
}

@book{isakov2017inverse,
  title = {Inverse Problems for Partial Differential Equations},
  author = {Isakov, V.},
  year = 2017,
  series = {Applied Mathematical Sciences},
  publisher = {Springer International Publishing}
}

@article{SirovichTurbulenceDynamicsCoherent1987a,
  title = {Turbulence and the {{Dynamics}} of {{Coherent Structures Part I}}: {{Coherent Structures}}},
  shorttitle = {Turbulence and the {{Dynamics}} of {{Coherent Structures Part I}}},
  author = {Sirovich, Lawrence},
  year = 1987,
  journal = {Quarterly of Applied Mathematics},
  volume = {45},
  number = {3},
  eprint = {43637457},
  eprinttype = {jstor},
  pages = {561--571},
  publisher = {Brown University}
}

@article{OhlbergerReducedBasisMethods2016,
  title = {Reduced {{Basis Methods}}: {{Success}}, {{Limitations}} and {{Future Challenges}}},
  shorttitle = {Reduced {{Basis Methods}}},
  author = {Ohlberger, Mario and Rave, Stephan},
  year = 2016,
  eprint = {1511.02021},
  primaryclass = {math},
  publisher = {arXiv},
  doi = {10.48550/arXiv.1511.02021},
  archiveprefix = {arXiv}
}

@article{RomorROMViscousIncompressible2025,
  title = {{{ROM}} for {{Viscous}}, {{Incompressible Flow}} in {{Polygons}} -- Exponential n-Width Bounds and Convergence Rate},
  author = {Romor, Francesco and Pichi, Federico and Stabile, Giovanni and Rozza, Gianluigi and Schwab, Christoph},
  year = 2025,
  eprint = {2512.22567},
  primaryclass = {math},
  publisher = {arXiv},
  doi = {10.48550/arXiv.2512.22567},
  archiveprefix = {arXiv}
}

@article{Lassila2013,
author = {Lassila, Toni and Manzoni, Andrea and Quarteroni, Alfio and Rozza, Gianluigi},
journal = {Bollettino dell'Unione Matematica Italiana},
language = {eng},
month = {2},
number = {1},
pages = {113-135},
publisher = {Unione Matematica Italiana},
title = {Generalized Reduced Basis Methods and n-width Estimates for the Approximation of the Solution Manifold of Parametric PDEs},
url = {http://eudml.org/doc/294013},
volume = {6},
year = {2013},
}

@article{GlasReducedBasisMethod2020,
  title = {A Reduced Basis Method for the Wave Equation},
  author = {Glas, Silke and Patera, Anthony T. and Urban, Karsten},
  year = 2020,
  journal = {International Journal of Computational Fluid Dynamics},
  volume = {34},
  number = {2},
  pages = {139--146},
  publisher = {IAHR Website},
  doi = {10.1080/10618562.2019.1686486},
  urldate = {2025-07-28}
}

@article{ArbesKolmogorovNwidthLinear2025,
  title = {The {{Kolmogorov N-width}} for Linear Transport: Exact Representation and the Influence of the Data},
  shorttitle = {The {{Kolmogorov N-width}} for Linear Transport},
  author = {Arbes, Florian and Greif, Constantin and Urban, Karsten},
  year = 2025,
  journal = {Advances in Computational Mathematics},
  volume = {51},
  number = {2},
  pages = {13},
  doi = {10.1007/s10444-025-10224-0},
  urldate = {2025-04-07}
}

@article{BuffaPrioriConvergenceGreedy2012,
  title = {A Priori Convergence of the {{Greedy}} Algorithm for the Parametrized Reduced Basis Method},
  author = {Buffa, Annalisa and Maday, Yvon and Patera, Anthony T. and Prud'homme, Christophe and Turinici, Gabriel},
  year = 2012,
  journal = {ESAIM: Mathematical Modelling and Numerical Analysis},
  volume = {46},
  number = {3},
  pages = {595--603},
  publisher = {EDP Sciences},
  doi = {10.1051/m2an/2011056}
}

@article{HesthavenNonlinearModelReduction2026,
  title = {Nonlinear Model Reduction for Transport-Dominated Problems},
  author = {Hesthaven, Jan S. and Peherstorfer, Benjamin and Unger, Benjamin},
  year = 2026,
  eprint = {2602.01397},
  primaryclass = {math},
  publisher = {arXiv},
  doi = {10.48550/arXiv.2602.01397},
  archiveprefix = {arXiv}
}

@article{PeherstorferBreakingKolmogorovBarrier2022,
  title = {Breaking the {{Kolmogorov Barrier}} with {{Nonlinear Model Reduction}}},
  author = {Peherstorfer, Benjamin},
  year = 2022,
  journal = {Notices of the American Mathematical Society},
  volume = {69},
  number = {05},
  doi = {10.1090/noti2475}
}

@article{NoninoCalibrationBasedALEModel2024a,
  title = {Calibration-{{Based ALE Model Order Reduction}} for {{Hyperbolic Problems}} with {{Self-Similar Travelling Discontinuities}}},
  author = {Nonino, Monica and Torlo, Davide},
  year = 2024,
  journal = {Journal of Scientific Computing},
  volume = {101},
  number = {3},
  pages = {60},
  doi = {10.1007/s10915-024-02694-z}
}

@article{PichiGraphConvolutionalAutoencoder2024,
  title = {A Graph Convolutional Autoencoder Approach to Model Order Reduction for Parametrized {{PDEs}}},
  author = {Pichi, Federico and Moya, Beatriz and Hesthaven, Jan S.},
  year = 2024,
  journal = {Journal of Computational Physics},
  volume = {501},
  pages = {112762},
  doi = {10.1016/j.jcp.2024.112762}
}

@article{BattistiWassersteinModelReduction2023,
  title = {Wasserstein Model Reduction Approach for Parametrized Flow Problems in Porous Media},
  author = {Battisti, Beatrice and Blickhan, Tobias and Enchery, Guillaume and Ehrlacher, Virginie and Lombardi, Damiano and Mula, Olga},
  year = 2023,
  journal = {ESAIM: Proceedings and Surveys},
  volume = {73},
  pages = {28--47},
  publisher = {EDP Sciences},
  doi = {10.1051/proc/202373028}
}

@article{EhrlacherNonlinearModelReduction2020,
  title = {Nonlinear Model Reduction on Metric Spaces. {{Application}} to One-Dimensional Conservative {{PDEs}} in {{Wasserstein}} Spaces},
  author = {Ehrlacher, Virginie and Lombardi, Damiano and Mula, Olga and Vialard, Fran{\c c}ois-Xavier},
  year = 2020,
  journal = {ESAIM: Mathematical Modelling and Numerical Analysis},
  volume = {54},
  number = {6},
  pages = {2159--2197},
  publisher = {EDP Sciences},
  doi = {10.1051/m2an/2020013}
}

@article{DoSparseWassersteinBarycenters2025,
  title = {Sparse {{Wasserstein Barycenters}} and {{Application}} to {{Reduced Order Modeling}}},
  author = {Do, Minh-Hieu and Feydy, Jean and Mula, Olga},
  year = 2025,
  journal = {Journal of Scientific Computing},
  volume = {102},
  number = {3},
  pages = {64},
  doi = {10.1007/s10915-024-02766-0}
}

@article{CucchiaraModelOrderReduction2024,
  title = {Model Order Reduction by Convex Displacement Interpolation},
  author = {Cucchiara, Simona and Iollo, Angelo and Taddei, Tommaso and Telib, Haysam},
  year = 2024,
  journal = {Journal of Computational Physics},
  volume = {514},
  pages = {113230},
  doi = {10.1016/j.jcp.2024.113230}
}

@article{IolloMappingCoherentStructures2022,
  title = {Mapping of Coherent Structures in Parameterized Flows by Learning Optimal Transportation with {{Gaussian}} Models},
  author = {Iollo, Angelo and Taddei, Tommaso},
  year = 2022,
  journal = {Journal of Computational Physics},
  volume = {471},
  pages = {111671},
  doi = {10.1016/j.jcp.2022.111671}
}

@book{SantambrogioOptimalTransportApplied2015,
  title = {Optimal {{Transport}} for {{Applied Mathematicians}}: {{Calculus}} of {{Variations}}, {{PDEs}}, and {{Modeling}}},
  shorttitle = {Optimal {{Transport}} for {{Applied Mathematicians}}},
  author = {Santambrogio, Filippo},
  year = 2015,
  series = {Progress in {{Nonlinear Differential Equations}} and {{Their Applications}}},
  volume = {87},
  publisher = {Springer International Publishing},
  address = {Cham},
  doi = {10.1007/978-3-319-20828-2},
  urldate = {2026-02-24},
  isbn = {978-3-319-20827-5 978-3-319-20828-2}
}

@book{villani2021topics,
  title = {Topics in Optimal Transportation},
  author = {Villani, C.},
  year = 2021,
  series = {Graduate Studies in Mathematics},
  publisher = {American Mathematical Society}
}

@article{ContiMultifidelityReducedorderSurrogate2023,
  title = {Multi-Fidelity Reduced-Order Surrogate Modelling},
  author = {Conti, Paolo and Guo, Mengwu and Manzoni, Andrea and Frangi, Attilio and Brunton, Steven L. and Nathan Kutz, J.},
  year = 2024,
  journal = {Proceedings of the Royal Society A: Mathematical, Physical and Engineering Sciences},
  volume = {480},
  number = {2283},
  pages = {20230655},
  doi = {10.1098/rspa.2023.0655}
}

@article{HowardMultifidelityDeepOperator2023,
  title = {Multifidelity Deep Operator Networks for Data-Driven and Physics-Informed Problems},
  author = {Howard, Amanda A. and Perego, Mauro and Karniadakis, George Em and Stinis, Panos},
  year = 2023,
  journal = {Journal of Computational Physics},
  volume = {493},
  pages = {112462},
  doi = {10.1016/j.jcp.2023.112462}
}

@article{Monge1781,
  author = {Monge, Gaspard},
  publisher = {Académie Royale des Sciences},
  year = {1781},
  journal = {Histoire de l'Académie Royale des Sciences},
  title = {{{Mémoire sur la théorie des déblais et des remblais}}}
}

@article{Kantorovich1942,
  author = {Kantorovich, L. V.},
  publisher = {Springer Science and Business Media LLC},
  year = {2006},
  doi = {10.1007/s10958-006-0049-2},
  issn = {1573-8795},
  journal = {Journal of Mathematical Sciences},
  number = {4},
  pages = {1381--1382},
  title = {On the Translocation of Masses},
  volume = {133}
}

@article{saurel2018diffuse,
  title={Diffuse-interface capturing methods for compressible two-phase flows},
  author={Saurel, Richard and Pantano, Carlos},
  journal={Annual Review of Fluid Mechanics},
  volume={50},
  pages={105--130},
  year={2018},
  doi={10.1146/annurev-fluid-122316-050109},
  publisher={Annual Reviews}
}

@article{roccon2025boiling,
  title={Boiling heat transfer by phase-field method: A. Roccon},
  author={Roccon, Alessio},
  journal={Acta Mechanica},
  volume={236},
  number={9},
  pages={5623--5638},
  year={2025},
  doi={10.1007/s00707-024-04122-7},
  publisher={Springer}
}

@article{weber2026consistent,
  title={A consistent and scalable framework suitable for boiling flows using the conservative diffuse interface method},
  author={Weber, Lorenz and Mukherjee, Aritra and Class, Andreas G and Brandt, Luca},
  journal={Journal of Computational Physics},
  pages={114680},
  year={2026},
  doi={10.1016/j.jcp.2026.114680},
  publisher={Elsevier}
}

@article{mirjalili2020conservative,
  title={A conservative diffuse interface method for two-phase flows with provable boundedness properties},
  author={Mirjalili, Shahab and Ivey, Christopher B and Mani, Ali},
  journal={Journal of Computational Physics},
  volume={401},
  pages={109006},
  year={2020},
  doi={10.1016/j.jcp.2019.109006},
  publisher={Elsevier}
}

@article{mirjalili2022computational,
  title={A computational model for interfacial heat and mass transfer in two-phase flows using a phase field method},
  author={Mirjalili, Shahab and Jain, Suhas S and Mani, Ali},
  journal={International Journal of Heat and Mass Transfer},
  volume={197},
  pages={123326},
  year={2022},
  doi={10.1016/j.ijheatmasstransfer.2022.123326},
  publisher={Elsevier}
}

@article{jain2022accurate,
  title={Accurate conservative phase-field method for simulation of two-phase flows},
  author={Jain, Suhas S},
  journal={Journal of Computational Physics},
  volume={469},
  pages={111529},
  year={2022},
  doi={10.1016/j.jcp.2022.111529},
  publisher={Elsevier}
}

@article{soligo2019breakage,
  title={Breakage, coalescence and size distribution of surfactant-laden droplets in turbulent flow},
  author={Soligo, Giovanni and Roccon, Alessio and Soldati, Alfredo},
  journal={Journal of Fluid Mechanics},
  volume={881},
  pages={244--282},
  year={2019},
  doi={10.1017/jfm.2019.772},
  publisher={Cambridge University Press}
}

@article{soligo2019mass,
  title={Mass-conservation-improved phase field methods for turbulent multiphase flow simulation: G. Soligo et al.},
  author={Soligo, Giovanni and Roccon, Alessio and Soldati, Alfredo},
  journal={Acta Mechanica},
  volume={230},
  number={2},
  pages={683--696},
  year={2019},
  doi={10.1007/s00707-018-2304-2},
  publisher={Springer}
}

@article{roccon2023phase,
  title={Phase-field modeling of complex interface dynamics in drop-laden turbulence},
  author={Roccon, Alessio and Zonta, Francesco and Soldati, Alfredo},
  journal={Physical Review Fluids},
  volume={8},
  number={9},
  pages={090501},
  year={2023},
  doi={10.1103/PhysRevFluids.8.090501},
  publisher={APS}
}

@article{huang2023consistent,
  title={A consistent and conservative Phase-Field method for compressible multiphase flows with shocks},
  author={Huang, Ziyang and Johnsen, Eric},
  journal={Journal of Computational Physics},
  volume={488},
  pages={112195},
  year={2023},
  doi={10.1016/j.jcp.2023.112195},
  publisher={Elsevier}
}

@article{collis2025thermodynamically,
  title={A thermodynamically consistent and robust four-equation model for multi-phase multi-component compressible flows using ENO-type schemes including interface regularization},
  author={Collis, Henry and Bezgin, Deniz A and Mirjalili, Shahab and Mani, Ali},
  journal={arXiv preprint arXiv:2504.14063},
  doi={10.48550/arXiv.2504.14063},
  year={2025}
}

@article{tonicello2024high,
  title={A high-order diffused-interface approach for two-phase compressible flow simulations using a discontinuous Galerkin framework},
  author={Tonicello, Niccol{\`o} and Ihme, Matthias},
  journal={Journal of Computational Physics},
  volume={508},
  pages={112983},
  year={2024},
  doi={10.1016/j.jcp.2024.112983},
  publisher={Elsevier}
}

@article{ibrahim2025conservative,
  title={A conservative phase-field lattice Boltzmann method for boiling heat transfer at high density ratios},
  author={Ibrahim, Mohammed and Rajamuni, Methma and Zhang, Chuangde and Chen, Li and Young, John and Tian, Fang-Bao},
  journal={Physics of Fluids},
  volume={37},
  number={5},
  year={2025},
  doi={10.1063/5.0261024},
  publisher={AIP Publishing}
}

@article{salimi2025low,
  title={A low Mach number diffuse-interface model for multicomponent two-phase flows with phase change},
  author={Salimi, Salar Zamani and Mukherjee, Aritra and Pelanti, Marica and Brandt, Luca},
  journal={Journal of Computational Physics},
  volume={523},
  pages={113683},
  year={2025},
  doi={10.1016/j.jcp.2024.113683},
  publisher={Elsevier}
}

@article{PichiArtificialNeuralNetwork2023,
  title = {An Artificial Neural Network Approach to Bifurcating Phenomena in Computational Fluid Dynamics},
  author = {Pichi, Federico and Ballarin, Francesco and Rozza, Gianluigi and Hesthaven, Jan S.},
  year = 2023,
  journal = {Computers \& Fluids},
  volume = {254},
  pages = {105813},
  doi = {10.1016/j.compfluid.2023.105813}
}

@article{GeelenLocalizedNonintrusiveReducedorder2022,
  title = {Localized Non-Intrusive Reduced-Order Modelling in the Operator Inference Framework},
  author = {Geelen, Rudy and Willcox, Karen},
  year = 2022,
  journal = {Philosophical Transactions of the Royal Society A: Mathematical, Physical and Engineering Sciences},
  volume = {380},
  number = {2229},
  pages = {20210206},
  publisher = {Royal Society},
  doi = {10.1098/rsta.2021.0206}
}

@article{JinAdaptiveHybridReduced2025,
  title = {Adaptive and Hybrid Reduced Order Models to Mitigate {{Kolmogorov}} Barrier in a Multiscale Kinetic Transport Equation},
  author = {Jin, Tianyu and Peng, Zhichao and Xiang, Yang},
  year = 2025,
  eprint = {2505.08214},
  primaryclass = {math},
  publisher = {arXiv},
  doi = {10.48550/arXiv.2505.08214},
  urldate = {2025-05-14},
  archiveprefix = {arXiv}
}

@article{AmsallemNonlinearModelOrder2012,
  title = {Nonlinear Model Order Reduction Based on Local Reduced-Order Bases},
  author = {Amsallem, David and Zahr, Matthew J. and Farhat, Charbel},
  year = 2012,
  journal = {International Journal for Numerical Methods in Engineering},
  volume = {92},
  number = {10},
  pages = {891--916},
  doi = {10.1002/nme.4371}
}

@article{HessLocalizedReducedorderModeling2019,
  title = {A Localized Reduced-Order Modeling Approach for {{PDEs}} with Bifurcating Solutions},
  author = {Hess, Martin and Alla, Alessandro and Quaini, Annalisa and Rozza, Gianluigi and Gunzburger, Max},
  year = 2019,
  journal = {Computer Methods in Applied Mechanics and Engineering},
  volume = {351},
  pages = {379--403},
  doi = {10.1016/j.cma.2019.03.050}
}

@article{PeherstorferModelReductionTransportDominated2020,
author = {Peherstorfer, Benjamin},
title = {Model Reduction for Transport-Dominated Problems via Online Adaptive Bases and Adaptive Sampling},
journal = {SIAM Journal on Scientific Computing},
volume = {42},
number = {5},
pages = {A2803-A2836},
year = {2020},
doi = {10.1137/19M1257275}
}

@article{RimModelReductionTransportDominated2023,
author = {Rim, Donsub and Peherstorfer, Benjamin and Mandli, Kyle T.},
title = {Manifold Approximations via Transported Subspaces: Model Reduction for Transport-Dominated Problems},
journal = {SIAM Journal on Scientific Computing},
volume = {45},
number = {1},
pages = {A170-A199},
year = {2023},
doi = {10.1137/20M1316998}
}

@techreport{harlow1971fluid,
  title={Fluid Dynamics. {A} {LASL} Monograph},
  author={Harlow, Francis H and Amsden, Anthony A},
  year={1971},
  institution={Los Alamos National Lab (LANL), Los Alamos, NM (United States)}
}

@article{jain2020conservative,
  title={A conservative diffuse-interface method for compressible two-phase flows},
  author={Jain, Suhas S and Mani, Ali and Moin, Parviz},
  journal={Journal of Computational Physics},
  volume={418},
  pages={109606},
  year={2020},
  doi={10.1016/j.jcp.2020.109606},
  publisher={Elsevier}
}

@article{pirozzoli2025efficient,
  title={Efficient implementation of the Allen--Cahn phase-field method for material interface tracking},
  author={Pirozzoli, Sergio and Di Giorgio, Simone and Rossi, Daniele},
  journal={Computers \& Fluids},
  pages={106768},
  year={2025},
  doi={10.1016/j.compfluid.2025.106768},
  publisher={Elsevier}
}

@article{white2025high,
  title={A high-order discontinuous Galerkin method for compressible interfacial flows with consistent and conservative Phase Fields},
  author={White, William J and Huang, Ziyang and Johnsen, Eric},
  journal={Journal of Computational Physics},
  volume={527},
  pages={113830},
  year={2025},
  doi={10.1016/j.jcp.2025.113830},
  publisher={Elsevier}
}

@article{huang2025bound,
  title={Bound preservation for the consistent and conservative phase-field method for compressible single-, two-, and N-phase flows},
  author={Huang, Ziyang and Johnsen, Eric},
  journal={Journal of Computational Physics},
  volume={526},
  pages={113783},
  year={2025},
  doi={10.1016/j.jcp.2025.113783},
  publisher={Elsevier}
}

@article{huang2022consistent,
  title={A consistent and conservative phase-field model for thermo-gas-liquid-solid flows including liquid-solid phase change},
  author={Huang, Ziyang and Lin, Guang and Ardekani, Arezoo M},
  journal={Journal of Computational Physics},
  volume={449},
  pages={110795},
  year={2022},
  doi={10.1016/j.jcp.2021.110795},
  publisher={Elsevier}
}

@article{haghani2021phase,
  title={Phase-change modeling based on a novel conservative phase-field method},
  author={Haghani-Hassan-Abadi, Reza and Fakhari, Abbas and Rahimian, Mohammad-Hassan},
  journal={Journal of Computational Physics},
  volume={432},
  pages={110111},
  year={2021},
  doi={10.1016/j.jcp.2021.110111},
  publisher={Elsevier}
}

@article{cutforth2025convolutional,
  title={Convolutional autoencoders for the reconstruction of three-dimensional interfacial multiphase flows},
  author={Cutforth, Murray and Mirjalili, Shahab},
  journal={AI Thermal Fluids},
  pages={100035},
  year={2026},
  doi={10.1016/j.aitf.2026.100035},
  publisher={Elsevier}
}

@inproceedings{wiewel2019latent,
  title={Latent space physics: Towards learning the temporal evolution of fluid flow},
  author={Wiewel, Steffen and Becher, Moritz and Thuerey, Nils},
  booktitle={Computer graphics forum},
  volume={38},
  number={2},
  pages={71--82},
  year={2019},
  organization={Wiley Online Library}
}

@article{kani2019reduced,
  title={Reduced-order modeling of subsurface multi-phase flow models using deep residual recurrent neural networks},
  author={Kani, J Nagoor and Elsheikh, Ahmed H},
  journal={Transport in Porous Media},
  volume={126},
  number={3},
  pages={713--741},
  year={2019},
  publisher={Springer}
}

@article{haas2020bubcnn,
  title={BubCNN: Bubble detection using Faster RCNN and shape regression network},
  author={Haas, Tim and Schubert, Christian and Eickhoff, Moritz and Pfeifer, Herbert},
  journal={Chemical Engineering Science},
  volume={216},
  pages={115467},
  year={2020},
  doi={10.1016/j.ces.2019.115467},
  publisher={Elsevier}
}

@article{KhamlichEfficientNumericalStrategies2026,
  title = {Efficient Numerical Strategies for Entropy-Regularized Semi-Discrete Optimal Transport},
  author = {Khamlich, Moaad and Romor, Francesco and Rozza, Gianluigi},
  year = 2026,
  journal = {Computer Methods in Applied Mechanics and Engineering},
  volume = {453},
  pages = {118821},
  doi = {10.1016/j.cma.2026.118821}
}

@article{LevyNumericalAlgorithmL22015,
  title = {A {{Numerical Algorithm}} for {{L2 Semi-Discrete Optimal Transport}} in {{3D}}},
  author = {L{\'e}vy, Bruno},
  year = 2015,
  journal = {ESAIM: Mathematical Modelling and Numerical Analysis},
  volume = {49},
  number = {6},
  pages = {1693--1715},
  publisher = {EDP Sciences},
  doi = {10.1051/m2an/2015055}
}

@article{BenamouComputationalFluidMechanics2000,
  title = {A Computational Fluid Mechanics Solution to the {{Monge-Kantorovich}} Mass Transfer Problem},
  author = {Benamou, Jean-David and Brenier, Yann},
  year = 2000,
  journal = {Numerische Mathematik},
  volume = {84},
  number = {3},
  pages = {375--393},
  doi = {10.1007/s002110050002}
}

@article{KastNonintrusiveMultifidelityMethod2020,
  title = {A Non-Intrusive Multifidelity Method for the Reduced Order Modeling of Nonlinear Problems},
  author = {Kast, Mariella and Guo, Mengwu and Hesthaven, Jan S.},
  year = 2020,
  journal = {Computer Methods in Applied Mechanics and Engineering},
  volume = {364},
  pages = {112947},
  doi = {10.1016/j.cma.2020.112947},
}

@article{MorrisonGFNGraphFeedforward2024a,
  title = {{{GFN}}: {{A}} Graph Feedforward Network for Resolution-Invariant Reduced Operator Learning in Multifidelity Applications},
  shorttitle = {{{GFN}}},
  author = {Morrison, Ois{\'i}n M. and Pichi, Federico and Hesthaven, Jan S.},
  year = 2024,
  journal = {Computer Methods in Applied Mechanics and Engineering},
  volume = {432},
  pages = {117458},
  doi = {10.1016/j.cma.2024.117458}
}

@article{baer1986two,
  title={A two-phase mixture theory for the deflagration-to-detonation transition (DDT) in reactive granular materials},
  author={Baer, Melvin R and Nunziato, Jace W},
  journal={International journal of multiphase flow},
  volume={12},
  number={6},
  pages={861--889},
  year={1986},
  doi={10.1016/0301-9322(86)90033-9},
  publisher={Elsevier}
}

@article{tryggvason1988numerical,
  title={Numerical simulations of the Rayleigh-Taylor instability},
  author={Tryggvason, Gr{\'e}tar},
  journal={Journal of Computational Physics},
  volume={75},
  number={2},
  pages={253--282},
  year={1988},
  doi={10.1016/0021-9991(88)90112-X},
  publisher={Elsevier}
}

\end{document}